\documentclass[12pt]{amsart}
\usepackage{amsmath}
\usepackage{rotating}
\usepackage{graphics}
\usepackage[all]{xy}
\xyoption{v2} \xyoption{knot}

\let\url=\undefined
\usepackage
[bookmarks,colorlinks,linkcolor=blue,citecolor=blue,urlcolor=blue] {hyperref}

\title{Classification of braids which give rise to interchange}
\author {Stefan Forcey}
\address{Department of Physics and Mathematics\\
Tennessee State University\\
Nashville, TN 37209 \\
USA\\
}
 \email{sforcey@tnstate.edu}
  \urladdr{http://faculty.tnstate.edu/sforcey/}

\author { Felita Humes}
\address{Department of Physics and Mathematics\\
Tennessee State University\\
Nashville, TN 37209 \\
USA\\
}
 \email{felitahume@hotmail.com}

\thanks{Thanks to {\Xy-pic} for the diagrams. }


\swapnumbers \theoremstyle{plain}
\newtheorem{theorem}{Theorem}[section]
\newtheorem{lemma}[theorem]{Lemma}
\newtheorem{corollary}[theorem]{Corollary}

\theoremstyle{definition}
\newtheorem{definition}[theorem]{Definition}
\newtheorem{example}[theorem]{Example}

\theoremstyle{remark}
\newtheorem{remark}[theorem]{Remark}

\def\cal#1{\mathcal{#1}}

\def\th{^\text{th}}

\begin{document}
\begin{abstract}
It is well known that the existence of a braiding in a monoidal category ${\cal V}$
 allows many
 higher structures to be built upon that foundation. These include a monoidal 2-category ${\cal V}$-Cat
 of enriched
 categories and functors over ${\cal V}$, a monoidal bicategory ${\cal V}$-Mod of enriched
 categories and modules,
 a category of operads in ${\cal V}$ and
 a 2-fold monoidal category structure on ${\cal V}$. These all rely on the braiding to provide the
 existence of an interchange morphism $\eta$ necessary for either their structure or its properties.
 We ask, given a braiding on ${\cal V}$,
 what  non-equal structures of a given kind from this list exist
 which are based upon the braiding. For  example, what non-equal monoidal structures are available on ${\cal V}$-Cat, or what non-equal
 operad structures are available which base their associative structure on the braiding in ${\cal V}$.
 The basic question is the same as asking what non-equal 2-fold monoidal structures exist on a given braided
 category.
The main results are that the possible 2-fold monoidal structures are classified by a particular
set of four strand braids which we completely characterize, and that these 2-fold monoidal
categories are divided into two equivalence classes by the relation of 2-fold monoidal equivalence.
 \end{abstract}

\maketitle

\section{Introduction}
There are several levels of connection between the categorical concepts of braiding and
interchange. The first study of these ideas was by Joyal and Street in \cite{JS}. They point out
that a second tensor product in a given category which is a monoidal functor with respect to the
first gives rise to a braiding, and vice-versa. Most recently the work of Balteanu,  Fiedorowicz,
Schw${\rm \ddot a}$nzl, and Vogt in \cite{Balt} includes description of the same correspondence in
the context of lax monoidal functors. The connection between the $n$-fold monoidal categories in
\cite{Balt} and the theory of higher categories is through the  periodic table as laid out in
\cite{Baez1}.  Here Baez and Dolan organize the $k$-tuply monoidal $n$-categories, by which
terminology they refer to $(n+k)$-categories  that are trivial below dimension $k.$ The triviality
of lower cells allows the higher ones to compose freely, and thus these special cases of
$(n+k)$-categories are viewed as $n$-categories with $k$ multiplications.

 A $k$-tuply monoidal $n$-category is
a special $k$-fold monoidal $n$-category. The specialization results from the definition(s) of
$n$-category, all of which seem to include the axiom that the interchange transformation between
two ways of composing four higher morphisms along two different lower dimensions is required to be
an isomorphism. In \cite{Balt} the $k$-fold monoidal categories have interchange transformations
that are not isomorphisms. If those transformations are indeed isomorphisms then the $k$-fold
monoidal 1-categories do reduce to the braided and symmetric 1-categories of the periodic table.
Whether this continues for higher dimensions, yielding for example the sylleptic monoidal
2-categories of the periodic table as 3-fold monoidal 2-categories with interchange isomorphisms,
is an open question.

The key requirement of a 2-fold monoidal structure on a category is that a second tensor product
(in the sense of \cite{JS}) must be a functor which preserves the structure of the first tensor
product. Technically we say that the second is a monoidal functor with respect to the first. When
the two tensor products are identical, this translates into the existence of a coherent interchange
transformation $\eta_{ABCD}: (A\otimes B)\otimes (C\otimes D) \to (A\otimes C)\otimes (B\otimes
D)$. The chief goal of this paper is to study and classify the braids on four strands which can
play the role of an interchange transformation in a braided category seen as a 2-fold monoidal
category. To be precise, given a braided category  $({\cal V},\otimes,\alpha,c,I)$ (with strict
units, a strong associator $\alpha$, and braiding $c$), we ask the central question: For which
four-strand braids $b$ does the category ${\cal V}$ have in general a coherent 2-fold monoidal
structure, when that structure has $\otimes_1 = \otimes_2 = \otimes$ as functors, has $\alpha^1 =
\alpha^2 = \alpha$ as natural transformations, has strict unit $I$ for both identical tensor
products, and has $x$ as the underlying braid of $\eta$?

For example, the standard choice of $\eta = 1\otimes c_{BC} \otimes 1$ (where $\otimes$ is
associative) corresponds to the braid $x = \xy (0,0) *{ \xy 0;/r1pc/: +(-2,0)
   ,{\xcapv-@(0)}+(1,1),{\xoverv}+(2,1),{\xcapv-@(0)}
   \endxy}\endxy~.~$
There is a canonical epimorphism
  $\sigma : B_n \to S_n$ of the braid group on $n$ strands onto the permutation group. The permutation given by
  $\sigma$ is that
  given by the strands of the braid on the $n$ original positions. For instance on a standard generator of $B_n$,
  $\sigma_i$, we have $\sigma(\sigma_i) = (i\text{ }\text{ }i+1)$.
  Candidates for interchange would seem to be those defined using any braid $x \in B_4$
  such that $\sigma(x) = (2~3)$. However, it will be seen that many braids which at first seem to accomplish the same
   interchange (transposing the middle two terms and nothing else)
   do not in fact correspond to any $\eta$  which in general satisfies the axioms making
    ${\cal V}$ into a 2-fold monoidal category.
     For contrast, here is a braid
   which turns out to exemplify this failure: $x = \xy 0;/r1pc/:
   ,{\xcapv-@(0)}+(1,1),{\xoverv}+(2,1),{\xcapv-@(0)}+(-3,0)
   ,{\xcapv-@(0)}+(1,1),{\xoverv}+(2,1),{\xcapv-@(0)}+(-3,0)
   ,{\xcapv-@(0)}+(1,1),{\xoverv}+(2,1),{\xcapv-@(0)}+(-3,0)
   \endxy~.$

Under the hypotheses in this central question, there are four more requirements on the braid $x$
which must be filled for the answer to be affirmative. We refer to them, in parallel to the axioms
in \cite{Balt}, as
 internal and external unit conditions, and internal
and external associativity conditions. This nomenclature is the same as for the corresponding
axioms of a 2-fold monoidal category, which we will give in full in Section~\ref{2fold}. We refer
to the strands of a braid by their initial positions.  A sub-braid will refer to the braid
resulting from the deletion of a subset of the strands of a braid.

The internal unit condition on the four-strand  braid $x$ is that the sub-braids resulting from
deleting either the first two or the last two strands are both the identity (trivial braid) on two
strands. The external unit condition is that the sub-braids resulting from deleting
      either the first and third strand or the second and fourth
    strand are again the identity on two strands.

Internal associativity is pictured as an equality in $B_6$ of two six-strand braids derived from
the original four strand braid.
 We call the two derived six strand
braids $Lx$ and $Rx$. $Lx$ is algorithmically described as: (a) performing a copy of $x$ on the
first 4 strands; (b) grouping the new first and second pairs as if the edges of two ribbons--the
two pairs are actually (1,3) and (2,4); and (c) performing a copy of $x$ on the four new
``strands''-- the two ribbons along with the remaining two strands 5 and 6. $Rx$ is described by
three similar steps, but the initial
    copy of $x$ is on the last 4 strands, and the following ribbon edge pairing is on the pairs (3,5) and
    (4,6). The required equality, where
 $\sigma(x) = (2~3)$, is pictured here. Shading between two strands represents the ribbon edge pairing:

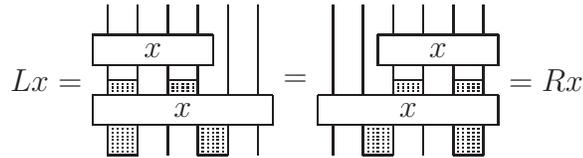
\[\xy (0,0) *{ Lx=~} \endxy
\xy (0,0) *{ \xy
  (-10,16)*{}="b1"; (-6,16)*{}="b2";
  (-2,16)*{}="b3"; (2,16)*{}="b4";
  (6,16)*{}="b5"; (10,16)*{}="b6";
(-12,12)*{}="l0"; (-10,12)*{}="l1";
  (-6,12)*{}="l2"; (-2,12)*{}="l3";
  (2,12)*{}="l4"; (4,12)*{}="l5";
 (-12,8)*{}="ll0"; (-10,8)*{}="ll1";
  (-6,8)*{}="ll2"; (-2,8)*{}="ll3";
  (2,8)*{}="ll4"; (4,8)*{}="ll5";
  (-10,6)*{}="s1"; (-6,6)*{}="s2";
  (-2,6)*{}="s3"; (2,6)*{}="s4";
(-10,5)*{}="ss1"; (-6,5)*{}="ss2";
  (-2,5)*{}="ss3"; (2,5)*{}="ss4";
  (-10,4.5)*{}="sa1"; (-6,4.5)*{}="sa2";
  (-2,4.5)*{}="sa3"; (2,4.5)*{}="sa4";
(-10,5.5)*{}="ssa1"; (-6,5.5)*{}="ssa2";
  (-2,5.5)*{}="ssa3"; (2,5.5)*{}="ssa4";
  (-12,4)*={}="bb0";
  (-10,4)*{}="bb1"; (-6,4)*{}="bb2";
  (-2,4)*{}="bb3"; (2,4)*{}="bb4";
  (6,4)*{}="bb5"; (10,4)*{}="bb6";
  (12,4)*={}="bb7";
(-12,0)*={}="Tb0";
  (-10,0)*{}="Tb1"; (-6,0)*{}="Tb2";
  (-2,0)*{}="Tb3"; (2,0)*{}="Tb4";
  (6,0)*{}="Tb5"; (10,0)*{}="Tb6";
  (12,0)*={}="Tb7";
  (-10,-1)*{}="sT1"; (-6,-1)*{}="sT2";
  (2,-1)*{}="sT4"; (6,-1)*{}="sT5";
  (-10,-2)*{}="ssT1"; (-6,-2)*{}="ssT2";
   (2,-2)*{}="ssT4"; (6,-2)*{}="ssT5";
  (-10,-3)*{}="sssT1"; (-6,-3)*{}="sssT2";
  (2,-3)*{}="sssT4"; (6,-3)*{}="sssT5";
  (-10,-.5)*{}="aT1"; (-6,-.5)*{}="aT2";
  (2,-.5)*{}="aT4"; (6,-.5)*{}="aT5";
  (-10,-1.5)*{}="saT1"; (-6,-1.5)*{}="saT2";
  (2,-1.5)*{}="saT4"; (6,-1.5)*{}="saT5";
  (-10,-2.5)*{}="ssaT1"; (-6,-2.5)*{}="ssaT2";
   (2,-2.5)*{}="ssaT4"; (6,-2.5)*{}="ssaT5";
  (-10,-3.5)*{}="sssaT1"; (-6,-3.5)*{}="sssaT2";
  (2,-3.5)*{}="sssaT4"; (6,-3.5)*{}="sssaT5";
(-10,-4)*{}="T1"; (-6,-4)*{}="T2";
  (-2,-4)*{}="T3"; (2,-4)*{}="T4";
  (6,-4)*{}="T5"; (10,-4)*{}="T6";
{\ar@{-} "b1";"l1"}; {\ar@{-} "b2";"l2"}; {\ar@{-} "b3";"l3"}; {\ar@{-} "b4";"l4"};
  {\ar@{-} "b5";"bb5"}; {\ar@{-} "b6";"bb6"};
  {\ar@{-} "bb1";"ll1"}; {\ar@{-} "bb2";"ll2"}; {\ar@{-} "bb3";"ll3"}; {\ar@{-} "bb4";"ll4"};
  {\ar@{-} "Tb1";"T1"}; {\ar@{-} "Tb2";"T2"}; {\ar@{-} "Tb3";"T3"}; {\ar@{-} "Tb4";"T4"};
  {\ar@{-} "Tb5";"T5"}; {\ar@{-} "Tb6";"T6"};
{\ar@{-} "l0";"ll0"}; {\ar@{-} "l5";"ll5"}; {\ar@{-} "Tb0";"bb0"}; {\ar@{-} "Tb7";"bb7"};
 {\ar@{-} "l0";"l5"}; {\ar@{-} "ll0";"ll5"}; {\ar@{-} "bb0";"bb7"}; {\ar@{-} "Tb0";"Tb7"};
 {\ar@{-} "s1";"s2"}; {\ar@{-} "s3";"s4"};
  {\ar@{.} "ss1";"ss2"}; {\ar@{.} "ss3";"ss4"};
   {\ar@{.} "sT1";"sT2"}; {\ar@{.} "sT4";"sT5"};
     {\ar@{.} "ssT1";"ssT2"}; {\ar@{.} "ssT4";"ssT5"};
       {\ar@{.} "sssT1";"sssT2"}; {\ar@{.} "sssT4";"sssT5"};
          {\ar@{-} "T1";"T2"}; {\ar@{-} "T4";"T5"};
{\ar@{.} "aT1";"aT2"}; {\ar@{.} "aT4";"aT5"};
 {\ar@{.} "sa1";"sa2"}; {\ar@{.} "sa3";"sa4"};
  {\ar@{.} "ssa1";"ssa2"}; {\ar@{.} "ssa3";"ssa4"};
   {\ar@{.} "saT1";"saT2"}; {\ar@{.} "saT4";"saT5"};
     {\ar@{.} "ssaT1";"ssaT2"}; {\ar@{.} "ssaT4";"ssaT5"};
       {\ar@{.} "sssaT1";"sssaT2"}; {\ar@{.} "sssaT4";"sssaT5"};
  (-4,10)*{x}; (0,2)*{x};
\endxy } \endxy
\xy (0,0)*{~=~} \endxy
 \xy (0,0) *{ \xy
  (10,16)*{}="b1"; (6,16)*{}="b2";
  (2,16)*{}="b3"; (-2,16)*{}="b4";
  (-6,16)*{}="b5"; (-10,16)*{}="b6";
(12,12)*{}="l0"; (10,12)*{}="l1";
  (6,12)*{}="l2"; (2,12)*{}="l3";
  (-2,12)*{}="l4"; (-4,12)*{}="l5";
 (12,8)*{}="ll0"; (10,8)*{}="ll1";
  (6,8)*{}="ll2"; (2,8)*{}="ll3";
  (-2,8)*{}="ll4"; (-4,8)*{}="ll5";
  (10,6)*{}="s1"; (6,6)*{}="s2";
  (2,6)*{}="s3"; (-2,6)*{}="s4";
(10,5)*{}="ss1"; (6,5)*{}="ss2";
  (2,5)*{}="ss3"; (-2,5)*{}="ss4";
(10,5.5)*{}="ssa1"; (6,5.5)*{}="ssa2";
  (2,5.5)*{}="ssa3"; (-2,5.5)*{}="ssa4";
  (10,4.5)*{}="sa1"; (6,4.5)*{}="sa2";
  (2,4.5)*{}="sa3"; (-2,4.5)*{}="sa4";
  (12,4)*={}="bb0";
  (10,4)*{}="bb1"; (6,4)*{}="bb2";
  (2,4)*{}="bb3"; (-2,4)*{}="bb4";
  (-6,4)*{}="bb5"; (-10,4)*{}="bb6";
  (-12,4)*={}="bb7";
(12,0)*={}="Tb0";
  (10,0)*{}="Tb1"; (6,0)*{}="Tb2";
  (2,0)*{}="Tb3"; (-2,0)*{}="Tb4";
  (-6,0)*{}="Tb5"; (-10,0)*{}="Tb6";
  (-12,0)*={}="Tb7";
  (10,-1)*{}="sT1"; (6,-1)*{}="sT2";
  (-2,-1)*{}="sT4"; (-6,-1)*{}="sT5";
  (10,-2)*{}="ssT1"; (6,-2)*{}="ssT2";
   (-2,-2)*{}="ssT4"; (-6,-2)*{}="ssT5";
  (10,-3)*{}="sssT1"; (6,-3)*{}="sssT2";
  (-2,-3)*{}="sssT4"; (-6,-3)*{}="sssT5";
  (10,-.5)*{}="aT1"; (6,-.5)*{}="aT2";
  (-2,-.5)*{}="aT4"; (-6,-.5)*{}="aT5";
  (10,-1.5)*{}="saT1"; (6,-1.5)*{}="saT2";
  (-2,-1.5)*{}="saT4"; (-6,-1.5)*{}="saT5";
  (10,-2.5)*{}="ssaT1"; (6,-2.5)*{}="ssaT2";
   (-2,-2.5)*{}="ssaT4"; (-6,-2.5)*{}="ssaT5";
  (10,-3.5)*{}="sssaT1"; (6,-3.5)*{}="sssaT2";
  (-2,-3.5)*{}="sssaT4"; (-6,-3.5)*{}="sssaT5";
(10,-4)*{}="T1"; (6,-4)*{}="T2";
  (2,-4)*{}="T3"; (-2,-4)*{}="T4";
  (-6,-4)*{}="T5"; (-10,-4)*{}="T6";
{\ar@{-} "b1";"l1"}; {\ar@{-} "b2";"l2"}; {\ar@{-} "b3";"l3"}; {\ar@{-} "b4";"l4"};
  {\ar@{-} "b5";"bb5"}; {\ar@{-} "b6";"bb6"};
  {\ar@{-} "bb1";"ll1"}; {\ar@{-} "bb2";"ll2"}; {\ar@{-} "bb3";"ll3"}; {\ar@{-} "bb4";"ll4"};
  {\ar@{-} "Tb1";"T1"}; {\ar@{-} "Tb2";"T2"}; {\ar@{-} "Tb3";"T3"}; {\ar@{-} "Tb4";"T4"};
  {\ar@{-} "Tb5";"T5"}; {\ar@{-} "Tb6";"T6"};
{\ar@{-} "l0";"ll0"}; {\ar@{-} "l5";"ll5"}; {\ar@{-} "Tb0";"bb0"}; {\ar@{-} "Tb7";"bb7"};
 {\ar@{-} "l0";"l5"}; {\ar@{-} "ll0";"ll5"}; {\ar@{-} "bb0";"bb7"}; {\ar@{-} "Tb0";"Tb7"};
 {\ar@{-} "s1";"s2"}; {\ar@{-} "s3";"s4"};
  {\ar@{.} "ss1";"ss2"}; {\ar@{.} "ss3";"ss4"};
   {\ar@{.} "sT1";"sT2"}; {\ar@{.} "sT4";"sT5"};
     {\ar@{.} "ssT1";"ssT2"}; {\ar@{.} "ssT4";"ssT5"};
       {\ar@{.} "sssT1";"sssT2"}; {\ar@{.} "sssT4";"sssT5"};
          {\ar@{-} "T1";"T2"}; {\ar@{-} "T4";"T5"};
{\ar@{.} "aT1";"aT2"}; {\ar@{.} "aT4";"aT5"};
 {\ar@{.} "sa1";"sa2"}; {\ar@{.} "sa3";"sa4"};
  {\ar@{.} "ssa1";"ssa2"}; {\ar@{.} "ssa3";"ssa4"};
   {\ar@{.} "saT1";"saT2"}; {\ar@{.} "saT4";"saT5"};
     {\ar@{.} "ssaT1";"ssaT2"}; {\ar@{.} "ssaT4";"ssaT5"};
       {\ar@{.} "sssaT1";"sssaT2"}; {\ar@{.} "sssaT4";"sssaT5"};
  (4,10)*{x}; (0,2)*{x};
\endxy} \endxy
\xy (0,0) *{~ =Rx} \endxy \]

External associativity is pictured as an equality in $B_6$ of two six-strand braids derived from
the original four strand braid.
 We call the two derived six strand
braids $L'x$ and $R'x$. $L'x$ is algorithmically described as: (a) pairing sets of strands (2,3)
and (5,6) as if the edges of two ribbons; (b) performing a copy of $x$  on the four new
``strands''-- the two ribbons along with the remaining two strands 1 and 4; and (c) performing a
copy of $x$  on the new  first four strands which are actually (2,3,5,6).  $R'b$ is similarly
described as (a) pairing sets of strands (1,2) and (4,5) as if the edges of two ribbons; (b)
performing a copy of $x$  on the four new ``strands''-- the two ribbons along with the remaining
two strands 3 and 6; and (c) performing a copy of $x$  on the new  last four strands which are
actually (1,2,4,5). The required equality is:

\[\xy (0,0) *{ L'x =~} \endxy
\begin{turn}{180}
\xy (0,0) *{ \xy
  (-10,16)*{}="b1"; (-6,16)*{}="b2";
  (-2,16)*{}="b3"; (2,16)*{}="b4";
  (6,16)*{}="b5"; (10,16)*{}="b6";
(-12,12)*{}="l0"; (-10,12)*{}="l1";
  (-6,12)*{}="l2"; (-2,12)*{}="l3";
  (2,12)*{}="l4"; (4,12)*{}="l5";
 (-12,8)*{}="ll0"; (-10,8)*{}="ll1";
  (-6,8)*{}="ll2"; (-2,8)*{}="ll3";
  (2,8)*{}="ll4"; (4,8)*{}="ll5";
  (-10,6)*{}="s1"; (-6,6)*{}="s2";
  (-2,6)*{}="s3"; (2,6)*{}="s4";
(-10,5)*{}="ss1"; (-6,5)*{}="ss2";
  (-2,5)*{}="ss3"; (2,5)*{}="ss4";
  (-10,4.5)*{}="sa1"; (-6,4.5)*{}="sa2";
  (-2,4.5)*{}="sa3"; (2,4.5)*{}="sa4";
(-10,5.5)*{}="ssa1"; (-6,5.5)*{}="ssa2";
  (-2,5.5)*{}="ssa3"; (2,5.5)*{}="ssa4";
  (-12,4)*={}="bb0";
  (-10,4)*{}="bb1"; (-6,4)*{}="bb2";
  (-2,4)*{}="bb3"; (2,4)*{}="bb4";
  (6,4)*{}="bb5"; (10,4)*{}="bb6";
  (12,4)*={}="bb7";
(-12,0)*={}="Tb0";
  (-10,0)*{}="Tb1"; (-6,0)*{}="Tb2";
  (-2,0)*{}="Tb3"; (2,0)*{}="Tb4";
  (6,0)*{}="Tb5"; (10,0)*{}="Tb6";
  (12,0)*={}="Tb7";
  (-10,-1)*{}="sT1"; (-6,-1)*{}="sT2";
  (2,-1)*{}="sT4"; (6,-1)*{}="sT5";
  (-10,-2)*{}="ssT1"; (-6,-2)*{}="ssT2";
   (2,-2)*{}="ssT4"; (6,-2)*{}="ssT5";
  (-10,-3)*{}="sssT1"; (-6,-3)*{}="sssT2";
  (2,-3)*{}="sssT4"; (6,-3)*{}="sssT5";
  (-10,-.5)*{}="aT1"; (-6,-.5)*{}="aT2";
  (2,-.5)*{}="aT4"; (6,-.5)*{}="aT5";
  (-10,-1.5)*{}="saT1"; (-6,-1.5)*{}="saT2";
  (2,-1.5)*{}="saT4"; (6,-1.5)*{}="saT5";
  (-10,-2.5)*{}="ssaT1"; (-6,-2.5)*{}="ssaT2";
   (2,-2.5)*{}="ssaT4"; (6,-2.5)*{}="ssaT5";
  (-10,-3.5)*{}="sssaT1"; (-6,-3.5)*{}="sssaT2";
  (2,-3.5)*{}="sssaT4"; (6,-3.5)*{}="sssaT5";
(-10,-4)*{}="T1"; (-6,-4)*{}="T2";
  (-2,-4)*{}="T3"; (2,-4)*{}="T4";
  (6,-4)*{}="T5"; (10,-4)*{}="T6";
{\ar@{-} "b1";"l1"}; {\ar@{-} "b2";"l2"}; {\ar@{-} "b3";"l3"}; {\ar@{-} "b4";"l4"};
  {\ar@{-} "b5";"bb5"}; {\ar@{-} "b6";"bb6"};
  {\ar@{-} "bb1";"ll1"}; {\ar@{-} "bb2";"ll2"}; {\ar@{-} "bb3";"ll3"}; {\ar@{-} "bb4";"ll4"};
  {\ar@{-} "Tb1";"T1"}; {\ar@{-} "Tb2";"T2"}; {\ar@{-} "Tb3";"T3"}; {\ar@{-} "Tb4";"T4"};
  {\ar@{-} "Tb5";"T5"}; {\ar@{-} "Tb6";"T6"};
{\ar@{-} "l0";"ll0"}; {\ar@{-} "l5";"ll5"}; {\ar@{-} "Tb0";"bb0"}; {\ar@{-} "Tb7";"bb7"};
 {\ar@{-} "l0";"l5"}; {\ar@{-} "ll0";"ll5"}; {\ar@{-} "bb0";"bb7"}; {\ar@{-} "Tb0";"Tb7"};
 {\ar@{-} "s1";"s2"}; {\ar@{-} "s3";"s4"};
  {\ar@{.} "ss1";"ss2"}; {\ar@{.} "ss3";"ss4"};
   {\ar@{.} "sT1";"sT2"}; {\ar@{.} "sT4";"sT5"};
     {\ar@{.} "ssT1";"ssT2"}; {\ar@{.} "ssT4";"ssT5"};
       {\ar@{.} "sssT1";"sssT2"}; {\ar@{.} "sssT4";"sssT5"};
          {\ar@{-} "T1";"T2"}; {\ar@{-} "T4";"T5"};
{\ar@{.} "aT1";"aT2"}; {\ar@{.} "aT4";"aT5"};
 {\ar@{.} "sa1";"sa2"}; {\ar@{.} "sa3";"sa4"};
  {\ar@{.} "ssa1";"ssa2"}; {\ar@{.} "ssa3";"ssa4"};
   {\ar@{.} "saT1";"saT2"}; {\ar@{.} "saT4";"saT5"};
     {\ar@{.} "ssaT1";"ssaT2"}; {\ar@{.} "ssaT4";"ssaT5"};
       {\ar@{.} "sssaT1";"sssaT2"}; {\ar@{.} "sssaT4";"sssaT5"};
  (-4,10)*{x}; (0,2)*{x};
\endxy } \endxy \end{turn}
\xy (0,0)*{~=~} \endxy
\begin{turn}{180}
 \xy (0,0) *{ \xy
  (10,16)*{}="b1"; (6,16)*{}="b2";
  (2,16)*{}="b3"; (-2,16)*{}="b4";
  (-6,16)*{}="b5"; (-10,16)*{}="b6";
(12,12)*{}="l0"; (10,12)*{}="l1";
  (6,12)*{}="l2"; (2,12)*{}="l3";
  (-2,12)*{}="l4"; (-4,12)*{}="l5";
 (12,8)*{}="ll0"; (10,8)*{}="ll1";
  (6,8)*{}="ll2"; (2,8)*{}="ll3";
  (-2,8)*{}="ll4"; (-4,8)*{}="ll5";
  (10,6)*{}="s1"; (6,6)*{}="s2";
  (2,6)*{}="s3"; (-2,6)*{}="s4";
(10,5)*{}="ss1"; (6,5)*{}="ss2";
  (2,5)*{}="ss3"; (-2,5)*{}="ss4";
(10,5.5)*{}="ssa1"; (6,5.5)*{}="ssa2";
  (2,5.5)*{}="ssa3"; (-2,5.5)*{}="ssa4";
  (10,4.5)*{}="sa1"; (6,4.5)*{}="sa2";
  (2,4.5)*{}="sa3"; (-2,4.5)*{}="sa4";
  (12,4)*={}="bb0";
  (10,4)*{}="bb1"; (6,4)*{}="bb2";
  (2,4)*{}="bb3"; (-2,4)*{}="bb4";
  (-6,4)*{}="bb5"; (-10,4)*{}="bb6";
  (-12,4)*={}="bb7";
(12,0)*={}="Tb0";
  (10,0)*{}="Tb1"; (6,0)*{}="Tb2";
  (2,0)*{}="Tb3"; (-2,0)*{}="Tb4";
  (-6,0)*{}="Tb5"; (-10,0)*{}="Tb6";
  (-12,0)*={}="Tb7";
  (10,-1)*{}="sT1"; (6,-1)*{}="sT2";
  (-2,-1)*{}="sT4"; (-6,-1)*{}="sT5";
  (10,-2)*{}="ssT1"; (6,-2)*{}="ssT2";
   (-2,-2)*{}="ssT4"; (-6,-2)*{}="ssT5";
  (10,-3)*{}="sssT1"; (6,-3)*{}="sssT2";
  (-2,-3)*{}="sssT4"; (-6,-3)*{}="sssT5";
  (10,-.5)*{}="aT1"; (6,-.5)*{}="aT2";
  (-2,-.5)*{}="aT4"; (-6,-.5)*{}="aT5";
  (10,-1.5)*{}="saT1"; (6,-1.5)*{}="saT2";
  (-2,-1.5)*{}="saT4"; (-6,-1.5)*{}="saT5";
  (10,-2.5)*{}="ssaT1"; (6,-2.5)*{}="ssaT2";
   (-2,-2.5)*{}="ssaT4"; (-6,-2.5)*{}="ssaT5";
  (10,-3.5)*{}="sssaT1"; (6,-3.5)*{}="sssaT2";
  (-2,-3.5)*{}="sssaT4"; (-6,-3.5)*{}="sssaT5";
(10,-4)*{}="T1"; (6,-4)*{}="T2";
  (2,-4)*{}="T3"; (-2,-4)*{}="T4";
  (-6,-4)*{}="T5"; (-10,-4)*{}="T6";
{\ar@{-} "b1";"l1"}; {\ar@{-} "b2";"l2"}; {\ar@{-} "b3";"l3"}; {\ar@{-} "b4";"l4"};
  {\ar@{-} "b5";"bb5"}; {\ar@{-} "b6";"bb6"};
  {\ar@{-} "bb1";"ll1"}; {\ar@{-} "bb2";"ll2"}; {\ar@{-} "bb3";"ll3"}; {\ar@{-} "bb4";"ll4"};
  {\ar@{-} "Tb1";"T1"}; {\ar@{-} "Tb2";"T2"}; {\ar@{-} "Tb3";"T3"}; {\ar@{-} "Tb4";"T4"};
  {\ar@{-} "Tb5";"T5"}; {\ar@{-} "Tb6";"T6"};
{\ar@{-} "l0";"ll0"}; {\ar@{-} "l5";"ll5"}; {\ar@{-} "Tb0";"bb0"}; {\ar@{-} "Tb7";"bb7"};
 {\ar@{-} "l0";"l5"}; {\ar@{-} "ll0";"ll5"}; {\ar@{-} "bb0";"bb7"}; {\ar@{-} "Tb0";"Tb7"};
 {\ar@{-} "s1";"s2"}; {\ar@{-} "s3";"s4"};
  {\ar@{.} "ss1";"ss2"}; {\ar@{.} "ss3";"ss4"};
   {\ar@{.} "sT1";"sT2"}; {\ar@{.} "sT4";"sT5"};
     {\ar@{.} "ssT1";"ssT2"}; {\ar@{.} "ssT4";"ssT5"};
       {\ar@{.} "sssT1";"sssT2"}; {\ar@{.} "sssT4";"sssT5"};
          {\ar@{-} "T1";"T2"}; {\ar@{-} "T4";"T5"};
{\ar@{.} "aT1";"aT2"}; {\ar@{.} "aT4";"aT5"};
 {\ar@{.} "sa1";"sa2"}; {\ar@{.} "sa3";"sa4"};
  {\ar@{.} "ssa1";"ssa2"}; {\ar@{.} "ssa3";"ssa4"};
   {\ar@{.} "saT1";"saT2"}; {\ar@{.} "saT4";"saT5"};
     {\ar@{.} "ssaT1";"ssaT2"}; {\ar@{.} "ssaT4";"ssaT5"};
       {\ar@{.} "sssaT1";"sssaT2"}; {\ar@{.} "sssaT4";"sssaT5"};
  (4,10)*{x}; (0,2)*{x};
\endxy} \endxy \end{turn}
\xy (0,0) *{~ =R'x} \endxy \]

Example~\ref{LbRb} and Example~\ref{Lbprime} contain pictures of original and derived braids.
Notice that both $L'x = Lx$ and $R'x = Rx$ if $x$ is equal to its own 180 degree rotation.

 The classification in
 Theorem~\ref{main} of  four-strand braids which obey all four requirements turns out to be fairly simple.
   These braids, which can underlie a coherent
  interchange in a 2-fold monoidal category, we designate as \emph{(unital) interchanging braids.} The first main result
is that unital interchanging braids are precisely those given in terms of standard generators by
$$b_{n^{\pm}}=(\sigma_2\sigma_1\sigma_3\sigma_2)^{\pm n}\sigma_2^{\pm 1}(\sigma_1\sigma_3)^{\mp n}
$$  for $n$ a non-negative integer. For example $b_{2^+}$ appears as:
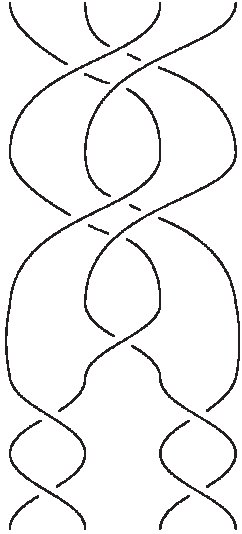
\[ \xy
  (15,10)*{}="b1"; (5,10)*{}="b2";
  (-5,10)*{}="b3"; (-15,10)*{}="b4";
  (15,-10)*{}="T1"; (5,-10)*{}="T2";
  (-5,-10)*{}="T3"; (-15,-10)*{}="T4";
  "b2"; "T4" **\crv{(5,2) & (-15,2.5)}
  \POS?(.35)*{\hole}="2x" \POS?(.54)*{\hole}="2y";
  "b4"; "2y" **\crv{(-15,6) };
  "b3"; "2x" **\crv{(-5,5) };
  "b1"; "T3" **\crv{(15,3.5) & (-5,3.5)}
  \POS?(.55)*{\hole}="4x" \POS?(.7)*{\hole}="4y";
    "4y";"2y" **\crv{(-2,-1)};
    "2x";"4x" **\crv{(3.5,1.75)};
  "4y";"T2" **\crv{(5,-3)};
    "T1";"4x" **\crv{(15,-2)};
  (15,-10)*{}="b1"; (5,-10)*{}="b2";
  (-5,-10)*{}="b3"; (-15,-10)*{}="b4";
  (15,-30)*{}="T1"; (5,-30)*{}="T2";
  (-5,-30)*{}="T3"; (-15,-30)*{}="T4";
  "b2"; "T4" **\crv{(5,-18) & (-14,-17.5)}
  \POS?(.35)*{\hole}="2x" \POS?(.54)*{\hole}="2y";
  "b4"; "2y" **\crv{(-15,-14) };
  "b3"; "2x" **\crv{(-5,-15) };
  "b1"; "T3" **\crv{(15,-16.5) & (-5,-16.5)}
  \POS?(.55)*{\hole}="4x" \POS?(.7)*{\hole}="4y";
    "4y";"2y" **\crv{(-2,-21)};
    "2x";"4x" **\crv{(3.5,-18.25)};
  "4y";"T2" **\crv{(5,-23)};
    "T1";"4x" **\crv{(14.25,-22)};
(15,-30)*{}="b1"; (5,-30)*{}="b2";
  (-5,-30)*{}="b3"; (-15,-30)*{}="b4";
  (15,-40)*{}="T1"; (5,-40)*{}="T2";
  (-5,-40)*{}="T3"; (-15,-40)*{}="T4";
  "b2"; "T3" **\crv{(5,-34) & (-5,-36)}
   \POS?(.5)*{\hole}="2x";
  "b4"; "T4" **\crv{(-16,-35) };
  "b3"; "2x" **\crv{(-5,-32.5) };
  "b1"; "T1" **\crv{(16,-35) };
     "2x"; "T2" **\crv{(5,-37.5)};
  (15,-40)*{}="b4"; (5,-40)*{}="b3";
  (-5,-40)*{}="b2"; (-15,-40)*{}="b1";
  (15,-50)*{}="T4"; (5,-50)*{}="T3";
  (-5,-50)*{}="T2"; (-15,-50)*{}="T1";
  "b1"; "T2" **\crv{(-14.25,-44) & (-5,-46)}
   \POS?(.5)*{\hole}="2x";
  "b3"; "T4" **\crv{(5,-44) & (15,-46)}
  \POS?(.5)*{\hole}="2y";
  "b2"; "2x" **\crv{(-5,-42.5) };
  "b4"; "2y" **\crv{(14.25,-42.5) };
  "2x"; "T1" **\crv{(-15,-47.5)};
  "2y"; "T3" **\crv{(5,-47.5)};
  (15,-50)*{}="b4"; (5,-50)*{}="b3";
  (-5,-50)*{}="b2"; (-15,-50)*{}="b1";
  (15,-60)*{}="T4"; (5,-60)*{}="T3";
  (-5,-60)*{}="T2"; (-15,-60)*{}="T1";
  "b1"; "T2" **\crv{(-15,-54) & (-5,-56)}
   \POS?(.5)*{\hole}="2x";
  "b3"; "T4" **\crv{(5,-54) & (15,-56)}
  \POS?(.5)*{\hole}="2y";
  "b2"; "2x" **\crv{(-5,-52.5) };
  "b4"; "2y" **\crv{(15,-52.5) };
  "2x"; "T1" **\crv{(-15,-57.5)};
  "2y"; "T3" **\crv{(5,-57.5)};
\endxy \]

Several geometrical facts can be observed about these unital interchanging braids. First, the braid
$b_{n^{\pm}}$ is a special
 element from the double coset $H\sigma_2^{\pm 1}K$ of the braid group on four strands. Here
 $H$ is the subgroup generated by the braid $\sigma_2\sigma_1\sigma_3\sigma_2$ and  $K$ is the
subgroup generated by $\sigma_1\sigma_3.$ Second, the braids $b_{n^{\pm}}$ are each equal to their
own 180 degree rotations. For example we have

$$ 
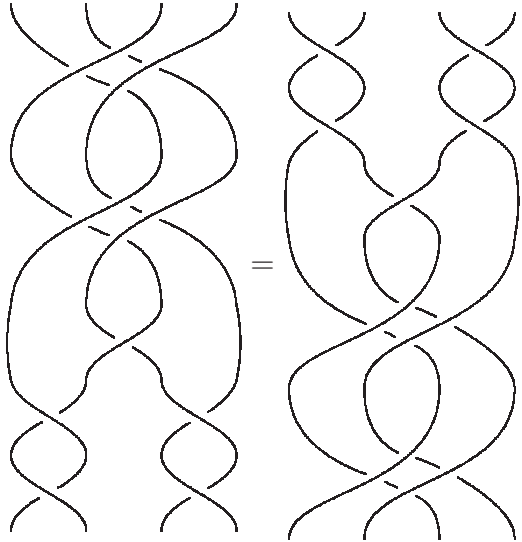
\xy (0,0) *{ \xy
  (15,10)*{}="b1"; (5,10)*{}="b2";
  (-5,10)*{}="b3"; (-15,10)*{}="b4";
  (15,-10)*{}="T1"; (5,-10)*{}="T2";
  (-5,-10)*{}="T3"; (-15,-10)*{}="T4";
  "b2"; "T4" **\crv{(5,2) & (-15,2.5)}
  \POS?(.35)*{\hole}="2x" \POS?(.54)*{\hole}="2y";
  "b4"; "2y" **\crv{(-15,6) };
  "b3"; "2x" **\crv{(-5,5) };
  "b1"; "T3" **\crv{(15,3.5) & (-5,3.5)}
  \POS?(.55)*{\hole}="4x" \POS?(.7)*{\hole}="4y";
    "4y";"2y" **\crv{(-2,-1)};
    "2x";"4x" **\crv{(3.5,1.75)};
  "4y";"T2" **\crv{(5,-3)};
    "T1";"4x" **\crv{(15,-2)};
  (15,-10)*{}="b1"; (5,-10)*{}="b2";
  (-5,-10)*{}="b3"; (-15,-10)*{}="b4";
  (15,-30)*{}="T1"; (5,-30)*{}="T2";
  (-5,-30)*{}="T3"; (-15,-30)*{}="T4";
  "b2"; "T4" **\crv{(5,-18) & (-14,-17.5)}
  \POS?(.35)*{\hole}="2x" \POS?(.54)*{\hole}="2y";
  "b4"; "2y" **\crv{(-15,-14) };
  "b3"; "2x" **\crv{(-5,-15) };
  "b1"; "T3" **\crv{(15,-16.5) & (-5,-16.5)}
  \POS?(.55)*{\hole}="4x" \POS?(.7)*{\hole}="4y";
    "4y";"2y" **\crv{(-2,-21)};
    "2x";"4x" **\crv{(3.5,-18.25)};
  "4y";"T2" **\crv{(5,-23)};
    "T1";"4x" **\crv{(14.25,-22)};
(15,-30)*{}="b1"; (5,-30)*{}="b2";
  (-5,-30)*{}="b3"; (-15,-30)*{}="b4";
  (15,-40)*{}="T1"; (5,-40)*{}="T2";
  (-5,-40)*{}="T3"; (-15,-40)*{}="T4";
  "b2"; "T3" **\crv{(5,-34) & (-5,-36)}
   \POS?(.5)*{\hole}="2x";
  "b4"; "T4" **\crv{(-16,-35) };
  "b3"; "2x" **\crv{(-5,-32.5) };
  "b1"; "T1" **\crv{(16,-35) };
     "2x"; "T2" **\crv{(5,-37.5)};
  (15,-40)*{}="b4"; (5,-40)*{}="b3";
  (-5,-40)*{}="b2"; (-15,-40)*{}="b1";
  (15,-50)*{}="T4"; (5,-50)*{}="T3";
  (-5,-50)*{}="T2"; (-15,-50)*{}="T1";
  "b1"; "T2" **\crv{(-14.25,-44) & (-5,-46)}
   \POS?(.5)*{\hole}="2x";
  "b3"; "T4" **\crv{(5,-44) & (15,-46)}
  \POS?(.5)*{\hole}="2y";
  "b2"; "2x" **\crv{(-5,-42.5) };
  "b4"; "2y" **\crv{(14.25,-42.5) };
  "2x"; "T1" **\crv{(-15,-47.5)};
  "2y"; "T3" **\crv{(5,-47.5)};
  (15,-50)*{}="b4"; (5,-50)*{}="b3";
  (-5,-50)*{}="b2"; (-15,-50)*{}="b1";
  (15,-60)*{}="T4"; (5,-60)*{}="T3";
  (-5,-60)*{}="T2"; (-15,-60)*{}="T1";
  "b1"; "T2" **\crv{(-15,-54) & (-5,-56)}
   \POS?(.5)*{\hole}="2x";
  "b3"; "T4" **\crv{(5,-54) & (15,-56)}
  \POS?(.5)*{\hole}="2y";
  "b2"; "2x" **\crv{(-5,-52.5) };
  "b4"; "2y" **\crv{(15,-52.5) };
  "2x"; "T1" **\crv{(-15,-57.5)};
  "2y"; "T3" **\crv{(5,-57.5)}; \endxy}
\endxy
\xy (0,0)*{~=~ } \endxy
\begin{turn}{180}
\xy (0,0)*{ \xy
  (15,10)*{}="b1"; (5,10)*{}="b2";
  (-5,10)*{}="b3"; (-15,10)*{}="b4";
  (15,-10)*{}="T1"; (5,-10)*{}="T2";
  (-5,-10)*{}="T3"; (-15,-10)*{}="T4";
  "b2"; "T4" **\crv{(5,2) & (-15,2.5)}
  \POS?(.35)*{\hole}="2x" \POS?(.54)*{\hole}="2y";
  "b4"; "2y" **\crv{(-15,6) };
  "b3"; "2x" **\crv{(-5,5) };
  "b1"; "T3" **\crv{(15,3.5) & (-5,3.5)}
  \POS?(.55)*{\hole}="4x" \POS?(.7)*{\hole}="4y";
    "4y";"2y" **\crv{(-2,-1)};
    "2x";"4x" **\crv{(3.5,1.75)};
  "4y";"T2" **\crv{(5,-3)};
    "T1";"4x" **\crv{(15,-2)};
  (15,-10)*{}="b1"; (5,-10)*{}="b2";
  (-5,-10)*{}="b3"; (-15,-10)*{}="b4";
  (15,-30)*{}="T1"; (5,-30)*{}="T2";
  (-5,-30)*{}="T3"; (-15,-30)*{}="T4";
  "b2"; "T4" **\crv{(5,-18) & (-14,-17.5)}
  \POS?(.35)*{\hole}="2x" \POS?(.54)*{\hole}="2y";
  "b4"; "2y" **\crv{(-15,-14) };
  "b3"; "2x" **\crv{(-5,-15) };
  "b1"; "T3" **\crv{(15,-16.5) & (-5,-16.5)}
  \POS?(.55)*{\hole}="4x" \POS?(.7)*{\hole}="4y";
    "4y";"2y" **\crv{(-2,-21)};
    "2x";"4x" **\crv{(3.5,-18.25)};
  "4y";"T2" **\crv{(5,-23)};
    "T1";"4x" **\crv{(14.25,-22)};
(15,-30)*{}="b1"; (5,-30)*{}="b2";
  (-5,-30)*{}="b3"; (-15,-30)*{}="b4";
  (15,-40)*{}="T1"; (5,-40)*{}="T2";
  (-5,-40)*{}="T3"; (-15,-40)*{}="T4";
  "b2"; "T3" **\crv{(5,-34) & (-5,-36)}
   \POS?(.5)*{\hole}="2x";
  "b4"; "T4" **\crv{(-16,-35) };
  "b3"; "2x" **\crv{(-5,-32.5) };
  "b1"; "T1" **\crv{(16,-35) };
     "2x"; "T2" **\crv{(5,-37.5)};
  (15,-40)*{}="b4"; (5,-40)*{}="b3";
  (-5,-40)*{}="b2"; (-15,-40)*{}="b1";
  (15,-50)*{}="T4"; (5,-50)*{}="T3";
  (-5,-50)*{}="T2"; (-15,-50)*{}="T1";
  "b1"; "T2" **\crv{(-14.25,-44) & (-5,-46)}
   \POS?(.5)*{\hole}="2x";
  "b3"; "T4" **\crv{(5,-44) & (15,-46)}
  \POS?(.5)*{\hole}="2y";
  "b2"; "2x" **\crv{(-5,-42.5) };
  "b4"; "2y" **\crv{(14.25,-42.5) };
  "2x"; "T1" **\crv{(-15,-47.5)};
  "2y"; "T3" **\crv{(5,-47.5)};
  (15,-50)*{}="b4"; (5,-50)*{}="b3";
  (-5,-50)*{}="b2"; (-15,-50)*{}="b1";
  (15,-60)*{}="T4"; (5,-60)*{}="T3";
  (-5,-60)*{}="T2"; (-15,-60)*{}="T1";
  "b1"; "T2" **\crv{(-15,-54) & (-5,-56)}
   \POS?(.5)*{\hole}="2x";
  "b3"; "T4" **\crv{(5,-54) & (15,-56)}
  \POS?(.5)*{\hole}="2y";
  "b2"; "2x" **\crv{(-5,-52.5) };
  "b4"; "2y" **\crv{(15,-52.5) };
  "2x"; "T1" **\crv{(-15,-57.5)};
  "2y"; "T3" **\crv{(5,-57.5)};
\endxy } \endxy
\end{turn}
 $$

The best way to visualize this equality is to draw a rectangle around the ``center'' portion of the
left braid, where a single copy of $\sigma_2$ divides the two ``double stranded positive
crossings'' above from the two ``negative crossings in tandem'' below. Now imagine rotating this
rectangle out of the plane of the page so as to uncross the  upper double stranded crossings. After
one full rotation (in order to completely undo the double stranded crossings) the right hand braid
is achieved. A good exercise would involve asking that this geometric argument be made into an
inductive proof. The braid equality in terms of generators is
$$b_{n^{\pm}}=(\sigma_2\sigma_1\sigma_3\sigma_2)^{\pm n}\sigma_2^{\pm 1}(\sigma_1\sigma_3)^{\mp n}
= (\sigma_1\sigma_3)^{\mp n}\sigma_2^{\pm 1}(\sigma_2\sigma_1\sigma_3\sigma_2)^{\pm n}.$$
 Thus the
third resulting fact is that the braid $b_{n^{\pm}}$ is in the intersection $H\sigma_2^{\pm 1}K
\cap K\sigma_2^{\pm 1}H.$ It would be interesting to know whether or not the only braids in this
intersection are the braids $b_{n^{\pm}}.$ It would also be interesting to know what connection, if
any, there is to the the similar braid equalities which arise in the theory of tortile categories,
as in \cite{Shum}.

The rest of this paper proceeds as follows:
 In the second section we begin with a review of the
category ${\cal V}$-Cat of enriched categories over a braided category ${\cal V}$. This is due to
the fact that when ${\cal V}$ is braided then ${\cal V}$-Cat can be equipped with  a monoidal
structure. It turns  out that the central question of which braids can underlie a coherent
interchange is equivalent to the question of which braids can underlie the middle four interchange
of the composition morphisms for a tensor product on ${\cal V}$-Cat. To be precise, given a braided
category  $({\cal V},\otimes,\alpha,c,I)$ (with strict unit, a strong associator $\alpha$, and
braiding $c$), we ask the new question: For which four-strand braids $b$ does the 2-category ${\cal
V}$-Cat have in general a coherent monoidal structure, given canonical choices for the objects,
hom-objects, and unit morphisms of the tensor product of two enriched categories, and the canonical
choice for the associator in ${\cal V}$-Cat, and given $b$ as the underlying braid of the middle
four interchange $\eta$?

In \cite{forcey1} it is shown that the external and internal unit conditions of a 2-fold monoidal
category ${\cal V}$
  imply the unital nature of ${\cal V}$-Cat and the unit axioms for a product of ${\cal V}$-categories respectively.
  The external and internal associativity conditions imply respectively
  the ${\cal V}$-functoriality of the associator in  ${\cal V}$-Cat and
  the associativity of the composition morphisms for tensor products of ${\cal V}$-categories.
Here we actually move in the opposite direction of implication: in order to find evidence of
sufficiency of the conditions which a  braid must meet to be interchanging we find tensor products
on ${\cal V}$-Cat which have the interchanging braids underlying their middle four interchange.

 In
the third section we review the axioms of a 2-fold monoidal category and demonstrate the necessity
of the conditions for our main result. In the fourth section we ask which of the 2-fold monoidal
categories we have described as arising from a certain braided category are equivalent as 2-fold
monoidal categories. Our result is that the relation of equivalence of 2-fold monoidal categories
splits our interchanging braids into into two equivalence classes, represented by $\sigma_2$ and
$\sigma_2^{-1}.$ The fifth section gives a list of easily detectable obstructions which prevent a
braid from having the interchanging property--i.e. which prevent it from being equivalent to one of
the braids described by our main result. The value of our classification is principally to provide
a solid framework for proofs about structures based upon a braiding. It turns out that only certain
braids can correspond to interchange transformations. Thus our results can be used to provide cases
for proofs, either by treating all the cases up to braid equivalence or more often just by treating
representative cases of categorical equivalence classes.

As an example in the sixth section we generalize results mentioned by Joyal and Street. They point
out that the category of enriched categories over a braided category is in general not braided and
that taking the opposite is not an involution. We give a proof which uses knot theory to
demonstrate non-existence in general for all possible interchanges. In questions of classification
of structures in a specific monoidal category, our result on interchanges may be necessary in order
to construct a complete picture. In the last section we relate our results to classification of
operads in a braided category. Throughout we work in monoidal categories with a coherent strong
associator, where ``strong'' implies that the natural transformation in question is an isomorphism.
The units will however be strict.

Thanks are due to many whose time was contributed to the development of this paper. Included are
the referee, who offered excellent suggestions for improvement of clarity, as well as Imre Tuba,
Jesse Siehler, and Ross Street.

\section{Braiding and Enrichment}
First we briefly review the definition of a category enriched over a monoidal category ${\cal V}$.
  Enriched functors and
    enriched natural transformations make the collection of enriched categories into a
    2-category ${\cal V}$-Cat.
     The definitions and proofs
    can be found in more or less detail in \cite{Kelly} and \cite{EK1} and of course
    in \cite{MacLane}. Some are included here for easy reference.

    \begin{definition} A {\it monoidal category} is a category ${\cal V}$
      together with a functor
      $\otimes: {\cal V}\times{\cal V}\to{\cal V}$  and an object $I$ such that
      \begin{enumerate}
      \item $\otimes$ is  associative up to the coherent natural isomorphisms
      $$\alpha_{ABC}: (A\otimes B)\otimes C \to A\otimes (B\otimes C) $$ called associators.
       The coherence
      axiom is given by the usual commuting pentagonal diagram as in \cite{MacLane}.


    \item In this paper, $I$ is a strict $2$-sided unit for $\otimes$.
    \end{enumerate}
    \end{definition}
    \begin {definition} \label{V:Cat} A (small) ${\cal V}$ {\it-Category} ${\cal A}$ is a set $\left|{\cal A}\right|$ of
    {\it objects},
    a {\it hom-object} ${\cal A}(A,B) \in \left|{\cal V}\right|$ for
    each pair of objects of ${\cal A}$, a family of {\it composition morphisms} $M_{ABC}:{\cal A}(B,C)
    \otimes{\cal A}(A,B)\to{\cal A}(A,C)$ for each triple of objects, and an {\it identity element} $j_{A}:I\to{\cal A}(A,A)$ for each object.
    The composition morphisms are subject to the associativity axiom which states that the following pentagon commutes

          \noindent
                      \begin{center}
                  \resizebox{5.5in}{!}{
          $$
          \xymatrix@C=-15pt{
          &({\cal A}(C,D)\otimes {\cal A}(B,C))\otimes {\cal A}(A,B)\text{ }\text{ }
          \ar[rr]^{\scriptstyle \alpha}
          \ar[dl]^{\scriptstyle M \otimes 1}
          &&\text{ }\text{ }{\cal A}(C,D)\otimes ({\cal A}(B,C)\otimes {\cal A}(A,B))
          \ar[dr]^{\scriptstyle 1 \otimes M}&\\
          {\cal A}(B,D)\otimes {\cal A}(A,B)
          \ar[drr]^{\scriptstyle M}
          &&&&{\cal A}(C,D)\otimes {\cal A}(A,C)
          \ar[dll]^{\scriptstyle M}
          \\&&{\cal A}(A,D))&&&
          }$$
          }
                                \end{center}

    and to the unit axioms which state that both the triangles in the following diagram commute

    $$
      \xymatrix{
      I\otimes {\cal A}(A,B)
      \ar[rrd]^{=}
      \ar[dd]_{j_{B}\otimes 1}
      &&&&{\cal A}(A,B)\otimes I
      \ar[dd]^{1\otimes j_{A}}
      \ar[lld]^{=}\\
      &&{\cal A}(A,B)\\
      {\cal A}(B,B)\otimes {\cal A}(A,B)
      \ar[rru]^{M_{ABB}}
      &&&&{\cal A}(A,B)\otimes {\cal A}(A,A)
      \ar[llu]^{M_{AAB}}
      }
   $$

   \end{definition}

    If ${\cal V} = \mathbf{Set}$ then these diagrams are the usual category axioms.
    Basically, composition
    of morphisms is replaced by tensoring
    and the resulting diagrams are required to commute. The next two definitions exhibit this
    principle and are important since they give
    us the setting in which to construct a category of ${\cal V}$-categories.

   \begin{definition} \label{enriched:funct} For ${\cal V}$-categories
   ${\cal A}$ and ${\cal B}$, a ${\cal V}$-$functor$ $T:{\cal A}\to{\cal B}$ is a function
    $T:\left| {\cal A} \right| \to \left| {\cal B} \right|$ and a family of
    morphisms $T_{AB}:{\cal A}(A,B) \to {\cal B}(TA,TB)$ in ${\cal V}$ indexed by
    pairs $A,B \in \left| {\cal A} \right|$.
    The usual rules for a functor that state $T(f \circ g) = Tf \circ Tg$
    and $T1_{A} = 1_{TA}$ become in the enriched setting, respectively, the commuting diagrams

   $$
    \xymatrix{
    &{\cal A}(B,C)\otimes {\cal A}(A,B)
    \ar[rr]^{\scriptstyle M}
    \ar[d]^{\scriptstyle T \otimes T}
    &&{\cal A}(A,C)
    \ar[d]^{\scriptstyle T}&\\
    &{\cal B}(TB,TC)\otimes {\cal B}(TA,TB)
    \ar[rr]^{\scriptstyle M}
    &&{\cal B}(TA,TC)
    }
   $$
  and
   $$
    \xymatrix{
    &&{\cal A}(A,A)
    \ar[dd]^{\scriptstyle T_{AA}}\\
    I
    \ar[rru]^{\scriptstyle j_{A}}
    \ar[rrd]_{\scriptstyle j_{TA}}\\
    &&{\cal B}(TA,TA)
    }
   $$
  ${\cal V}$-functors can be composed to form a category called ${\cal V}$-Cat. This category
  is actually enriched over $\mathbf{Cat}$, the category of (small) categories with Cartesian product.
     \end{definition}
%
%
%

  \begin{definition} A {\it braiding} for a monoidal category ${\cal V}$ is a family of natural
  isomorphisms $c_{XY}: X \otimes Y \to Y \otimes X$
  such that the following diagrams commute. They are drawn next to their underlying braids.
  Recall that by ``underlying braid'' of a composite of braidings and associators between two
  products of the same $n$ objects, we refer to the $n$-strand braid with crossings corresponding
  precisely to each instance of the braiding and its inverse.
  \begin{enumerate}
  \item
  $$
  \xy 0;/r1pc/:
  ,{\vtwist}+(2,1),{\xcapv-@(0)}
  +(-2,0),{\xcapv-@(0)}+(1,1),{\vtwist}
  \endxy
  \xymatrix{
  &(X \otimes Y) \otimes Z
  \ar[dl]^{c_{XY} \otimes 1}
  \ar[r]^{\alpha_{XYZ}}
  &X \otimes (Y \otimes Z)
  \ar[dr]^{c_{X(Y \otimes Z)}}
  \\(Y \otimes X) \otimes Z
  \ar[dr]^{\alpha_{YXZ}}
  &&&(Y \otimes Z) \otimes X
  \ar[dl]^{\alpha_{YZX}}
  \\&Y \otimes (X \otimes Z)
  \ar[r]^{1\otimes c_{XZ}}
  &Y \otimes (Z \otimes X)
  }
  $$

  \item
  $$
  \xy 0;/r1pc/:
  ,{\xcapv-@(0)}+(1,1),{\vtwist}+(-1,0)
  ,{\vtwist}+(2,1),{\xcapv-@(0)}
  \endxy
  \xymatrix{
  &X \otimes (Y \otimes Z)
  \ar[dl]^{1 \otimes c_{YZ}}
  \ar[r]^{\alpha_{XYZ}^{-1}}
  &(X \otimes Y) \otimes Z
  \ar[dr]^{c_{(X \otimes Y)Z}}
  \\X \otimes (Z \otimes Y)
  \ar[dr]^{\alpha_{XZY}^{-1}}
  &&&Z \otimes (X \otimes Y)
  \ar[dl]^{\alpha_{ZXY}^{-1}}
  \\&(X \otimes Z) \otimes Y
  \ar[r]^{c_{XZ} \otimes 1}
  &(Z \otimes X) \otimes Y
  }
  $$
  \end{enumerate}

  A braided category is a monoidal category with a chosen braiding. We will assume a strict unit in
  the monoidal categories considered here which implies a strict respect for units by the braiding.
  That is, $c_{IA} = c_{AI} = 1_A.$
  \end{definition}
  Joyal and Street
  proved the coherence theorem
  for braided categories in \cite{JS}, an immediate corollary of which is that in a free
  braided category generated by a set of
  objects, a diagram commutes in general
   if and only if all legs having the same source and target have the same underlying braid.

  \begin{definition}\label{braiding_cl}
  A {\it symmetry} is a braiding such that the following diagram commutes
  $$
  \xymatrix{
  X \otimes Y
  \ar[rr]^1
  \ar[dr]^{c_{XY}}
  &&X \otimes Y
  \\
  &Y \otimes X
  \ar[ur]^{c_{YX}}
  }
  $$
  In other words $c_{XY}^{-1} = c_{YX}$. A symmetric category is a monoidal category with a chosen symmetry.
  \end{definition}
  As pointed out by Joyal and Street, it is true that
$c^{-1}$ is a  braiding whenever $c$ is. These two braidings are equivalent  if and only if $c$ is
a symmetry; see Theorem~\ref{coolio} for the proof of this. It should be noted that there is
immediately an obstruction to other potential braidings based on the original. For sake of
efficiency we use notation
 $c^n_{AB} = c_{AB} \circ c_{BA} \circ c_{AB} \circ \dots \circ c_{AB} $ where there are $n$ instances of $c.$
 It appears at first that if $c_{AB}$ is a
braiding then $c' = c^{\pm(2n+1)}$ is potentially a braiding for any  $n,$ but actually we find
that:
\begin{lemma}\label{dos}
 for $n \ge 1,$ $c' = c^{\pm(2n+1)}$ is a braiding if and only if $c$ is a symmetry. (In that case
$c^{\pm(2n+1)}$ is also a symmetry.) \end{lemma}

\begin{proof} The obstruction arises from the
the braided coherence theorem applied to the hexagonal diagrams with $c^{\pm(2n+1)}$ in place of
the original instances of $c$.
 Observe that, when we test the potential
braiding for $n=1$, the hexagonal diagram (1) of Definition~\ref{braiding_cl} has legs with the
following two underlying braids. As denoted, this is an inequality:
 \noindent
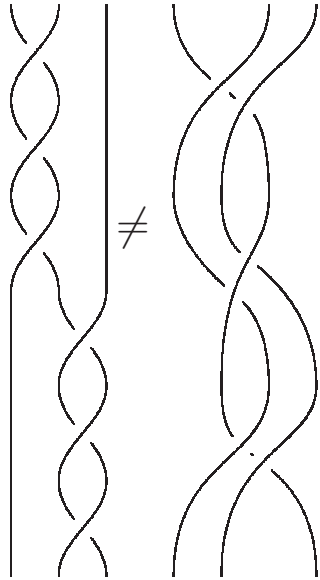
\begin{center}
$$ \xy 0;/r.15pc/:
  (15,10)*{}="b1"; (5,10)*{}="b2";
  (-5,10)*{}="b3";
  (15,-10)*{}="T1"; (5,-10)*{}="T2";
  (-5,-10)*{}="T3";
  "b2"; "T3" **\crv{(5,2) & (-5,-2)}
   \POS?(.5)*{\hole}="2x";
  "b3"; "2x" **\crv{(-5,5) };
  "b1"; "T1" **\crv{(15,0) };
     "2x"; "T2" **\crv{(5,-5)};
(15,-10)*{}="b1"; (5,-10)*{}="b2";
  (-5,-10)*{}="b3";
  (15,-30)*{}="T1"; (5,-30)*{}="T2";
  (-5,-30)*{}="T3";
  "b2"; "T3" **\crv{(5,-18) & (-5,-22)}
   \POS?(.5)*{\hole}="2x";
  "b3"; "2x" **\crv{(-5,-15) };
  "b1"; "T1" **\crv{(15,-20) };
     "2x"; "T2" **\crv{(5,-25)};
(15,-30)*{}="b1"; (5,-30)*{}="b2";
  (-5,-30)*{}="b3";
  (15,-50)*{}="T1"; (5,-50)*{}="T2";
  (-5,-50)*{}="T3";
  "b2"; "T3" **\crv{(5,-38) & (-5,-42)}
   \POS?(.5)*{\hole}="2x";
  "b3"; "2x" **\crv{(-5,-35) };
  "b1"; "T1" **\crv{(15,-40) };
     "2x"; "T2" **\crv{(5,-45)};
(15,-50)*{}="b1"; (5,-50)*{}="b2";
  (-5,-50)*{}="b3";
  (15,-70)*{}="T1"; (5,-70)*{}="T2";
  (-5,-70)*{}="T3";
  "b1"; "T2" **\crv{(15,-58) & (5,-62)}
   \POS?(.5)*{\hole}="2x";
  "b2"; "2x" **\crv{(5,-55) };
  "b3"; "T3" **\crv{(-5,-60) };
     "2x"; "T1" **\crv{(15,-65)};
(15,-70)*{}="b1"; (5,-70)*{}="b2";
  (-5,-70)*{}="b3";
  (15,-90)*{}="T1"; (5,-90)*{}="T2";
  (-5,-90)*{}="T3";
  "b1"; "T2" **\crv{(15,-78) & (5,-82)}
   \POS?(.5)*{\hole}="2x";
  "b2"; "2x" **\crv{(5,-75) };
  "b3"; "T3" **\crv{(-5,-80) };
     "2x"; "T1" **\crv{(15,-85)};
(15,-90)*{}="b1"; (5,-90)*{}="b2";
  (-5,-90)*{}="b3";
  (15,-110)*{}="T1"; (5,-110)*{}="T2";
  (-5,-110)*{}="T3";
  "b1"; "T2" **\crv{(15,-98) & (5,-102)}
   \POS?(.5)*{\hole}="2x";
  "b2"; "2x" **\crv{(5,-95) };
  "b3"; "T3" **\crv{(-5,-100) };
     "2x"; "T1" **\crv{(15,-105)};
  \endxy
\text{  }
 \xy 0;/r.11pc/: (0,-50)*{\text{\Large $\neq$ }}="T";
\endxy
\text{  }
   \xy 0;/r.15pc/:
  (15,10)*{}="b1"; (5,10)*{}="b2";
  (-15,10)*{}="b3"; (-15,10)*{}="b4";
  (5,-30)*{}="T1"; (5,-30)*{}="T2";
  (-5,-30)*{}="T3"; (-15,-30)*{}="T4";
  "b2"; "T4" **\crv{(5,-6) & (-15,-5)}
  \POS?(.51)*{\hole}="2x" ;
    "b3"; "2x" **\crv{(-15,0) };
  "b1"; "T3" **\crv{(15,-3) & (-5,-3)}
  \POS?(.68)*{\hole}="4x";
       "2x";"4x" **\crv{(-4.5,-7.5)};
     "T1";"4x" **\crv{(5,-14)};
  (5,-30)*{}="b1";
  (-5,-30)*{}="b3"; (-15,-30)*{}="b4";
  (15,-70)*{}="T1"; (5,-70)*{}="T2";
  (-5,-70)*{}="T3"; (-15,-70)*{}="T4";
    "b1"; "T3" **\crv{(5,-43) & (-5,-43)}
  \POS?(.45)*{\hole}="4x" \POS?(.67)*{\hole}="4y";
  "b4"; "4y" **\crv{(-15,-40) };
  "b3"; "4x" **\crv{(-5,-40) };
    "4y";"T2" **\crv{(5,-56)};
    "T1";"4x" **\crv{(15,-54)};
(15,-70)*{}="b1"; (5,-70)*{}="b2";
  (-5,-70)*{}="b3"; (-15,-70)*{}="b4";
  (15,-110)*{}="T1"; (5,-110)*{}="T2";
  (-5,-110)*{}="T3"; (-15,-110)*{}="T4";
  "b2"; "T4" **\crv{(5,-86) & (-15,-85)}
  \POS?(.35)*{\hole}="2x" ;
    "b3"; "2x" **\crv{(-5,-80) };
  "b1"; "T3" **\crv{(15,-83) & (-5,-83)}
  \POS?(.55)*{\hole}="4x" ;
       "2x";"4x" **\crv{(3.5,-86.5)};
     "T1";"4x" **\crv{(15,-94)};
\endxy $$
  \end{center}
Indeed we have that the required equality of braids for the first hexagonal axiom can never hold
for $n\ge 1.$ We check the positive powers of $c$ and note that the case for the negative powers is
shown similarly. For $c' = c^{2n+1}$ the braid inequality underlying the legs of the hexagonal
diagram, in terms of the standard braid generators, is $\sigma_1^{2n+1}\sigma_2^{2n+1} \neq
\sigma_1\sigma_2(\sigma_2\sigma_1\sigma_1\sigma_2)^{n }$. It is easy to see this inequality since
the semigroup of positive braids embeds into the braid group of the same number of strands, as
shown in \cite{Gar}. Thus any two positive braids  are equivalent in the braid group if and only if
they are equivalent in the positive semigroup, i.e. related by a chain of braid relations. For
three strand braids
 the only possible braid relation is the standard
$\sigma_1\sigma_2\sigma_1=\sigma_2\sigma_1\sigma_2.$ Note that in the braid words representing the
three strand braids in question there are no instances of either side of this relation, and so both
are in a unique positive form, and so clearly not equal. \end{proof}

  If ${\cal V}$ is braided then we can define additional structure on ${\cal V}$-Cat. The two classic structures
   are duality and tensor product.
\begin{definition}\label{canonop}
  First
  there is a left opposite of a ${\cal V}$-category which has
  $\left|{\cal A}^{op} \right| = \left|{\cal A}\right|$ and ${\cal A}^{op}(A,A') = {\cal A}(A',A).$ The
  composition morphisms are given by
  $$
    \xymatrix{
    {\cal A}^{op}(A',A'')\otimes {\cal A}^{op}(A,A')
    \ar@{=}[d]\\
    {\cal A}(A'',A')\otimes {\cal A}(A',A)
    \ar[d]_{c_{{\cal A}(A'',A')\otimes {\cal A}(A',A)}}\\
    {\cal A}(A',A)\otimes {\cal A}(A'',A')
    \ar[d]_{M_{AA'A''}}\\
    {\cal A}(A'',A)
    \ar@{=}[d]\\
    {\cal A}^{op}(A,A'')
    }
  $$
  The axiom for associativity of the composition morphisms in ${\cal A}^{op}$ holds due to the naturality of the
  braiding, the axiom for $M$ in ${\cal A}$, and the commutativity of a pentagonal diagram. This latter
  commutes since the braids underlying
  its legs are the two sides of the braid relation, also known as the Yang-Baxter equation.
  The unit morphisms in ${\cal A}^{op}$ are the same as the original $j_A:I\to {\cal A}(A,A)  = {\cal A}^{op}(A,A).$
  The unit axioms are obeyed due to the fact that $c_{IA} = c_{AI} = 1_A.$

  The right opposite denoted ${\cal A}^{po}$ is given by the same definition of composition and unit morphisms,
   but using $c^{-1}.$

   The two opposites take a ${\cal V}$-functor $F$ to its own function on objects
   but with $F^{(op)}_{AA'} = F^{(po)}_{AA'} = F_{A'A}.$ It is easy to check that thus they are functorial.
   That the image of a ${\cal V}$-functor under the opposites is still a ${\cal V}$-functor is due to the naturality of
   $c.$

   \end{definition}
  It is clear that $({\cal A}^{po})^{op}= ({\cal A}^{op})^{po} ={\cal A}.$
It is  also clear from this definition that $({\cal A}^{op})^{op}\ne {\cal A}$ in general unless
$c$ is a
  symmetry, and the same is true for the right opposite.
\begin{definition}\label{canontimes}
  The second structure is a
  tensor product for ${\cal V}$-Cat, that is, a 2-functor
  $$\otimes : {\cal V} \text{-Cat} \times {\cal V} \text{-Cat} \to {\cal V} \text{-Cat}.$$
  (In previous papers we have denoted the product(s) in ${\cal V}$-Cat
  with a superscript (1) in parentheses, but here it will be understood by context. The superscript (1) will
  still be used to denote that a given natural transformation is in ${\cal V}$-Cat.)
  The product of two ${\cal V}$-categories ${\cal A}$ and ${\cal B}$ has
  $\left|{\cal A} \otimes {\cal B}\right| = \left|{\cal A}\right| \times \left|{\cal B}\right|$ and
  $({\cal A} \otimes {\cal B})((A,B),(A',B')) = {\cal A}(A,A') \otimes {\cal B}(B,B').$

  The unit morphisms for the product ${\cal V}$-categories are the composites
  $$\xymatrix{I \cong I \otimes I \ar[r]_<<<<<{j_{A} \otimes j_{B}} & {\cal A}(A,A)\otimes {\cal B}(B,B) }$$


  The composition morphisms
  \begin{footnotesize}
  $$M_{(A,B)(A',B')(A'',B'')} : \Bigl({\cal A}\otimes {\cal B}\Bigr)\Bigl((A',B'),(A'',B'')\Bigr)\otimes \Bigl({\cal A}\otimes {\cal B}\Bigr)\Bigl((A,B),(A',B')\Bigr)\to \Bigl({\cal A}\otimes {\cal B}\Bigr)\Bigl((A,
  B),(A'',B'')\Bigr)$$
  \end{footnotesize}
  may be given canonically by
  $$
  \xymatrix{
  \Bigl({\cal A}\otimes {\cal B}\Bigr)\Bigl((A',B'),(A'',B'')\Bigr)\otimes \Bigl({\cal A}\otimes {\cal B}\Bigr)\Bigl((A,B),(A',B')\Bigr)
  \ar@{=}[d]\\
  \Bigl({\cal A}(A',A'')\otimes {\cal B}(B',B'')\Bigr)\otimes \Bigl({\cal A}(A,A')\otimes {\cal B}(B,B')\Bigr)
  \ar[d]_{(1 \otimes \alpha^{-1}) \circ \alpha}\\
  {\cal A}(A',A'')\otimes \Bigl(\bigl({\cal B}(B',B'')\otimes {\cal A}(A,A')\bigr)\otimes {\cal B}(B,B')\Bigr)
  \ar[d]_{1 \otimes (c_{{\cal B}(B',B''){\cal A}(A,A')} \otimes 1)}\\
  {\cal A}(A',A'')\otimes \Bigl(\bigl({\cal A}(A,A')\otimes {\cal B}(B',B'')\bigr) \otimes {\cal B}(B,B')\Bigr)
  \ar[d]_{\alpha^{-1} \circ (1 \otimes \alpha)}\\
  \Bigl({\cal A}(A',A'')\otimes {\cal A}(A,A')\Bigr)\otimes \Bigl({\cal B}(B',B'') \otimes {\cal B}(B,B')\Bigr)
  \ar[d]_{M_{AA'A''}\otimes M_{BB'B''}}\\
  {\cal A}(A,A'')\otimes {\cal B}(B,B'')
  \ar@{=}[d]\\
  \Bigl({\cal A}\otimes {\cal B}\Bigr)\Bigl((A,B),(A'',B'')\Bigr)
  }
  $$
\end{definition}
  That, in general, $({\cal A}\otimes {\cal B})^{op} \ne {\cal A}^{op}\otimes {\cal B}^{op}$  unless $c$ is a symmetry
  follows from the
  following braid inequality:
  $$
  \begin{turn}{180}
    \xy (0,0)*{ \xy 0;/r1pc/:
    ,{\xoverv}+(2,1),{\xoverv}+(-2,0)
    ,{\xcapv-@(0)}+(1,1),{\xoverv}+(2,1),{\xcapv-@(0)}+(-3,0)
    \endxy} \endxy
    \xy (0,0)*{~\ne~} \endxy
    \xy (0,0) *{ \xy 0;/r1pc/:
    ,{\xcapv-@(0)}+(1,1),{\xoverv}+(2,1),{\xcapv-@(0)}+(-3,0)
    ,{\xcapv-@(0)}+(1,1),{\xoverv}+(2,1),{\xcapv-@(0)}+(-3,0)
    ,{\xoverv}+(2,1),{\xoverv}+(-2,0)
    ,{\xcapv-@(0)}+(1,1),{\xoverv}+(2,1),{\xcapv-@(0)}+(-3,0)
    \endxy} \endxy
   \end{turn}
  $$

  Now consider more carefully the morphisms of ${\cal V}$ which make up the composition morphism for a tensor product
  of enriched categories, especially those
  which accomplish the ``middle four interchange'' \cite{Kelly} of the interior hom-objects, that is, all but the last
   pair of
  instances of the original composition $M$.
  In the symmetric case, any other combination of instances of $\alpha$ and $c$
  with the same domain and range would be equal, due to symmetric coherence.
  In the merely braided case, there at first seems to be a much larger range of available choices.
   Candidates for composition morphisms would seem to be those defined using any braid $b \in B_4$
  such that $\sigma(b) = (2~3)$.

  Thus a candidate for a new monoidal structure on ${\cal V}$-Cat
  could be given by the same canonical choices for objects, hom-objects,
  and unit morphisms as in Definition~\ref{canontimes} but with  alternate composition morphisms.
  The composition morphisms would be defined as above, but
  with the middle four interchange denoted $\eta_{(b)}$
  given by a
   series of instances of $\alpha$ and $c$
  such that the underlying braid is $b.$ Thus we might define
   $M_{(A,B)(A',B')(A'',B'')} = (M_{AA'A''}\otimes M_{BB'B''})\circ \eta_{(b)}.$
    That $M_{AA'A''}\otimes M_{BB'B''}$ will have the correct domain on which to operate is guaranteed by the
  permutation condition
  on $b$.

 Two important axioms that must hold for a proposed alternate monoidal structure on ${\cal V}$-Cat
  are associativity of composition $M$ (inside the proposed tensor product of two ${\cal V}$-categories)
  and ${\cal V}$-functoriality of the associator $\alpha$
  (so that there exists an associator for the proposed tensor product).
  For the associativity of composition to hold the following diagram must
    commute, where the first vertex is:
     $$
    \Bigl(\bigl({\cal A}\otimes {\cal B}\bigl)\bigl((A'',B''),(A''',B''')\bigr)\otimes \bigl({\cal A}\otimes {\cal B}\bigr)\bigl((A',B'),(A'',B'')\bigr)\Bigr)\otimes \bigl({\cal A}\otimes {\cal B}\bigr)\bigl((A,B),(A',B')\bigr)
    $$
    and the last bullet represents
    $\bigl({\cal A}\otimes {\cal B}\bigr)\bigl((A,B),(A''',B''')\bigr).$

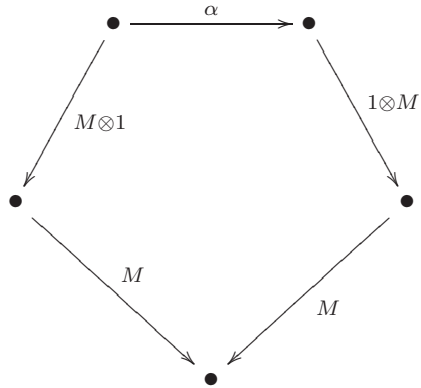
    $$
      \xymatrix{
      &\bullet
      \ar[rr]^{ \alpha}
      \ar[ddl]^{ M \otimes 1}
      &&\bullet
      \ar[ddr]^{ 1 \otimes M}&\\\\
      \bullet
      \ar[ddrr]^{ M}
      &&&&\bullet
      \ar[ddll]^{ M}
      \\\\&&\bullet
      }$$

    This means that the exterior of the following expanded diagram is required to commute.
    We leave out some parentheses
    for clarity and denote the middle four interchange by $\eta_{(b)}$ (perhaps composed with some associators).
    Also for convenience we write
    $X={\cal A}(A,A')$, $X'={\cal A}(A',A'')$, $X''={\cal A}(A'',A''')$, $Y={\cal B}(B,B')$,
     $Y'={\cal B}(B',B'')$ and $Y''={\cal B}(B'',B''').$

    \noindent
                          \begin{center}
                  \resizebox{6in}{!}{
    $$
    \xymatrix@C=-45pt{
    &[X''\otimes Y''\otimes X'\otimes Y']\otimes (X\otimes Y)
    \ar[dr]^{\alpha}
    \ar[dl]_{\eta_{(b)}}\\
    [X''\otimes X'\otimes Y''\otimes Y']\otimes (X\otimes Y)
    \ar[d]_{\alpha}
    &&(X''\otimes Y'')\otimes [X'\otimes Y'\otimes X\otimes Y]
    \ar[d]^{\eta_{(b)}}\\
    (X''\otimes X')\otimes (Y''\otimes Y')\otimes X\otimes Y
    \ar[d]_{\eta_{(b)}}
    &&(X''\otimes Y'')\otimes [X'\otimes X\otimes Y'\otimes Y]
    \ar[d]^{\alpha}
    \\
    [(X''\otimes X')\otimes X]\otimes [(Y''\otimes Y')\otimes Y]
        \ar[drr]^{\alpha \otimes \alpha}
    \ar[d]_{(M\otimes 1)\otimes (M \otimes 1)}
    &&X''\otimes Y''\otimes (X'\otimes X)\otimes (Y'\otimes Y)
    \ar[d]^{\eta_{(b)}}\\
    [{\cal A}(A',A''')\otimes X]\otimes [{\cal B}(B',B''')\otimes Y]
    \ar[d]_{M\otimes M}
    &&[X''\otimes (X'\otimes X)]\otimes [Y''\otimes (Y'\otimes Y)]
    \ar[d]^{(1\otimes M)\otimes (1 \otimes M)}\\
    {\cal A}(A,A''')\otimes {\cal B}(B,B''')
    &&[X''\otimes {\cal A}(A,A'')]\otimes [Y''\otimes {\cal B}(B,B'')]
    \ar[ll]^{M\otimes M}
    }
    $$
    }
    \end{center}

    The bottom region commutes by the associativity axioms for ${\cal A}$ and ${\cal B}.$
    We are left needing to show that the underlying braids are equal
    for the two legs of the upper region.
     In Example~\ref{LbRb} we give some examples of the underlying braids of the left and right
     legs
    for various choices of $b$. By inspection of the diagram these left and right underlying braids are the
    six-strand braids we denote respectively  $Lb$ and $Rb$.
    Recall from the introduction that  the requirement that $Lb=Rb$ is called internal associativity.
    The first example for $b$ is the one used in the original definition  of ${\cal A}\otimes {\cal B}$ given above.
    \begin{example}\label{LbRb}
    ~
    \end{example}
     $$
     \xy (0,0)*{
     b_{(1)}= b_{0^+}=~} \endxy
     \begin{turn}{180}
     \xy (0,0) *{
     \xy 0;/r1pc/: +(2,0)
     ,{\xcapv-@(0)}+(1,1),{\xoverv}+(2,1),{\xcapv-@(0)}
     \endxy } \endxy
     \end{turn}
     \xy (0,0)*{\text{ ; $Lb_{(1)} =~$}} \endxy
     \begin{turn}{180}
     \xy (0,0) *{
     \xy 0;/r1pc/: +(3,0)
        ,{\xcapv-@(0)}+(1,1),{\xcapv-@(0)}+(1,1),{\xoverv}+(2,1),{\xcapv-@(0)}+(1,1),{\xcapv-@(0)}+(-5,0)
        ,{\xcapv-@(0)}+(1,1),{\xoverv}+(2,1),{\xcapv-@(0)}+(1,1),{\xcapv-@(0)}+(1,1),{\xcapv-@(0)}+(-5,0)
        ,{\xcapv-@(0)}+(1,1),{\xcapv-@(0)}+(1,1),{\xcapv-@(0)}+(1,1),{\xoverv}+(2,1),{\xcapv-@(0)}
     \endxy }\endxy
     \end{turn}
     \xy (0,0) *{~ = ~} \endxy
     \begin{turn}{180}
     \xy (0,0) *{
     \xy 0;/r1pc/: +(3,0)
        ,{\xcapv-@(0)}+(1,1),{\xcapv-@(0)}+(1,1),{\xoverv}+(2,1),{\xcapv-@(0)}+(1,1),{\xcapv-@(0)}+(-5,0)
        ,{\xcapv-@(0)}+(1,1),{\xcapv-@(0)}+(1,1),{\xcapv-@(0)}+(1,1),{\xoverv}+(2,1),{\xcapv-@(0)}+(-5,0)
        ,{\xcapv-@(0)}+(1,1),{\xoverv}+(2,1),{\xcapv-@(0)}+(1,1),{\xcapv-@(0)}+(1,1),{\xcapv-@(0)}
     \endxy }\endxy
     \end{turn}
     \xy (0,0) *{~ = Rb_{(1)}.} \endxy
     $$

     $$
     \xy (0,0) *{
     b_{(2)}=~} \endxy
     \begin{turn}{180}
     \xy (0,0) *{
     \xy 0;/r1pc/:
     ,{\xcapv-@(0)}+(1,1),{\xoverv}+(2,1),{\xcapv-@(0)}+(-3,0)
     ,{\xoverv}+(2,1),{\xcapv-@(0)}+(1,1),{\xcapv-@(0)}+(-3,0)
     ,{\xoverv}+(2,1),{\xcapv-@(0)}+(1,1),{\xcapv-@(0)}
     \endxy } \endxy
     \end{turn}
     \xy (0,0) *{
     \text{ ; $Lb_{(2)} =~$ }}\endxy
     \begin{turn}{180}
     \xy (0,0) *{
     \xy 0;/r1pc/:
        ,{\xcapv-@(0)}+(1,1),{\xcapv-@(0)}+(1,1),{\xoverv}+(2,1),{\xcapv-@(0)}+(1,1),{\xcapv-@(0)}+(-5,0)
        ,{\xcapv-@(0)}+(1,1),{\xoverv}+(2,1),{\xcapv-@(0)}+(1,1),{\xcapv-@(0)}+(1,1),{\xcapv-@(0)}+(-5,0)
        ,{\xoverv}+(2,1),{\xcapv-@(0)}+(1,1),{\xcapv-@(0)}+(1,1),{\xcapv-@(0)}+(1,1),{\xcapv-@(0)}+(-5,0)
        ,{\xoverv}+(2,1),{\xcapv-@(0)}+(1,1),{\xcapv-@(0)}+(1,1),{\xcapv-@(0)}+(1,1),{\xcapv-@(0)}+(-5,0)
        ,{\xcapv-@(0)}+(1,1),{\xcapv-@(0)}+(1,1),{\xcapv-@(0)}+(1,1),{\xoverv}+(2,1),{\xcapv-@(0)}+(-5,0)
        ,{\xcapv-@(0)}+(1,1),{\xcapv-@(0)}+(1,1),{\xoverv}+(2,1),{\xcapv-@(0)}+(1,1),{\xcapv-@(0)}+(-5,0)
        ,{\xcapv-@(0)}+(1,1),{\xcapv-@(0)}+(1,1),{\xoverv}+(2,1),{\xcapv-@(0)}+(1,1),{\xcapv-@(0)}+(-5,0)
     \endxy} \endxy
     \end{turn}
     \xy (0,0) *{~ \ne ~}\endxy
     \begin{turn}{180}
     \xy (0,0) *{
     \xy 0;/r1pc/:
        ,{\xcapv-@(0)}+(1,1),{\xcapv-@(0)}+(1,1),{\xoverv}+(2,1),{\xcapv-@(0)}+(1,1),{\xcapv-@(0)}+(-5,0)
        ,{\xcapv-@(0)}+(1,1),{\xoverv}+(2,1),{\xoverv}+(2,1),{\xcapv-@(0)}+(-5,0)
        ,{\xoverv}+(2,1),{\xoverv}+(2,1),{\xcapv-@(0)}+(1,1),{\xcapv-@(0)}+(-5,0)
        ,{\xcapv-@(0)}+(1,1),{\xoverv}+(2,1),{\xcapv-@(0)}+(1,1),{\xcapv-@(0)}+(1,1),{\xcapv-@(0)}+(-5,0)
        ,{\xcapv-@(0)}+(1,1),{\xoverv}+(2,1),{\xcapv-@(0)}+(1,1),{\xcapv-@(0)}+(1,1),{\xcapv-@(0)}+(-5,0)
        ,{\xoverv}+(2,1),{\xoverv}+(2,1),{\xcapv-@(0)}+(1,1),{\xcapv-@(0)}+(-5,0)
        ,{\xcapv-@(0)}+(1,1),{\xoverv}+(2,1),{\xcapv-@(0)}+(1,1),{\xcapv-@(0)}+(1,1),{\xcapv-@(0)}+(-5,0)
        ,{\xcapv-@(0)}+(1,1),{\xoverv}+(2,1),{\xcapv-@(0)}+(1,1),{\xcapv-@(0)}+(1,1),{\xcapv-@(0)}+(-5,0)
        ,{\xoverv}+(2,1),{\xcapv-@(0)}+(1,1),{\xcapv-@(0)}+(1,1),{\xcapv-@(0)}+(1,1),{\xcapv-@(0)}+(-5,0)
        ,{\xoverv}+(2,1),{\xcapv-@(0)}+(1,1),{\xcapv-@(0)}+(1,1),{\xcapv-@(0)}+(1,1),{\xcapv-@(0)}+(-5,0)
     \endxy }\endxy
     \end{turn}
     \xy (0,0) *{~ = Rb_{(2)}.} \endxy
     $$

     $$
     \xy (0,0) *{
     b_{(3)}=~} \endxy
     \begin{turn}{180}
     \xy (0,0) *{
     \xy 0;/r1pc/:
     ,{\xcapv-@(0)}+(1,1),{\xcapv-@(0)}+(1,1),{\xoverv}+(-2,0)
     ,{\xcapv-@(0)}+(1,1),{\xcapv-@(0)}+(1,1),{\xoverv}+(-2,0)
     ,{\xcapv-@(0)}+(1,1),{\xoverv}+(2,1),{\xcapv-@(0)}
     \endxy }\endxy
     \end{turn}
     \xy (0,0) *{
     \text{ ; $Lb_{(3)} = ~$} }\endxy
     \begin{turn}{180}
     \xy (0,0) *{
     \xy 0;/r1pc/:
        ,{\xcapv-@(0)}+(1,1),{\xcapv-@(0)}+(1,1),{\xcapv-@(0)}+(1,1),{\xoverv}+(2,1),{\xcapv-@(0)}+(-5,0)
        ,{\xcapv-@(0)}+(1,1),{\xcapv-@(0)}+(1,1),{\xcapv-@(0)}+(1,1),{\xcapv-@(0)}+(1,1),{\xoverv}+(-4,0)
        ,{\xcapv-@(0)}+(1,1),{\xcapv-@(0)}+(1,1),{\xcapv-@(0)}+(1,1),{\xcapv-@(0)}+(1,1),{\xoverv}+(-4,0)
        ,{\xcapv-@(0)}+(1,1),{\xcapv-@(0)}+(1,1),{\xcapv-@(0)}+(1,1),{\xoverv}+(2,1),{\xcapv-@(0)}+(-5,0)
        ,{\xcapv-@(0)}+(1,1),{\xcapv-@(0)}+(1,1),{\xoverv}+(2,1),{\xcapv-@(0)}+(1,1),{\xcapv-@(0)}+(-5,0)
        ,{\xcapv-@(0)}+(1,1),{\xoverv}+(2,1),{\xcapv-@(0)}+(1,1),{\xcapv-@(0)}+(1,1),{\xcapv-@(0)}+(-5,0)
        ,{\xcapv-@(0)}+(1,1),{\xcapv-@(0)}+(1,1),{\xcapv-@(0)}+(1,1),{\xcapv-@(0)}+(1,1),{\xoverv}+(-4,0)
        ,{\xcapv-@(0)}+(1,1),{\xcapv-@(0)}+(1,1),{\xcapv-@(0)}+(1,1),{\xcapv-@(0)}+(1,1),{\xoverv}+(-4,0)
        ,{\xcapv-@(0)}+(1,1),{\xcapv-@(0)}+(1,1),{\xcapv-@(0)}+(1,1),{\xoverv}+(2,1),{\xcapv-@(0)}+(-5,0)
     \endxy } \endxy
     \end{turn}
     \xy (0,0) *{~=~}\endxy
     \begin{turn}{180}
     \xy (0,0) *{
     \xy 0;/r1pc/:
        ,{\xcapv-@(0)}+(1,1),{\xcapv-@(0)}+(1,1),{\xcapv-@(0)}+(1,1),{\xcapv-@(0)}+(1,1),{\xoverv}+(-4,0)
        ,{\xcapv-@(0)}+(1,1),{\xcapv-@(0)}+(1,1),{\xcapv-@(0)}+(1,1),{\xoverv}+(2,1),{\xcapv-@(0)}+(-5,0)
        ,{\xcapv-@(0)}+(1,1),{\xcapv-@(0)}+(1,1),{\xcapv-@(0)}+(1,1),{\xoverv}+(2,1),{\xcapv-@(0)}+(-5,0)
        ,{\xcapv-@(0)}+(1,1),{\xcapv-@(0)}+(1,1),{\xoverv}+(2,1),{\xoverv}+(-4,0)
        ,{\xcapv-@(0)}+(1,1),{\xcapv-@(0)}+(1,1),{\xcapv-@(0)}+(1,1),{\xoverv}+(2,1),{\xcapv-@(0)}+(-5,0)
        ,{\xcapv-@(0)}+(1,1),{\xcapv-@(0)}+(1,1),{\xoverv}+(2,1),{\xcapv-@(0)}+(1,1),{\xcapv-@(0)}+(-5,0)
        ,{\xcapv-@(0)}+(1,1),{\xcapv-@(0)}+(1,1),{\xoverv}+(2,1),{\xcapv-@(0)}+(1,1),{\xcapv-@(0)}+(-5,0)
        ,{\xcapv-@(0)}+(1,1),{\xoverv}+(2,1),{\xcapv-@(0)}+(1,1),{\xcapv-@(0)}+(1,1),{\xcapv-@(0)}
     \endxy } \endxy
     \end{turn}
     \xy (0,0) *{~ = Rb_{(3)}.} \endxy
     $$
     $$
     \xy (0,0) *{
     b_{(4)}=~ }\endxy
     \begin{turn}{180}
     \xy (0,0) *{
     \xy 0;/r1pc/:
     ,{\xcapv-@(0)}+(1,1),{\xoverv}+(2,1),{\xcapv-@(0)}+(-3,0)
     ,{\xcapv-@(0)}+(1,1),{\xoverv}+(2,1),{\xcapv-@(0)}+(-3,0)
     ,{\xcapv-@(0)}+(1,1),{\xoverv}+(2,1),{\xcapv-@(0)}+(-3,0)
     \endxy }\endxy
     \end{turn}
     \xy (0,0) *{
     \text{ ; $Lb_{(4)} =~$} }\endxy
     \begin{turn}{180}
     \xy (0,0) *{
     \xy 0;/r1pc/:
        ,{\xcapv-@(0)}+(1,1),{\xcapv-@(0)}+(1,1),{\xoverv}+(2,1),{\xcapv-@(0)}+(1,1),{\xcapv-@(0)}+(-5,0)
        ,{\xcapv-@(0)}+(1,1),{\xoverv}+(2,1),{\xcapv-@(0)}+(1,1),{\xcapv-@(0)}+(1,1),{\xcapv-@(0)}+(-5,0)
        ,{\xcapv-@(0)}+(1,1),{\xoverv}+(2,1),{\xcapv-@(0)}+(1,1),{\xcapv-@(0)}+(1,1),{\xcapv-@(0)}+(-5,0)
        ,{\xcapv-@(0)}+(1,1),{\xcapv-@(0)}+(1,1),{\xoverv}+(2,1),{\xcapv-@(0)}+(1,1),{\xcapv-@(0)}+(-5,0)
        ,{\xcapv-@(0)}+(1,1),{\xcapv-@(0)}+(1,1),{\xoverv}+(2,1),{\xcapv-@(0)}+(1,1),{\xcapv-@(0)}+(-5,0)
        ,{\xcapv-@(0)}+(1,1),{\xoverv}+(2,1),{\xcapv-@(0)}+(1,1),{\xcapv-@(0)}+(1,1),{\xcapv-@(0)}+(-5,0)
        ,{\xcapv-@(0)}+(1,1),{\xcapv-@(0)}+(1,1),{\xcapv-@(0)}+(1,1),{\xoverv}+(2,1),{\xcapv-@(0)}+(-5,0)
        ,{\xcapv-@(0)}+(1,1),{\xcapv-@(0)}+(1,1),{\xcapv-@(0)}+(1,1),{\xoverv}+(2,1),{\xcapv-@(0)}+(-5,0)
        ,{\xcapv-@(0)}+(1,1),{\xcapv-@(0)}+(1,1),{\xcapv-@(0)}+(1,1),{\xoverv}+(2,1),{\xcapv-@(0)}+(-5,0)
     \endxy }\endxy
     \end{turn}
     \xy (0,0) *{~\ne~}\endxy
     \begin{turn}{180}
     \xy (0,0) *{
     \xy 0;/r1pc/:
        ,{\xcapv-@(0)}+(1,1),{\xcapv-@(0)}+(1,1),{\xoverv}+(2,1),{\xcapv-@(0)}+(1,1),{\xcapv-@(0)}+(-5,0)
        ,{\xcapv-@(0)}+(1,1),{\xcapv-@(0)}+(1,1),{\xcapv-@(0)}+(1,1),{\xoverv}+(2,1),{\xcapv-@(0)}+(-5,0)
        ,{\xcapv-@(0)}+(1,1),{\xcapv-@(0)}+(1,1),{\xcapv-@(0)}+(1,1),{\xoverv}+(2,1),{\xcapv-@(0)}+(-5,0)
        ,{\xcapv-@(0)}+(1,1),{\xcapv-@(0)}+(1,1),{\xoverv}+(2,1),{\xcapv-@(0)}+(1,1),{\xcapv-@(0)}+(-5,0)
        ,{\xcapv-@(0)}+(1,1),{\xcapv-@(0)}+(1,1),{\xoverv}+(2,1),{\xcapv-@(0)}+(1,1),{\xcapv-@(0)}+(-5,0)
        ,{\xcapv-@(0)}+(1,1),{\xcapv-@(0)}+(1,1),{\xcapv-@(0)}+(1,1),{\xoverv}+(2,1),{\xcapv-@(0)}+(-5,0)
        ,{\xcapv-@(0)}+(1,1),{\xoverv}+(2,1),{\xcapv-@(0)}+(1,1),{\xcapv-@(0)}+(1,1),{\xcapv-@(0)}+(-5,0)
        ,{\xcapv-@(0)}+(1,1),{\xoverv}+(2,1),{\xcapv-@(0)}+(1,1),{\xcapv-@(0)}+(1,1),{\xcapv-@(0)}+(-5,0)
        ,{\xcapv-@(0)}+(1,1),{\xoverv}+(2,1),{\xcapv-@(0)}+(1,1),{\xcapv-@(0)}+(1,1),{\xcapv-@(0)}+(-5,0)
     \endxy}\endxy
     \end{turn}
    \xy (0,0) *{~ = Rb_{(4)}.} \endxy
     $$
     $$
     \xy (0,0) *{
b_{(5)}= b_{1^-} = ~}\endxy
 \xy (0,0) *{
      \xy 0;/r1pc/:
      ,{\xcapv-@(0)}+(1,1),{\xunderv}+(2,1),{\xcapv-@(0)}+(-3,0)
      ,{\xunderv}+(2,1),{\xunderv}+(-2,0)
      ,{\xcapv-@(0)}+(1,1),{\xunderv}+(2,1),{\xcapv-@(0)}+(-3,0)
      ,{\xcapv-@(0)}+(1,1),{\xunderv}+(2,1),{\xcapv-@(0)}+(-3,0)
      ,{\xoverv}+(2,1),{\xoverv}+(-3,0)
      \endxy } \endxy
      \xy (0,0) *{
      \text{ ; $Lb_{(5)} = ~$ }} \endxy
      \xy (0,0) *{
       \xy 0;/r1pc/:
         ,{\xcapv-@(0)}+(1,1),{\xunderv}+(2,1),{\xcapv-@(0)}+(1,1),{\xcapv-@(0)}+(1,1),{\xcapv-@(0)}+(-5,0)
      ,{\xunderv}+(2,1),{\xunderv}+(2,1),{\xcapv-@(0)}+(1,1),{\xcapv-@(0)}+(-5,0)
      ,{\xcapv-@(0)}+(1,1),{\xunderv}+(2,1),{\xcapv-@(0)}+(1,1),{\xcapv-@(0)}+(1,1),{\xcapv-@(0)}+(-5,0)
      ,{\xcapv-@(0)}+(1,1),{\xunderv}+(2,1),{\xcapv-@(0)}+(1,1),{\xcapv-@(0)}+(1,1),{\xcapv-@(0)}+(-5,0)
      ,{\xoverv}+(2,1),{\xoverv}+(2,1),{\xcapv-@(0)}+(1,1),{\xcapv-@(0)}+(-5,0)
      ,{\xcapv-@(0)}+(1,1),{\xcapv-@(0)}+(1,1),{\xcapv-@(0)}+(1,1),{\xunderv}+(2,1)                ,{\xcapv-@(0)}+(-5,0)
      ,{\xcapv-@(0)}+(1,1),{\xcapv-@(0)}+(1,1),{\xunderv}+(2,1)                  ,{\xunderv}+(-4,0)
      ,{\xcapv-@(0)}+(1,1),{\xunderv}+(2,1)                   ,{\xunderv}+(2,1)                     ,{\xcapv-@(0)}+(-5,0)
      ,{\xunderv}+(2,1)                    ,{\xunderv}+(2,1)                    ,{\xcapv-@(0)}+(1,1),{\xcapv-@(0)}+(-5,0)
      ,{\xcapv-@(0)}+(1,1),{\xunderv}+(2,1)                   ,{\xcapv-@(0)}+(1,1),{\xcapv-@(0)}+(1,1),{\xcapv-@(0)}+(-5,0)
      ,{\xcapv-@(0)}+(1,1),{\xunderv}+(2,1)                   ,{\xcapv-@(0)}+(1,1),{\xcapv-@(0)}+(1,1),{\xcapv-@(0)}+(-5,0)
      ,{\xoverv}+(2,1)                    ,{\xunderv}+(2,1)                    ,{\xcapv-@(0)}+(1,1),{\xcapv-@(0)}+(-5,0)
      ,{\xcapv-@(0)}+(1,1),{\xoverv}+(2,1)                      ,{\xoverv}+(2,1)                      ,{\xcapv-@(0)}+(-5,0)
      ,{\xcapv-@(0)}+(1,1),{\xcapv-@(0)}+(1,1),{\xcapv-@(0)}+(1,1),{\xcapv-@(0)}+(1,1),{\xoverv}+(-4,0)
      \endxy } \endxy
      \xy (0,0) *{~=~} \endxy
      \xy (0,0) *{
      \xy 0;/r1pc/:
      ,{\xcapv-@(0)}+(1,1),{\xcapv-@(0)}+(1,1),{\xcapv-@(0)}+(1,1),{\xunderv}+(2,1),{\xcapv-@(0)}+(-5,0)
       ,{\xcapv-@(0)}+(1,1),{\xcapv-@(0)}+(1,1),{\xunderv}+(2,1),{\xunderv}+(-4,0)
       ,{\xcapv-@(0)}+(1,1),{\xcapv-@(0)}+(1,1),{\xcapv-@(0)}+(1,1),{\xunderv}+(2,1),{\xcapv-@(0)}+(-5,0)
       ,{\xcapv-@(0)}+(1,1),{\xcapv-@(0)}+(1,1),{\xcapv-@(0)}+(1,1),{\xunderv}+(2,1),{\xcapv-@(0)}+(-5,0)
      ,{\xcapv-@(0)}+(1,1),{\xcapv-@(0)}+(1,1),{\xoverv}+(2,1),{\xoverv}+(-4,0)
         ,{\xcapv-@(0)}+(1,1),{\xunderv}+(2,1)           ,{\xcapv-@(0)}+(1,1),{\xcapv-@(0)}+(1,1),{\xcapv-@(0)}+(-5,0)
         ,{\xunderv}+(2,1)                    ,{\xunderv}+(2,1)                     ,{\xcapv-@(0)}+(1,1),{\xcapv-@(0)}+(-5,0)
         ,{\xcapv-@(0)}+(1,1),{\xunderv}+(2,1)           ,{\xunderv}+(2,1)                     ,{\xcapv-@(0)}+(-5,0)
         ,{\xcapv-@(0)}+(1,1),{\xcapv-@(0)}+(1,1),{\xunderv}+(2,1)                    ,{\xunderv}+(-4,0)
         ,{\xcapv-@(0)}+(1,1),{\xcapv-@(0)}+(1,1),{\xcapv-@(0)}+(1,1),{\xunderv}+(2,1)                   ,{\xcapv-@(0)}+(-5,0)
         ,{\xcapv-@(0)}+(1,1),{\xcapv-@(0)}+(1,1),{\xcapv-@(0)}+(1,1),{\xunderv}+(2,1)                   ,{\xcapv-@(0)}+(-5,0)
         ,{\xcapv-@(0)}+(1,1),{\xcapv-@(0)}+(1,1),{\xunderv}+(2,1)                    ,{\xoverv}+(-4,0)
         ,{\xcapv-@(0)}+(1,1),{\xoverv}+(2,1)           ,{\xoverv}+(2,1)                        ,{\xcapv-@(0)}+(-5,0)
         ,{\xoverv}+(2,1)                       ,{\xcapv-@(0)}+(1,1),{\xcapv-@(0)}+(1,1),{\xcapv-@(0)}+(1,1),{\xcapv-@(0)}+(-5,0)
          \endxy }\endxy
          \xy (0,0) *{~ = Rb_{(5)}.} \endxy
           $$

\begin{remark}\label{suff}
Before turning to check on ${\cal V}$-functoriality of the associator, we note that $b_{(3)}$ is
the braid underlying the composition morphism of the product category $({\cal A}^{op})^{op}\otimes
{\cal B}$ where the product is defined using $b_{(1)}.$ This provides the hint that the two derived
braids $Lb_{(3)}, Rb_{(3)} \in B_6$ are equal because of the fact that the opposite of a ${\cal
V}$-category is a valid ${\cal V}$-category. In fact we can describe sufficient conditions for $Lb$
to be equivalent to $Rb$ by describing the braids $b$ that underlie the composition morphism of a
product category given generally by $((({\cal A}^{op})^{...op}\otimes ({\cal
B}^{op})^{...op})^{op})^{...op}$ where the number of $op$ exponents is arbitrary in each position.

Those braids are alternately described as lying in $H\sigma_2K \subset B_4$ where $H$ is the cyclic
subgroup generated by the braid $\sigma_2\sigma_1\sigma_3\sigma_2$ and $K'$ is the subgroup
generated by the two generators $\{\sigma_1, \sigma_3\}.$ The latter subgroup $K'$ is isomorphic to
$Z\times Z.$ The first coordinate corresponds to the number of $op$ exponents on ${\cal A}$ and the
second component to the number of $op$ exponents on ${\cal B}.$ Negative integers correspond to the
right opposites, $po$. The power of the element of $H$ corresponds to the number of $op$ exponents
on the product of the two enriched categories, that is, the number of $op$ exponents outside the
parentheses. That $b\in H\sigma_2K'$ implies $Lb=Rb$ follows from the fact that the composition
morphisms belonging to the opposite of a ${\cal V}$-category obey the pentagon axiom. An exercise
of some value is to check consistency of the definitions by constructing an inductive proof of the
implication based on braid group generators. This is not a necessary condition for $Lb = Rb$, since
for example the equation holds for $b = (\sigma_2\sigma_1\sigma_3\sigma_2)^n,$ but it may be when
the additional requirement that $\sigma(b) = (2~3)$ is added. More work needs to be done to
determine the necessary conditions and to study the structure and properties of the braids that
meet these conditions. Of course we will see shortly that when certain unit conditions are obeyed
then there is a necessary and sufficient condition.
\end{remark}

  ${\cal V}$-functoriality of the associator is necessary because
  here we have a 2-natural transformation $\alpha^{(1)}$. This means we have a family of ${\cal V}$-functors
  indexed by triples of ${\cal V}$-categories. On objects

  $\alpha^{(1)}_{{\cal A}{\cal B}{\cal C}}((A,B),C) = (A,(B,C)).$

  In order to guarantee that $\alpha^{(1)}$ obey the coherence pentagon for hom-object morphisms,
  we define it to be {\it based upon}
  $\alpha$ in ${\cal V}.$ This means precisely that:
  \begin{footnotesize}
  $$\alpha^{(1)}_{{\cal A}{\cal B}{\cal C}_{((A,B),C)((A',B'),C')}}: [({\cal A}\otimes {\cal B})\otimes {\cal C}]\Bigl(\bigl((A,B),C\bigr),\bigl((A',B'),C'\bigr)\Bigr) \to [{\cal A}\otimes ({\cal B}\otimes {\cal C})]\Bigl(\bigl(A,(B,C)\bigr),\bigl(A',(B',C')\bigr)\Bigr)$$
  \end{footnotesize}
  is defined to be
  \begin{small}
  $$\alpha_{{\cal A}(A,A'){\cal B}(B,B'){\cal C}(C,C')}:\Bigl({\cal A}(A,A')\otimes {\cal B}(B,B')\Bigr)\otimes {\cal C}(C,C')\to {\cal A}(A,A') \otimes \Bigl({\cal B}(B,B')\otimes {\cal C}(C,C')\Bigr).$$
  \end{small}
  This definition guarantees that the $\alpha^{(1)}$ pentagons for objects and for hom-objects commute:
  the first trivially and the second by the fact that the
  $\alpha$ pentagon commutes in ${\cal V}.$
  We must also check for ${\cal V}$-functoriality. The unit axioms are trivial -- we consider the more interesting
  axiom. The following diagram must commute, where the first vertex is:
  $$[({\cal A}\otimes {\cal B})\otimes {\cal C}]\Bigl(\bigl((A',B'),C'\bigr),\bigl((A'',B''),C''\bigr)\Bigr)\otimes[({\cal A}\otimes {\cal B})\otimes {\cal C}]\Bigl(\bigl((A,B),C\bigr),\bigl((A',B'),C'\bigr)\Bigr)$$
  and the last vertex is:
  $$[{\cal A}\otimes ({\cal B}\otimes {\cal C})]\Bigl(\bigl(A,(B,C)\bigr),\bigl(A'',(B'',C'')\bigr)\Bigr).$$
  $$
  \xymatrix{
  \bullet
  \ar[rr]^{M}
  \ar[d]^{\alpha^{(1)} \otimes \alpha^{(1)}}
  &&\bullet
  \ar[d]^{\alpha^{(1)}}
  \\
  \bullet
  \ar[rr]_{M}
  &&\bullet
  }
  $$

   This means that the exterior of the following expanded diagram is required to commute.
    We leave out some parentheses
    for clarity and denote the middle four interchange by $\eta_{(b)}$ (perhaps composed with some associators).
    Also for convenience we write
   $X = {\cal A}(A',A'')$, $Y = {\cal B}(B',B'')$,
   $Z = {\cal C}(C',C'')$, $X' = {\cal A}(A,A')$, $Y' = {\cal B}(B,B')$ and $Z' = {\cal C}(C,C')$

  \noindent
                              \begin{center}
                  \resizebox{6.4in}{!}{
  $$
  \xymatrix{
  &(X\otimes Y)\otimes Z \otimes (X'\otimes Y')\otimes Z'
  \ar[dr]^{\eta_{(b)}}
  \ar[dl]_{\alpha}\\
  X\otimes (Y\otimes Z) \otimes X'\otimes (Y'\otimes Z')
  \ar[d]_{\eta_{(b)}}
  &&(X\otimes Y)\otimes (X'\otimes Y') \otimes Z \otimes Z'
  \ar[d]^\alpha\\
  X\otimes X' \otimes (Y\otimes Z)\otimes (Y'\otimes Z')
  \ar[d]_{\alpha}
  &&[X\otimes Y\otimes X'\otimes Y'] \otimes (Z\otimes Z')
  \ar[d]^{\eta_{(b)}}
  \\
  (X\otimes X') \otimes [Y\otimes Z\otimes Y'\otimes Z']
  \ar[d]_{\eta_{(b)}}
  &&[(X\otimes X')\otimes (Y\otimes Y')] \otimes (Z\otimes Z')
  \ar[dll]^{\alpha}
  \ar[d]^{(M\otimes M) \otimes M}\\
  (X\otimes X') \otimes [(Y\otimes Y')\otimes (Z\otimes Z')]
  \ar[dr]_{M\otimes (M \otimes M)~}
  &&[{\cal A}(A,A'')\otimes {\cal B}(B,B'')]\otimes {\cal C}(C,C'')
  \ar[dl]^{\alpha}\\
  &{\cal A}(A,A'')\otimes [{\cal B}(B,B'')\otimes {\cal C}(C,C'')]
  }
  $$
  }
  \end{center}

  The bottom quadrilateral commutes by naturality of $\alpha$.
   We are left needing to show that the underlying braids are equal
    for the two legs of the upper region.
     In Example~\ref{Lbprime} we give some examples of the underlying braids of the left and right
     legs
    for the same choices of $b$ as shown in Example~\ref{LbRb}. By inspection of the
    diagram these left and right underlying braids are the
    six-strand braids we denote respectively  $L'b$ and $R'b$.
    Recall from the introduction that  the requirement that $L'b=R'b$ is called external associativity.
  The first braid is the one used in the original definition of ${\cal A}\otimes {\cal B}$ given above.
  \begin{example}\label{Lbprime}
  ~
  \end{example}
   $$
   \xy (0,0) *{
   b_{(1)}=b_{0^+}=~} \endxy
   \xy (0,0) *{
   \xy 0;/r1pc/: +(-2,0)
   ,{\xcapv-@(0)}+(1,1),{\xoverv}+(2,1),{\xcapv-@(0)}
   \endxy } \endxy
   \xy (0,0) *{
   \text{ ; $~L'b_{(1)}=~$} }\endxy
   \xy (0,0) *{
   \xy 0;/r1pc/: +(-3,0)
      ,{\xcapv-@(0)}+(1,1),{\xcapv-@(0)}+(1,1),{\xoverv}+(2,1),{\xcapv-@(0)}+(1,1),{\xcapv-@(0)}+(-5,0)
      ,{\xcapv-@(0)}+(1,1),{\xoverv}+(2,1),{\xcapv-@(0)}+(1,1),{\xcapv-@(0)}+(1,1),{\xcapv-@(0)}+(-5,0)
      ,{\xcapv-@(0)}+(1,1),{\xcapv-@(0)}+(1,1),{\xcapv-@(0)}+(1,1),{\xoverv}+(2,1),{\xcapv-@(0)}
   \endxy }\endxy
   \xy (0,0) *{~=~} \endxy
   \xy (0,0) *{
   \xy 0;/r1pc/: +(-3,0)
      ,{\xcapv-@(0)}+(1,1),{\xcapv-@(0)}+(1,1),{\xoverv}+(2,1),{\xcapv-@(0)}+(1,1),{\xcapv-@(0)}+(-5,0)
      ,{\xcapv-@(0)}+(1,1),{\xcapv-@(0)}+(1,1),{\xcapv-@(0)}+(1,1),{\xoverv}+(2,1),{\xcapv-@(0)}+(-5,0)
      ,{\xcapv-@(0)}+(1,1),{\xoverv}+(2,1),{\xcapv-@(0)}+(1,1),{\xcapv-@(0)}+(1,1),{\xcapv-@(0)}
   \endxy } \endxy
   \xy (0,0) *{~=R'b_{(1)}} \endxy
   $$

   $$
   \xy (0,0) *{
   b_{(2)}=~}\endxy
   \xy (0,0) *{
   \xy 0;/r1pc/:
   ,{\xcapv-@(0)}+(1,1),{\xcapv-@(0)}+(1,1),{\xoverv}+(-2,0)
   ,{\xcapv-@(0)}+(1,1),{\xcapv-@(0)}+(1,1),{\xoverv}+(-2,0)
   ,{\xcapv-@(0)}+(1,1),{\xoverv}+(2,1),{\xcapv-@(0)}
   \endxy }\endxy
   \xy (0,0) *{
   \text{ ;  $~L'b_{(2)}=~$}}\endxy
   \xy (0,0) *{
   \xy 0;/r1pc/:
      ,{\xcapv-@(0)}+(1,1),{\xcapv-@(0)}+(1,1),{\xcapv-@(0)}+(1,1),{\xoverv}+(2,1),{\xcapv-@(0)}+(-5,0)
      ,{\xcapv-@(0)}+(1,1),{\xcapv-@(0)}+(1,1),{\xcapv-@(0)}+(1,1),{\xcapv-@(0)}+(1,1),{\xoverv}+(-4,0)
      ,{\xcapv-@(0)}+(1,1),{\xcapv-@(0)}+(1,1),{\xcapv-@(0)}+(1,1),{\xcapv-@(0)}+(1,1),{\xoverv}+(-4,0)
      ,{\xcapv-@(0)}+(1,1),{\xcapv-@(0)}+(1,1),{\xcapv-@(0)}+(1,1),{\xoverv}+(2,1),{\xcapv-@(0)}+(-5,0)
      ,{\xcapv-@(0)}+(1,1),{\xcapv-@(0)}+(1,1),{\xoverv}+(2,1),{\xcapv-@(0)}+(1,1),{\xcapv-@(0)}+(-5,0)
      ,{\xcapv-@(0)}+(1,1),{\xoverv}+(2,1),{\xcapv-@(0)}+(1,1),{\xcapv-@(0)}+(1,1),{\xcapv-@(0)}+(-5,0)
      ,{\xcapv-@(0)}+(1,1),{\xcapv-@(0)}+(1,1),{\xcapv-@(0)}+(1,1),{\xcapv-@(0)}+(1,1),{\xoverv}+(-4,0)
      ,{\xcapv-@(0)}+(1,1),{\xcapv-@(0)}+(1,1),{\xcapv-@(0)}+(1,1),{\xcapv-@(0)}+(1,1),{\xoverv}+(-4,0)
      ,{\xcapv-@(0)}+(1,1),{\xcapv-@(0)}+(1,1),{\xcapv-@(0)}+(1,1),{\xoverv}+(2,1),{\xcapv-@(0)}+(-5,0)
   \endxy }\endxy
   \xy (0,0) *{~=~}\endxy
   \xy (0,0) *{
   \xy 0;/r1pc/:
      ,{\xcapv-@(0)}+(1,1),{\xcapv-@(0)}+(1,1),{\xcapv-@(0)}+(1,1),{\xcapv-@(0)}+(1,1),{\xoverv}+(-4,0)
      ,{\xcapv-@(0)}+(1,1),{\xcapv-@(0)}+(1,1),{\xcapv-@(0)}+(1,1),{\xoverv}+(2,1),{\xcapv-@(0)}+(-5,0)
      ,{\xcapv-@(0)}+(1,1),{\xcapv-@(0)}+(1,1),{\xcapv-@(0)}+(1,1),{\xoverv}+(2,1),{\xcapv-@(0)}+(-5,0)
      ,{\xcapv-@(0)}+(1,1),{\xcapv-@(0)}+(1,1),{\xoverv}+(2,1),{\xoverv}+(-4,0)
      ,{\xcapv-@(0)}+(1,1),{\xcapv-@(0)}+(1,1),{\xcapv-@(0)}+(1,1),{\xoverv}+(2,1),{\xcapv-@(0)}+(-5,0)
      ,{\xcapv-@(0)}+(1,1),{\xcapv-@(0)}+(1,1),{\xoverv}+(2,1),{\xcapv-@(0)}+(1,1),{\xcapv-@(0)}+(-5,0)
      ,{\xcapv-@(0)}+(1,1),{\xcapv-@(0)}+(1,1),{\xoverv}+(2,1),{\xcapv-@(0)}+(1,1),{\xcapv-@(0)}+(-5,0)
      ,{\xcapv-@(0)}+(1,1),{\xoverv}+(2,1),{\xcapv-@(0)}+(1,1),{\xcapv-@(0)}+(1,1),{\xcapv-@(0)}
   \endxy}\endxy
   \xy (0,0) *{~=R'b_{(2)}}\endxy
   $$
   $$
   \xy (0,0) *{
   b_{(3)}=~}\endxy
   \xy (0,0) *{
   \xy 0;/r1pc/:
   ,{\xcapv-@(0)}+(1,1),{\xoverv}+(2,1),{\xcapv-@(0)}+(-3,0)
   ,{\xoverv}+(2,1),{\xcapv-@(0)}+(1,1),{\xcapv-@(0)}+(-3,0)
   ,{\xoverv}+(2,1),{\xcapv-@(0)}+(1,1),{\xcapv-@(0)}
   \endxy}\endxy
   \xy (0,0) *{
   \text{ ;  $~L'b_{(3)}=$~}}\endxy
   \xy (0,0) *{
   \xy 0;/r1pc/:
      ,{\xcapv-@(0)}+(1,1),{\xcapv-@(0)}+(1,1),{\xoverv}+(2,1),{\xcapv-@(0)}+(1,1),{\xcapv-@(0)}+(-5,0)
      ,{\xcapv-@(0)}+(1,1),{\xoverv}+(2,1),{\xcapv-@(0)}+(1,1),{\xcapv-@(0)}+(1,1),{\xcapv-@(0)}+(-5,0)
      ,{\xoverv}+(2,1),{\xcapv-@(0)}+(1,1),{\xcapv-@(0)}+(1,1),{\xcapv-@(0)}+(1,1),{\xcapv-@(0)}+(-5,0)
      ,{\xoverv}+(2,1),{\xcapv-@(0)}+(1,1),{\xcapv-@(0)}+(1,1),{\xcapv-@(0)}+(1,1),{\xcapv-@(0)}+(-5,0)
      ,{\xcapv-@(0)}+(1,1),{\xcapv-@(0)}+(1,1),{\xcapv-@(0)}+(1,1),{\xoverv}+(2,1),{\xcapv-@(0)}+(-5,0)
      ,{\xcapv-@(0)}+(1,1),{\xcapv-@(0)}+(1,1),{\xoverv}+(2,1),{\xcapv-@(0)}+(1,1),{\xcapv-@(0)}+(-5,0)
      ,{\xcapv-@(0)}+(1,1),{\xcapv-@(0)}+(1,1),{\xoverv}+(2,1),{\xcapv-@(0)}+(1,1),{\xcapv-@(0)}+(-5,0)
   \endxy}\endxy
   \xy (0,0) *{~\ne~}\endxy
   \xy (0,0) *{
   \xy 0;/r1pc/:
      ,{\xcapv-@(0)}+(1,1),{\xcapv-@(0)}+(1,1),{\xoverv}+(2,1),{\xcapv-@(0)}+(1,1),{\xcapv-@(0)}+(-5,0)
      ,{\xcapv-@(0)}+(1,1),{\xoverv}+(2,1),{\xoverv}+(2,1),{\xcapv-@(0)}+(-5,0)
      ,{\xoverv}+(2,1),{\xoverv}+(2,1),{\xcapv-@(0)}+(1,1),{\xcapv-@(0)}+(-5,0)
      ,{\xcapv-@(0)}+(1,1),{\xoverv}+(2,1),{\xcapv-@(0)}+(1,1),{\xcapv-@(0)}+(1,1),{\xcapv-@(0)}+(-5,0)
      ,{\xcapv-@(0)}+(1,1),{\xoverv}+(2,1),{\xcapv-@(0)}+(1,1),{\xcapv-@(0)}+(1,1),{\xcapv-@(0)}+(-5,0)
      ,{\xoverv}+(2,1),{\xoverv}+(2,1),{\xcapv-@(0)}+(1,1),{\xcapv-@(0)}+(-5,0)
      ,{\xcapv-@(0)}+(1,1),{\xoverv}+(2,1),{\xcapv-@(0)}+(1,1),{\xcapv-@(0)}+(1,1),{\xcapv-@(0)}+(-5,0)
      ,{\xcapv-@(0)}+(1,1),{\xoverv}+(2,1),{\xcapv-@(0)}+(1,1),{\xcapv-@(0)}+(1,1),{\xcapv-@(0)}+(-5,0)
      ,{\xoverv}+(2,1),{\xcapv-@(0)}+(1,1),{\xcapv-@(0)}+(1,1),{\xcapv-@(0)}+(1,1),{\xcapv-@(0)}+(-5,0)
      ,{\xoverv}+(2,1),{\xcapv-@(0)}+(1,1),{\xcapv-@(0)}+(1,1),{\xcapv-@(0)}+(1,1),{\xcapv-@(0)}+(-5,0)
   \endxy}\endxy
   \xy (0,0) *{~=R'b_{(3)}}\endxy
   $$

   $$
   \xy (0,0) *{
   b_{(4)}=~}\endxy
   \xy (0,0) *{
   \xy 0;/r1pc/:
   ,{\xcapv-@(0)}+(1,1),{\xoverv}+(2,1),{\xcapv-@(0)}+(-3,0)
   ,{\xcapv-@(0)}+(1,1),{\xoverv}+(2,1),{\xcapv-@(0)}+(-3,0)
   ,{\xcapv-@(0)}+(1,1),{\xoverv}+(2,1),{\xcapv-@(0)}+(-3,0)
   \endxy}\endxy
   \xy (0,0) *{
   \text{ ; $~L'b_{(4)}=~$}}\endxy
   \xy (0,0) *{
   \xy 0;/r1pc/:
      ,{\xcapv-@(0)}+(1,1),{\xcapv-@(0)}+(1,1),{\xoverv}+(2,1),{\xcapv-@(0)}+(1,1),{\xcapv-@(0)}+(-5,0)
      ,{\xcapv-@(0)}+(1,1),{\xoverv}+(2,1),{\xcapv-@(0)}+(1,1),{\xcapv-@(0)}+(1,1),{\xcapv-@(0)}+(-5,0)
      ,{\xcapv-@(0)}+(1,1),{\xoverv}+(2,1),{\xcapv-@(0)}+(1,1),{\xcapv-@(0)}+(1,1),{\xcapv-@(0)}+(-5,0)
      ,{\xcapv-@(0)}+(1,1),{\xcapv-@(0)}+(1,1),{\xoverv}+(2,1),{\xcapv-@(0)}+(1,1),{\xcapv-@(0)}+(-5,0)
      ,{\xcapv-@(0)}+(1,1),{\xcapv-@(0)}+(1,1),{\xoverv}+(2,1),{\xcapv-@(0)}+(1,1),{\xcapv-@(0)}+(-5,0)
      ,{\xcapv-@(0)}+(1,1),{\xoverv}+(2,1),{\xcapv-@(0)}+(1,1),{\xcapv-@(0)}+(1,1),{\xcapv-@(0)}+(-5,0)
      ,{\xcapv-@(0)}+(1,1),{\xcapv-@(0)}+(1,1),{\xcapv-@(0)}+(1,1),{\xoverv}+(2,1),{\xcapv-@(0)}+(-5,0)
      ,{\xcapv-@(0)}+(1,1),{\xcapv-@(0)}+(1,1),{\xcapv-@(0)}+(1,1),{\xoverv}+(2,1),{\xcapv-@(0)}+(-5,0)
      ,{\xcapv-@(0)}+(1,1),{\xcapv-@(0)}+(1,1),{\xcapv-@(0)}+(1,1),{\xoverv}+(2,1),{\xcapv-@(0)}+(-5,0)
   \endxy}\endxy
   \xy (0,0) *{~\ne~}\endxy
   \xy (0,0) *{
   \xy 0;/r1pc/:
      ,{\xcapv-@(0)}+(1,1),{\xcapv-@(0)}+(1,1),{\xoverv}+(2,1),{\xcapv-@(0)}+(1,1),{\xcapv-@(0)}+(-5,0)
      ,{\xcapv-@(0)}+(1,1),{\xcapv-@(0)}+(1,1),{\xcapv-@(0)}+(1,1),{\xoverv}+(2,1),{\xcapv-@(0)}+(-5,0)
      ,{\xcapv-@(0)}+(1,1),{\xcapv-@(0)}+(1,1),{\xcapv-@(0)}+(1,1),{\xoverv}+(2,1),{\xcapv-@(0)}+(-5,0)
      ,{\xcapv-@(0)}+(1,1),{\xcapv-@(0)}+(1,1),{\xoverv}+(2,1),{\xcapv-@(0)}+(1,1),{\xcapv-@(0)}+(-5,0)
      ,{\xcapv-@(0)}+(1,1),{\xcapv-@(0)}+(1,1),{\xoverv}+(2,1),{\xcapv-@(0)}+(1,1),{\xcapv-@(0)}+(-5,0)
      ,{\xcapv-@(0)}+(1,1),{\xcapv-@(0)}+(1,1),{\xcapv-@(0)}+(1,1),{\xoverv}+(2,1),{\xcapv-@(0)}+(-5,0)
      ,{\xcapv-@(0)}+(1,1),{\xoverv}+(2,1),{\xcapv-@(0)}+(1,1),{\xcapv-@(0)}+(1,1),{\xcapv-@(0)}+(-5,0)
      ,{\xcapv-@(0)}+(1,1),{\xoverv}+(2,1),{\xcapv-@(0)}+(1,1),{\xcapv-@(0)}+(1,1),{\xcapv-@(0)}+(-5,0)
      ,{\xcapv-@(0)}+(1,1),{\xoverv}+(2,1),{\xcapv-@(0)}+(1,1),{\xcapv-@(0)}+(1,1),{\xcapv-@(0)}+(-5,0)
   \endxy}\endxy
   \xy (0,0) *{~=R'b_{(4)}}\endxy
   $$
   $$
   \xy (0,0) *{
      b_{(5)}=b_{1^-}=~}\endxy
      \begin{turn}{180}
      \xy (0,0) *{
      \xy 0;/r1pc/:
      ,{\xcapv-@(0)}+(1,1),{\xunderv}+(2,1),{\xcapv-@(0)}+(-3,0)
      ,{\xunderv}+(2,1),{\xunderv}+(-2,0)
      ,{\xcapv-@(0)}+(1,1),{\xunderv}+(2,1),{\xcapv-@(0)}+(-3,0)
      ,{\xcapv-@(0)}+(1,1),{\xunderv}+(2,1),{\xcapv-@(0)}+(-3,0)
      ,{\xoverv}+(2,1),{\xoverv}+(-3,0)
      \endxy}\endxy
      \end{turn}
      \xy (0,0) *{
      \text{ ; $~L'b_{(5)}=~$}}\endxy
      \begin{turn}{180}
      \xy (0,0) *{
      \xy 0;/r1pc/:
         ,{\xcapv-@(0)}+(1,1),{\xunderv}+(2,1),{\xcapv-@(0)}+(1,1),{\xcapv-@(0)}+(1,1),{\xcapv-@(0)}+(-5,0)
      ,{\xunderv}+(2,1),{\xunderv}+(2,1),{\xcapv-@(0)}+(1,1),{\xcapv-@(0)}+(-5,0)
      ,{\xcapv-@(0)}+(1,1),{\xunderv}+(2,1),{\xcapv-@(0)}+(1,1),{\xcapv-@(0)}+(1,1),{\xcapv-@(0)}+(-5,0)
      ,{\xcapv-@(0)}+(1,1),{\xunderv}+(2,1),{\xcapv-@(0)}+(1,1),{\xcapv-@(0)}+(1,1),{\xcapv-@(0)}+(-5,0)
      ,{\xoverv}+(2,1),{\xoverv}+(2,1),{\xcapv-@(0)}+(1,1),{\xcapv-@(0)}+(-5,0)
      ,{\xcapv-@(0)}+(1,1),{\xcapv-@(0)}+(1,1),{\xcapv-@(0)}+(1,1),{\xunderv}+(2,1)                ,{\xcapv-@(0)}+(-5,0)
      ,{\xcapv-@(0)}+(1,1),{\xcapv-@(0)}+(1,1),{\xunderv}+(2,1)                  ,{\xunderv}+(-4,0)
      ,{\xcapv-@(0)}+(1,1),{\xunderv}+(2,1)                   ,{\xunderv}+(2,1)                     ,{\xcapv-@(0)}+(-5,0)
      ,{\xunderv}+(2,1)                    ,{\xunderv}+(2,1)                    ,{\xcapv-@(0)}+(1,1),{\xcapv-@(0)}+(-5,0)
      ,{\xcapv-@(0)}+(1,1),{\xunderv}+(2,1)                   ,{\xcapv-@(0)}+(1,1),{\xcapv-@(0)}+(1,1),{\xcapv-@(0)}+(-5,0)
      ,{\xcapv-@(0)}+(1,1),{\xunderv}+(2,1)                   ,{\xcapv-@(0)}+(1,1),{\xcapv-@(0)}+(1,1),{\xcapv-@(0)}+(-5,0)
      ,{\xoverv}+(2,1)                    ,{\xunderv}+(2,1)                    ,{\xcapv-@(0)}+(1,1),{\xcapv-@(0)}+(-5,0)
      ,{\xcapv-@(0)}+(1,1),{\xoverv}+(2,1)                      ,{\xoverv}+(2,1)                      ,{\xcapv-@(0)}+(-5,0)
      ,{\xcapv-@(0)}+(1,1),{\xcapv-@(0)}+(1,1),{\xcapv-@(0)}+(1,1),{\xcapv-@(0)}+(1,1),{\xoverv}+(-4,0)
      \endxy}\endxy
      \end{turn}
      \xy (0,0) *{~=~}\endxy
      \begin{turn}{180}
      \xy (0,0) *{
      \xy 0;/r1pc/:
      ,{\xcapv-@(0)}+(1,1),{\xcapv-@(0)}+(1,1),{\xcapv-@(0)}+(1,1),{\xunderv}+(2,1),{\xcapv-@(0)}+(-5,0)
       ,{\xcapv-@(0)}+(1,1),{\xcapv-@(0)}+(1,1),{\xunderv}+(2,1),{\xunderv}+(-4,0)
       ,{\xcapv-@(0)}+(1,1),{\xcapv-@(0)}+(1,1),{\xcapv-@(0)}+(1,1),{\xunderv}+(2,1),{\xcapv-@(0)}+(-5,0)
       ,{\xcapv-@(0)}+(1,1),{\xcapv-@(0)}+(1,1),{\xcapv-@(0)}+(1,1),{\xunderv}+(2,1),{\xcapv-@(0)}+(-5,0)
      ,{\xcapv-@(0)}+(1,1),{\xcapv-@(0)}+(1,1),{\xoverv}+(2,1),{\xoverv}+(-4,0)
         ,{\xcapv-@(0)}+(1,1),{\xunderv}+(2,1)           ,{\xcapv-@(0)}+(1,1),{\xcapv-@(0)}+(1,1),{\xcapv-@(0)}+(-5,0)
         ,{\xunderv}+(2,1)                    ,{\xunderv}+(2,1)                     ,{\xcapv-@(0)}+(1,1),{\xcapv-@(0)}+(-5,0)
         ,{\xcapv-@(0)}+(1,1),{\xunderv}+(2,1)           ,{\xunderv}+(2,1)                     ,{\xcapv-@(0)}+(-5,0)
         ,{\xcapv-@(0)}+(1,1),{\xcapv-@(0)}+(1,1),{\xunderv}+(2,1)                    ,{\xunderv}+(-4,0)
         ,{\xcapv-@(0)}+(1,1),{\xcapv-@(0)}+(1,1),{\xcapv-@(0)}+(1,1),{\xunderv}+(2,1)                   ,{\xcapv-@(0)}+(-5,0)
         ,{\xcapv-@(0)}+(1,1),{\xcapv-@(0)}+(1,1),{\xcapv-@(0)}+(1,1),{\xunderv}+(2,1)                   ,{\xcapv-@(0)}+(-5,0)
         ,{\xcapv-@(0)}+(1,1),{\xcapv-@(0)}+(1,1),{\xunderv}+(2,1)                    ,{\xoverv}+(-4,0)
         ,{\xcapv-@(0)}+(1,1),{\xoverv}+(2,1)           ,{\xoverv}+(2,1)                        ,{\xcapv-@(0)}+(-5,0)
         ,{\xoverv}+(2,1)                       ,{\xcapv-@(0)}+(1,1),{\xcapv-@(0)}+(1,1),{\xcapv-@(0)}+(1,1),{\xcapv-@(0)}+(-5,0)
          \endxy}\endxy
          \end{turn}
          \xy (0,0) *{~=R'b_{(5)}}\endxy
   $$

\begin{remark}
  A comparison with Example~\ref{LbRb} is of interest. Braids $b_{(2)}$ and $b_{(3)}$ are
  180 degree rotations of each other.
  Notice that the second braid in  Example~\ref{Lbprime}
  leads to an equality that is actually the same as for the third
  braid in Example~\ref{LbRb}. To see this the
  page can be rotated by 180 degrees. Similarly, the inequality preventing
  braid $b_{(2)}$ from yielding an associative composition morphism is the 180 degree
  rotation of the inequality preventing braid $b_{(3)}$ from yielding a ${\cal V}$-functorial associator.
  Braid $b_{(1)}$ and braid $b_{(5)}$ are each their own 180 degree rotation (we took advantage of
  the latter fact in drawing $L'b_{(5)}$ and $R'b_{(5)}$ above), and the two
  braids proving each to be the underlying braid of an associative composition
  morphism are the same two that show each to yield a ${\cal V}$-functorial associator.
  Braid $b_{(4)}$ is its own 180 degree rotation, and the two
  braids preventing it from being associative are the same two that obstruct it from being functorial.
  Thus there is a certain kind of duality
  between the requirements of associativity of the enriched composition and the
  functoriality of the associator.
  The full meaning of this duality becomes more clear in the study
  of (enrichment over) iterated monoidal categories as in \cite{forcey1}, where we see that in a braided category two
  potentially different tensor products have collapsed into one.

  If we were considering a strictly associative monoidal category ${\cal V}$ then
    the condition of a ${\cal V}$-functorial associator
    would become a condition of a well defined composition morphism.
\end{remark}

  The unit axioms required of the tensor product of two enriched categories
    are satisfied in general only if dropping either the first two or the last two strands of
    $b$
  leaves again the identity on two strands. (This is also
  due to the naturality of compositions of $\alpha$ and $c$ and the unit axioms obeyed by
  ${\cal A}$ and ${\cal B}.$) Recall that we refer to this as the internal unit condition.

  The canonical choice for the unit in ${\cal V}$-Cat is the enriched category ${\cal I}$
  which has only one object denoted $0$ and for which ${\cal I}(0,0)=I$, the
  unit in ${\cal V}$.
   For the unit ${\cal V}$-category ${\cal I}$ to be indeed a
     unit for the tensor product in question
    requires that in the underlying braid of the middle four interchange
    dropping either the first and third strand or the second and fourth
    strand leaves the identity on two strands. Recall that we refer to this as the external unit condition.
     For a careful demonstration of this see \cite{forcey1}, keeping in
  mind that the interchange $\eta$ described there corresponds to the middle four interchange here.
Note that the unit conditions are not met by $b_{(2)}$ and $b_{(3)}$ in the above examples.

  \begin{definition}
    An \emph{interchanging} \emph{unital} braid on four strands is one for which the permutation associated to the
    braid is $(2~3)$, for which both $Lb = Rb$  and $L'b = R'b$ in $B_6,$ and for which the
 unit conditions are satisfied:
    deleting any one of the pairs of strands
    (1, 2); (3, 4); (1, 3) , or (2, 4)
    results in the 2 strand identity braid. An \emph{interchange candidate} braid is an element of $B_4$ which has the
correct permutation and obeys the unit conditions.
          \end{definition}
Note that this definition describes precisely what needs to be true of a braid $b$ in order that it
 can arise as the underlying braid of the middle four interchange of the composition morphism of
the tensor product of enriched categories over an arbitrary braided category ${\cal V}.$ Here we
are restricting our attention to monoidal structures on  ${\cal V}$-Cat with the canonical choices
described in Definition~\ref{canontimes} for the objects, hom-objects and unit morphisms of a
tensor product of enriched categories. Also let the associator for that monoidal structure be based
upon the associator in the braided category ${\cal V}$, and the unit for that monoidal structure be
the canonical choice of the enriched category ${\cal I}.$

\begin{lemma}\label{underlies}
Given an arbitrary braided category ${\cal V},$ let a monoidal structure on  ${\cal V}$-Cat be
assumed to have the canonical choices described in Definition~\ref{canontimes} for the objects,
hom-objects and unit morphisms of a tensor product of enriched categories. Also let the associator
for that monoidal structure be based upon the associator in the braided category ${\cal V}$, and
the unit for that monoidal structure be the canonical choice of the enriched category ${\cal I}.$
Then a four-strand braid $b$ satisfies $L'b=R'b$ if it can arise as the underlying braid of the
middle four interchange of the composition morphism of the tensor product of enriched categories in
any such monoidal structure.
\end{lemma}
\begin{proof}
 Simply having a valid tensor product in ${\cal V}$-Cat for a  \emph{specific} braided
category ${\cal V}$, with a middle four interchange built out of instances of the the associators
and the braiding, does not imply that
 the underlying braid of the middle four interchange
is interchanging (e.g. the braiding might be a symmetry or the composition might be a coequalizer).
However if $b$ underlies a middle four interchange which gives a tensor product of ${\cal
V}$-categories which is valid for an arbitrary braided base, then we can get our result by choosing
the example of the free braided category on one object with duals, denoted $C_{1,2}$ as in
\cite{Baez}.

Notice that we need more structure than just the free braided category on one object. This is
because we are only given the equality in ${\cal V}$ implied by the diagram for ${\cal
V}$-functoriality of $\alpha^{(1)}$; this equality is that of two compositions of braidings each
with a tensor product of instances of $M$ attached.  Specifically $M\otimes M \otimes M$ follows
each of the compositions of braidings. Let the generating objects of $C_{1,2}$ be $x$ and its dual
$x^*.$ Recall that the objects are then strings of these generators and the morphisms are tangles
with the number of inputs the length of the domain and the number of outputs the length of the
range. The braiding is the same as described in \cite{JS} for the free braided category; the tangle
formed by crossing all the strands corresponding to an object $A$  (one strand for each generator
in the string) simultaneously with all those of $B$. We can find enriched categories over $C_{1,2}$
since there are monoids in $C_{1,2}$, recalling that monoids are one-object enriched categories. To
see the braid equality $L'b=R'b$ we can choose the monoid ${\cal X}=x\otimes x^*.$ Then the
composition morphism is given by:
$$M = 1_x\otimes e \otimes 1_{x^*} : x\otimes x^*\otimes x\otimes x^* \to x\otimes x^* $$
(where $e$ is the counit) which corresponds to the tangle:
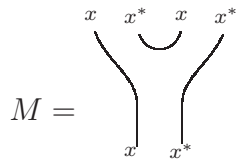
\[M = \xy
(-3,12)*{^{x^*}}="1"; 
(3,12)*{^x}="2"; 
(-9,12)*{^x}="A2";
(9,12)*{^{x^*}}="B2"; 
"1";"2" **\crv{(-3,7) & (3,7)};
(-3,-6)*{^{x^{~}}}="A";
(3,-6)*{^{x^*}}="B"; 
(-3,1)*{}="A1";
(3,1)*{}="B1"; 
"A";"A1" **\dir{-};
"B";"B1" **\dir{-}; 
"B2";"B1" **\crv{(8,7) & (3,5)}; "A2";"A1" **\crv{(-8,7) & (-3,5)};
\endxy \]
Now the equality implied by  ${\cal V}$-functoriality of the associator $\alpha^{(1)}_{{\cal
X}{\cal X}{\cal X}}$ is a tangle equality which is formed by starting with doubled versions of the
braids $L'b$ and $R'b$ (doubled since there is a strand for $x$ and $x^*$). Then both tangles are
finished with $M\otimes M \otimes M$, that is, three copies of the above tangle for $M$ attached.
For example, here is the left hand side (left leg) of the tangle equality for the braid $b_{(1)}$
that is implied by the ${\cal V}$-functoriality of the associator. Compare to $L'b_{(1)}$ above.

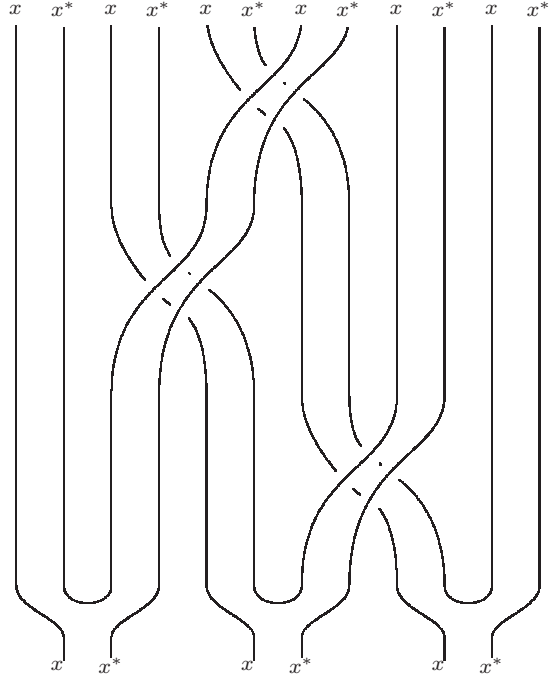
\[ \xy 0;/r.15pc/:
   (35,20)*{^{x^*}}="c3"; (25,20)*{^x}="c4";
  (55,20)*{^{x^*}}="c1"; (45,20)*{^x}="c2";
   (35,-20)*{}="Tc3"; (25,-20)*{}="Tc4";
  (55,-20)*{}="Tc1"; (45,-20)*{}="Tc2";
   (35,-60)*{}="TTc3"; (25,-60)*{}="TTc4";
  (55,-60)*{}="TTc1"; (45,-60)*{}="TTc2";
  (35,-100)*{}="TTTc3"; (25,-100)*{}="TTTc4";
  (55,-100)*{}="TTTc1"; (45,-100)*{}="TTTc2";
  (15,20)*{^{x^*}}="b1"; (5,20)*{^x}="b2";
  (-5,20)*{^{x^*}}="b3"; (-15,20)*{^x}="b4";
   (15,-20)*{}="T1"; (5,-20)*{}="T2";
  (-5,-20)*{}="T3"; (-15,-20)*{}="T4";
   (-5,-60)*{}="TT3"; (-15,-60)*{}="TT4";
  (15,-60)*{}="TT1"; (5,-60)*{}="TT2";
(-5,-100)*{}="TTT3"; (-15,-100)*{}="TTT4";
  (15,-100)*{}="TTT1"; (5,-100)*{}="TTT2";
   (-35,20)*{^x}="a1"; (-25,20)*{^{x^*}}="a2";
  (65,20)*{^x}="d1"; (75,20)*{^{x^*}}="d2";
    (-35,-100)*{}="Ta1"; (-25,-100)*{}="Ta2";
    (65,-100)*{}="Td1"; (75,-100)*{}="Td2";
"c4"; "T2" **\crv{(25,4) & (5,5)}
  \POS?(.35)*{\hole}="2cx" \POS?(.54)*{\hole}="2cy";
  "b2"; "2cy" **\crv{(5,12) };
  "b1"; "2cx" **\crv{(15,10) };
  "c3"; "T1" **\crv{(35,7) & (15,7)}
  \POS?(.55)*{\hole}="4cx" \POS?(.7)*{\hole}="4cy";
    "4cy";"2cy" **\crv{(18,-2)};
    "2cx";"4cx" **\crv{(23.5,3.5)};
  "4cy";"Tc4" **\crv{(25,-6)};
    "Tc3";"4cx" **\crv{(35,-4)};
"T3";"b3" **\dir{-};
 "T4";"b4" **\dir{-};
"Tc1";"c1" **\dir{-};
 "Tc2";"c2" **\dir{-};
  "T2"; "TT4" **\crv{(5,-36) & (-15,-35)}
  \POS?(.35)*{\hole}="2x" \POS?(.54)*{\hole}="2y";
  "T4"; "2y" **\crv{(-15,-28) };
  "T3"; "2x" **\crv{(-5,-30) };
  "T1"; "TT3" **\crv{(15,-33) & (-5,-33)}
  \POS?(.55)*{\hole}="4x" \POS?(.7)*{\hole}="4y";
    "4y";"2y" **\crv{(-2,-42)};
    "2x";"4x" **\crv{(3.5,-36.5)};
  "4y";"TT2" **\crv{(5,-46)};
    "TT1";"4x" **\crv{(15,-44)};
"Tc1";"TTc1" **\dir{-};
 "Tc2";"TTc2" **\dir{-};
"Tc3";"TTc3" **\dir{-};
 "Tc4";"TTc4" **\dir{-};
"TTc2"; "TTTc4" **\crv{(45,-76) & (25,-75)}
  \POS?(.35)*{\hole}="2cx" \POS?(.54)*{\hole}="2cy";
  "TTc4"; "2cy" **\crv{(25,-68) };
  "TTc3"; "2cx" **\crv{(35,-70) };
  "TTc1"; "TTTc3" **\crv{(55,-73) & (35,-73)}
  \POS?(.55)*{\hole}="4cx" \POS?(.7)*{\hole}="4cy";
    "4cy";"2cy" **\crv{(38,-82)};
    "2cx";"4cx" **\crv{(43.5,-76.5)};
  "4cy";"TTTc2" **\crv{(45,-86)};
    "TTTc1";"4cx" **\crv{(55,-84)};
    "TTT1";"TT1" **\dir{-};
 "TTT2";"TT2" **\dir{-};
"TTT3";"TT3" **\dir{-};
 "TTT4";"TT4" **\dir{-};
"a1";"Ta1" **\dir{-};
 "a2";"Ta2" **\dir{-};
"d1";"Td1" **\dir{-};
 "d2";"Td2" **\dir{-};
 (-25,-100)*{}="1"; %
(-15,-100)*{}="2"; %
(-35,-100)*{}="A2";%
(-5,-100)*{}="B2"; %
"1";"2" **\crv{(-25,-105) & (-15,-105)};%
(-25,-118)*{^{x^{~}}}="A";%
(-15,-118)*{^{x^*}}="B"; %
(-25,-111)*{}="A1";%
(-15,-111)*{}="B1"; %
"A";"A1" **\dir{-};
"B";"B1" **\dir{-}; 
"B2";"B1" **\crv{(-5,-105) & (-15,-107)}; "A2";"A1" **\crv{(-35,-105) & (-25,-107)};
 (25,-100)*{}="1"; %
(15,-100)*{}="2"; %
(35,-100)*{}="A2";%
(5,-100)*{}="B2"; %
"1";"2" **\crv{(25,-105) & (15,-105)};%
(25,-118)*{^{x^*}}="A";%
(15,-118)*{^{x^{~}}}="B"; %
(25,-111)*{}="A1";%
(15,-111)*{}="B1"; %
"A";"A1" **\dir{-};
"B";"B1" **\dir{-}; 
"B2";"B1" **\crv{(5,-105) & (15,-107)}; "A2";"A1" **\crv{(35,-105) & (25,-107)};
 (65,-100)*{}="1"; %
(55,-100)*{}="2"; %
(75,-100)*{}="A2";%
(45,-100)*{}="B2"; %
"1";"2" **\crv{(65,-105) & (55,-105)};%
(65,-118)*{^{x^*}}="A";%
(55,-118)*{^{x^{~}}}="B"; %
(65,-111)*{}="A1";%
(55,-111)*{}="B1"; %
"A";"A1" **\dir{-};
"B";"B1" **\dir{-}; 
"B2";"B1" **\crv{(45,-105) & (55,-107)}; "A2";"A1" **\crv{(75,-105) & (65,-107)};
\endxy \]

The right hand side is similarly drawn, with a doubled version of $R'b_{(1)}$ followed by three
copies of $M$.

That the two tangles are equal implies that their corresponding sub-tangles are equal; specifically
that their sub-tangles formed by deleting all but the input strands 1,3,5,8,10, and 12 (which
comprise all six output strands) are equal. These sub-tangles are $L'b$ and $R'b$ respectively.
\end{proof}

The question now is whether there are braids underlying the composition of a product of enriched
  categories besides the braids $b_{(1)}$ and $b_{(5)}$ above (and their inverses) which fulfill all obligations.
  The answer is yes.
  To find interchanging braids we iteratively build new monoidal structures from the standard ones,
  using the duality structure that exists on ${\cal V}$-Cat.
  By ${\cal A}^{op^n}$ we denote the $n^{th}$ (left) opposite of ${\cal A}$.
  By $\otimes$ and $\otimes'$ we denote the standard tensor products
  defined respectively with braid $b_{(1)}$ and its inverse underlying the middle four interchange.
  \begin{theorem}\label{newtimes}
  The tensor product of enriched categories given by
  $${\cal A}\otimes_{1^-}{\cal B} = ({\cal A}^{op}\otimes'{\cal B}^{op})^{po}$$
  is a valid monoidal product on ${\cal V}$-Cat. Furthermore, so are the tensor products
  $${\cal A}\otimes_{n^-}{\cal B} = ({\cal A}^{op^n}\otimes'{\cal B}^{op^n})^{po^n}$$
  as well as those with underlying braids that are the inverses of these, denoted
  $${\cal A}\otimes_{n^+}{\cal B} = ({\cal A}^{po^n}\otimes{\cal B}^{po^n})^{op^n}.$$
  \end{theorem}
  \begin{proof} The first tensor product is mentioned alone since the middle four interchange in its composition morphism
   has the underlying braid
  shown above as braid $b_{(5)}.$ Thus we have already demonstrated its fitness as a monoidal product.
  However this can be more efficiently shown just by noting that the category given by
  the product is certainly a valid enriched category, and that for three operands
  we have an associator
  from the isomorphism given by the following:

  \begin{align*}
\alpha^{1^-}_{{\cal A}{\cal B}{\cal C}} &: ({\cal A}\otimes_{1^-}{\cal B})\otimes_{1^-}{\cal C}\\
&~= ((({\cal A}^{op}\otimes'{\cal B}^{op})^{po})^{op} \otimes' {\cal C}^{op})^{po}\\
&~=(({\cal A}^{op}\otimes'{\cal B}^{op})\otimes' {\cal C}^{op})^{po}\\
&~\cong ({\cal A}^{op}\otimes'({\cal B}^{op}\otimes' {\cal C}^{op}))^{po}\\
&~=({\cal A}^{op}\otimes'(({\cal B}^{op}\otimes' {\cal C}^{op})^{po})^{op} )^{po}\\
&~= {\cal A}\otimes_{1^-}({\cal B}\otimes_{1^-}{\cal C})
\end{align*}

  The associator  implicit in this isomorphism
is constructed by taking the right opposite of instances of the standard associator for $\otimes'$;
 $$\alpha^{1^-}_{{\cal A}{\cal B}{\cal C}}=(\alpha'^{(1)}_{{\cal A}^{op}{\cal B}^{op}{\cal
 C}^{op}})^{po}.$$
The standard associator for $\otimes'$ is identical to the one for $\otimes.$  Thus $\alpha^{1^-}$
is based upon
  $\alpha$ in ${\cal V}$, since the object sets of the domain and range are the usual cartesian products and
  since
\begin{small}
  $$\alpha^{1^-}_{{\cal A}{\cal B}{\cal C}_{((A,B),C)((A',B'),C')}} =
   \alpha_{{\cal A}(A,A'){\cal B}(B,B'){\cal C}(C,C')} . $$\end{small}

This new associator is guaranteed to have ${\cal V}$-functorial instances since they are the images
(under the right opposite) of ${\cal V}$-functors.

  Inductively this process can be repeated with all the left opposites and right opposites raised to the
  $n^{th}$ degree.
    Recall that the unit ${\cal V}$-category ${\cal I}$ has only one object $0$ and ${\cal I}(0,0)=I,$ the
  unit in ${\cal V}$. That  ${\cal I}$ is indeed a unit for the tensor products in question
  follows from the facts  that ${\cal I}^{op} =  {\cal I} =  {\cal I}^{po}$ which are in turn evident from
  facts
  $c_{IA} = c_{AI} = 1_A.$
Thus we have that, using any of the above
  tensor products including the standard ones $\otimes = \otimes_{0^{+}}$ and $\otimes' = \otimes_{0^{-}}$
  defined respectively with braid $b_{(1)}$ and its inverse, ${\cal V}$-Cat
  is a monoidal 2-category. \end{proof}

The braids underlying these new tensor products
 are not hard
to describe directly. Suppressing the associators, the instances of the braiding forming the middle
four interchange for $\otimes_{1^+}$  are as follows, using $X',Y',X,Y$ to stand for hom objects as
above:
$$(c_{XX'}\otimes c_{YY'})\circ (1_{X}\otimes c^{-1}_{YX'} \otimes 1_{Y'})\circ(c^{-1}_{(X'\otimes Y')(X\otimes Y)}),$$
 with underlying braid $b_{(5)}.$ Note that $b_{(5)}$  has also been denoted $b_{1^-}.$
  Another is the following braid which underlies $\otimes_{2^+}$:
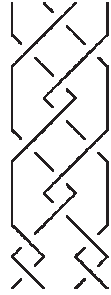
  $$
  \xy 0;/r1pc/:
        ,{\xcapv-@(0)}+(1,1),{\xoverv}+(2,1),{\xcapv-@(0)}+(-3,0)
        ,{\xoverv}+(2,1),{\xoverv}+(-2,0)
        ,{\xcapv-@(0)}+(1,1),{\xoverv}+(2,1),{\xcapv-@(0)}+(-3,0)
        ,{\xcapv-@(0)}+(1,1),{\xoverv}+(2,1),{\xcapv-@(0)}+(-3,0)
        ,{\xoverv}+(2,1),{\xoverv}+(-2,0)
        ,{\xcapv-@(0)}+(1,1),{\xoverv}+(2,1),{\xcapv-@(0)}+(-3,0)
        ,{\xcapv-@(0)}+(1,1),{\xoverv}+(2,1),{\xcapv-@(0)}+(-3,0)
        ,{\xunderv}+(2,1),{\xunderv}+(-2,0)
        ,{\xunderv}+(2,1),{\xunderv}+(-2,0)
      \endxy
  $$
Note that this is precisely the braid $b_{2^+}$ shown in the introduction. In fact the construction
of the new products leads to the observation that the braid underlying the middle four interchange
in the composition for the product ${\cal A}\otimes_{n^{\pm}}{\cal B}$  is the previously defined
braid $b_{n^{\pm}}.$
\begin{remark}\label{foop}
  Note that if in the definition of $\otimes_{n^{\pm}}$
  we replace $\otimes'$ with $\otimes$ or vice versa, then we have
  another valid tensor product, but with a braid underlying the middle four interchange in its composition
   morphisms equivalent to that found in $\otimes_{(n-1)^{\pm}}.$
   For example:
   $$
   \xy (0,0) *{
   \xy 0;/r1pc/:
      ,{\xcapv-@(0)}+(1,1),{\xunderv}+(2,1),{\xcapv-@(0)}+(-3,0)
      ,{\xunderv}+(2,1),{\xunderv}+(-2,0)
      ,{\xcapv-@(0)}+(1,1),{\xunderv}+(2,1),{\xcapv-@(0)}+(-3,0)
      ,{\xcapv-@(0)}+(1,1),{\xoverv}+(2,1),{\xcapv-@(0)}+(-3,0)
      ,{\xoverv}+(2,1),{\xoverv}+(-3,0)
      \endxy} \endxy
      \xy (0,0) *{~ =~} \endxy
      \xy (0,0) *{
      \xy 0;/r1pc/: +(-2,0)
   ,{\xcapv-@(0)}+(1,1),{\xunderv}+(2,1),{\xcapv-@(0)}
   \endxy} \endxy
   $$
  \end{remark}

Next we will show that this condition of being equivalent to some $b_{n^{\pm}}$
 is necessary for a braid to be interchanging, offer some quick checks to determine
when this condition holds, and  investigate when the resulting monoidal categories are equivalent.
All these steps are best taken in the context of iterated monoidal categories.



  \section{$2$-fold Monoidal Categories}\label{2fold}

    In this section we closely follow the authors of \cite{Balt} in defining a notion of iterated monoidal category.
    For those readers familiar with that source, note that we vary from their definition only
    by including associativity
    up to natural coherent isomorphisms. Thus we begin by reviewing the definition of lax monoidal
    functor. In our examples using a braided category, however, the natural transformations will
    all  be isomorphisms.


    \begin{definition}
    A {\it lax monoidal functor} $(F,\eta) :{\cal C}\to{\cal D}$ between monoidal categories consists
    of a functor $F$ such that $F(I)=I$ together with a natural transformation
    $$
    \eta_{AB}:F(A)\otimes F(B)\to F(A\otimes B),
    $$
    which satisfies the following conditions
    \begin{enumerate}
    \item Internal Associativity: The following diagram commutes
    $$
    \diagram
    (F(A)\otimes F(B))\otimes F(C)
    \rrto^{\eta_{AB}\otimes 1_{F(C)}}
    \dto^{\alpha}
    &&F(A\otimes B)\otimes F(C)
    \dto^{\eta_{(A\otimes B)C}}\\
    F(A)\otimes (F(B)\otimes F(C))
    \ar[d]^{1_{F(A)}\otimes \eta_{BC}}
    &&F((A\otimes B)\otimes C)
    \ar[d]^{F\alpha}\\
    F(A)\otimes F(B\otimes C)
    \rrto^{\eta_{A(B\otimes C)}}
    &&F(A\otimes (B\otimes C))
    \enddiagram
    $$
    \item Internal Unit Conditions: $\eta_{AI}=\eta_{IA}=1_{F(A)}.$
    \end{enumerate}
    \end{definition}
    Given two monoidal functors $(F,\eta) :{\cal C}\to{\cal D}$ and $(G,\zeta)
    :{\cal D}\to{\cal E}$,
    we define their composite to be the monoidal functor $(GF,\xi) :
    {\cal C}\to{\cal E}$, where
    $\xi$ denotes the composite
    $$
    \diagram
    GF(A)\otimes GF(B)\rrto^{\zeta_{F(A)F(B)}}
    && G\bigl(F(A)\otimes F(B)\bigr)\rrto^{G(\eta_{AB})}
    &&GF(A\otimes B).
    \enddiagram
    $$
    It is easy to verify that $\xi$ satisfies the internal associativity condition above by subdividing the
    necessary commuting diagram into two regions that commute by the axioms for $\eta$ and $\zeta$ respectively
    and two that commute due to their naturality.
     $\mathbf{MonCat}$ is the monoidal category of monoidal categories and monoidal
    functors, with the usual Cartesian product as in ${\mathbf{Cat}}$.

    A {\it monoidal natural transformation} $\theta:(F, \eta) \to (G, \zeta):{\cal D}\to{\cal E}$  is a
    natural transformation $\theta: F\to G$ between the underlying ordinary functors of $F$ and $G$ such that the
    following diagram commutes
    $$
    \xymatrix{
    F(A)\otimes F(B)
    \ar[r]^{\eta}
    \ar[d]^{\theta_A \otimes \theta_B}
    &F(A\otimes B)
    \ar[d]^{\theta_{A \otimes B}}
    \\
    G(A)\otimes G(B)
    \ar[r]^{\zeta}
    &G(A\otimes B)
    }
    $$



    \begin{definition} A  $2${\it -fold monoidal category} (with strong associators)
     is a monoidal category $({\cal V},\otimes_1,\alpha^1,I)$ and a
    monoidal functor
    $(\otimes_2,\eta):{\cal V}\times{\cal V}\to{\cal V}$ which satisfies
    \begin{enumerate}
    \item External Associativity: the following diagram describes a monoidal natural isomorphism
    $\alpha^2$ in $\mathbf{MonCat}.$

    $$
    \xymatrix{
    {\cal V}\times{\cal V}\times{\cal V}
    \rrto^{(\otimes_2,\eta)\times 1_{\cal V}}
    \dto_(0.4){1_{\cal V}\times(\otimes_2,\eta)}
    && {\cal V}\times{\cal V}
    \dto^{(\otimes_2,\eta)}
    \ar@{=>}[dll]^{\alpha^2}
    \\
    {\cal V}\times{\cal V}
    \rrto_{(\otimes_2,\eta)}
    &&
    {\cal V}
    }
    $$

    \item External Unit Conditions: the following diagram commutes in
    $\mathbf{MonCat}$

    $$
    \diagram
    {\cal V}\times I
    \rto^{\subseteq}
    \ddrto^{\cong}
    & {\cal V}\times{\cal V}
    \ddto^{(\otimes_2,\eta)}
    & I\times{\cal V}
    \lto_{\supseteq}
    \ddlto^{\cong}\\\\
    &{\cal V}
    \enddiagram
    $$
    \item Coherence: The underlying natural transformation $\alpha^2$ satisfies
    the usual coherence pentagon.

    \end{enumerate}
    \end{definition}

    Explicitly this means that we are given a second associative binary operation
    $\otimes_2:{\cal V}\times{\cal V}\to{\cal V}$, for which $I$ is also a two-sided unit.
   We are also given a natural transformation called the \emph{interchange} which is the
   functoriality constraint for $\otimes_2:$
    $$
    \eta_{ABCD}: (A\otimes_2 B)\otimes_1 (C\otimes_2 D)\to
    (A\otimes_1 C)\otimes_2(B\otimes_1 D).
    $$
    The internal unit conditions for $\otimes_2$ as a monoidal functor give $\eta_{ABII}=\eta_{IIAB}=1_{A\otimes_2 B}$,
    while the external unit conditions give $\eta_{AIBI}=\eta_{IAIB}=1_{A\otimes_1 B}$.
    The internal associativity condition for $\otimes_2$ as a monoidal functor gives the commutative
    diagram:
    \begin{footnotesize}
    $$
    \diagram
    ((U\otimes_2 V)\otimes_1 (W\otimes_2 X))\otimes_1 (Y\otimes_2 Z)
    \xto[rrr]^{\eta_{UVWX}\otimes_1 1_{Y\otimes_2 Z}}
    \ar[d]^{\alpha^1}
    &&&\bigl((U\otimes_1 W)\otimes_2(V\otimes_1 X)\bigr)\otimes_1 (Y\otimes_2 Z)
    \dto^{\eta_{(U\otimes_1 W)(V\otimes_1 X)YZ}}\\
    (U\otimes_2 V)\otimes_1 ((W\otimes_2 X)\otimes_1 (Y\otimes_2 Z))
    \dto^{1_{U\otimes_2 V}\otimes_1 \eta_{WXYZ}}
    &&&((U\otimes_1 W)\otimes_1 Y)\otimes_2((V\otimes_1 X)\otimes_1 Z)
    \ar[d]^{\alpha^1 \otimes_2 \alpha^1}
    \\
    (U\otimes_2 V)\otimes_1 \bigl((W\otimes_1 Y)\otimes_2(X\otimes_1 Z)\bigr)
    \xto[rrr]^{\eta_{UV(W\otimes_1 Y)(X\otimes_1 Z)}}
    &&& (U\otimes_1 (W\otimes_1 Y))\otimes_2(V\otimes_1 (X\otimes_1 Z))
    \enddiagram
    $$
    \end{footnotesize}
    The external associativity condition ($\alpha^2$ must be a monoidal natural transformation)
     gives the commutative diagram:
    \begin{footnotesize}
    $$
    \diagram
    ((U\otimes_2 V)\otimes_2 W)\otimes_1 ((X\otimes_2 Y)\otimes_2 Z)
    \xto[rrr]^{\eta_{(U\otimes_2 V)W(X\otimes_2 Y)Z}}
    \ar[d]^{\alpha^2 \otimes_1 \alpha^2}
    &&& \bigl((U\otimes_2 V)\otimes_1 (X\otimes_2 Y)\bigr)\otimes_2(W\otimes_1 Z)
    \dto^{\eta_{UVXY}\otimes_2 1_{W\otimes_1 Z}}\\
    (U\otimes_2 (V\otimes_2 W))\otimes_1 (X\otimes_2 (Y\otimes_2 Z))
    \dto^{\eta_{U(V\otimes_2 W)X(Y\otimes_2 Z)}}
    &&&((U\otimes_1 X)\otimes_2(V\otimes_1 Y))\otimes_2(W\otimes_1 Z)
    \ar[d]^{\alpha^2}
    \\
    (U\otimes_1 X)\otimes_2\bigl((V\otimes_2 W)\otimes_1 (Y\otimes_2 Z)\bigr)
    \xto[rrr]^{1_{U\otimes_1 X}\otimes_2\eta_{VWYZ}}
    &&& (U\otimes_1 X)\otimes_2((V\otimes_1 Y)\otimes_2(W\otimes_1 Z))
    \enddiagram
    $$
    \end{footnotesize}

Just as in \cite{Balt} we now define a 2-fold monoidal functor $(F,\lambda^1, \lambda^2)$
    between 2-fold monoidal categories. It is a functor $F:{\cal V}\to{\cal W}$ together
    with two natural transformations:
    $$\lambda^1_{AB}:F(A)\otimes_1 F(B)\to F(A\otimes_1 B)$$
    $$\lambda^2_{AB}:F(A)\otimes_2 F(B)\to F(A\otimes_2 B)$$
    satisfying the same associativity and unit conditions as in the case of monoidal functors.
    In addition we require that the following hexagonal interchange diagram commutes:
    $$
    \diagram
    (F(A)\otimes_2 F(B))\otimes_1(F(C)\otimes_2 F(D))
    \xto[rrr]^{ {\eta_{F(A)F(B)F(C)F(D)}}}
    \dto^{ {\lambda^2_{AB}\otimes_1\lambda^2_{CD}}}
    &&&(F(A)\otimes_1 F(C))\otimes_2(F(B)\otimes_1 F(D))
    \dto^{ \lambda^1_{AC}\otimes_2\lambda^1_{BD}}\\
    F(A\otimes_2 B)\otimes_1 F(C\otimes_2 D)
    \dto^{ {\lambda^1_{(A\otimes_2 B)(C\otimes_2 D)}}}
    &&&F(A\otimes_1 C)\otimes_2 F(B\otimes_1 D)
    \dto^{ {\lambda^2_{(A\otimes_1 C)(B\otimes_1 D)}}}\\
    F((A\otimes_2 B)\otimes_1(C\otimes_2 D))
    \xto[rrr]^{ F(\eta_{ABCD})}
    &&& F((A\otimes_1 C)\otimes_2(B\otimes_1 D))
    \enddiagram
    $$

    We can now refer to the category {\textbf{2}$\mathbf{ -MonCat}$ of 2-fold monoidal categories and
    2-fold monoidal functors.

    The authors of \cite{Balt} remark that we have natural transformations
    $$
    \eta_{AIIB}:A\otimes_1 B\to A\otimes_2 B\qquad\mbox{ and }\qquad
    \eta_{IABI}:A\otimes_1 B\to B\otimes_2 A.
    $$
    If they had insisted a 2-fold monoidal category be a tensor object in the category of monoidal categories
    and {\it strictly monoidal\/} functors, this would be equivalent to requiring that $\eta=1$.  In view
    of the above, they note that this would imply $A\otimes_1 B = A\otimes_2 B = B\otimes_1 A$ and similarly for morphisms.
    This is shown by what is usually referred to as the Eckmann-Hilton argument.

     Joyal and Street \cite{JS} considered a
    similar concept to Balteanu, Fiedorowicz, Schw${\rm \ddot a}$nzl and Vogt's idea of 2-fold monoidal category.
    The former pair required the natural transformation $\eta_{ABCD}$
    to be an isomorphism and showed that the resulting category is a braided monoidal category.
    As explained in \cite{Balt}, given such a category one
    obtains an equivalent braided monoidal category by ignoring one of the two
    operations, say $\otimes_2$, and defining the braiding for the
    remaining operation $\otimes_1$ to be the composite
    $$
    \diagram
    A\otimes_1 B\rrto^{\eta_{IABI}}
    && B\otimes_2 A\rrto^{\eta_{BIIA}^{-1}}
    && B\otimes_1 A.
    \enddiagram
    $$

    In \cite{Balt} it is shown that a 2-fold monoidal category with
    $\otimes_1 = \otimes_2 =\otimes$, $\eta$ an isomorphism and
    $$\eta_{AIBC} = \eta_{ABIC} = 1_{A \otimes B \otimes C}$$
    is a braided monoidal category with the braiding $c_{BC} = \eta_{IBCI}.$

 Also note that for ${\cal V}$ braided the interchange given by
 $\eta_{ABCD} = 1_A\otimes c_{BC} \otimes 1_D$ gives a 2-fold monoidal category where
 $\otimes_1 = \otimes_2 =\otimes.$ This interchange has the underlying braid $\sigma_2 \in B_4.$
In this setting we ask whether, given a braiding, there
  are alternate 2-fold monoidal structures on ${\cal V}$, with $\otimes_1 = \otimes_2 =\otimes$.
This is the same question as asking whether there are other interchanging unital braids besides
$b_{0^+} = b_{(1)} = \sigma_2$ and its inverse. To be precise, given a braided category  $({\cal
V},\otimes,\alpha,c,I)$ (with strict units, a strong associator $\alpha$, and braiding $c$), we ask
the central question: For which four-strand braids $b$ does the category ${\cal V}$ have in general
a coherent 2-fold monoidal structure, when that structure has $\otimes_1 = \otimes_2 = \otimes$ as
functors, has $\alpha^1 = \alpha^2 = \alpha$ as natural transformations, has strict unit $I$ for
both identical tensor products, and has $b$ as the underlying braid of $\eta$?

\begin{lemma}\label{y}
Given an arbitrary braided category ${\cal V}$, let the 2-fold structure of ${\cal V}$ be given by
$\otimes_1 = \otimes_2 =\otimes$. Then a four-strand braid $b$ is interchanging and unital if and
only if any interchange $\eta$ with underlying braid $b$ obeys the axioms of a 2-fold monoidal
category.
\end{lemma}
\begin{proof}
$Lb = Rb$ implies the internal associativity axiom of a 2-fold monoidal category,
 and $L'b = R'b$ implies the external associativity axiom, as we have foreshadowed with the naming of these braid
 equalities. This is seen by the coherence theorem for braided categories.  The unit axioms for the interchange
are also implied by the unit conditions on the braid, described by the fact that deleting certain
pairs of strands yields the identity braid. The converse implication is found by letting ${\cal V}$
be the free braided category. Then the axioms of a 2-fold monoidal category become precisely the
desired braid equalities.
\end{proof}

Now we are almost ready to state and prove the main result. First there are a couple of geometric
observations to be made about the braids $b_{n^{\pm}}= (\sigma_2\sigma_1\sigma_3\sigma_2)^{\pm
n}\sigma_2^{\pm 1}(\sigma_1\sigma_3)^{\mp n}.$
 Recall that we refer to the strands of a
braid by their initial positions. A sub-braid will refer to the braid resulting from the deletion
of  a subset of the strands of a braid.

\begin{lemma}\label{uno}
  If $n$ is odd  then deleting the outer two strands in the braid $b_{n^{\pm}}$  leaves the two strand sub-braid
$\sigma_1^{\pm n},$ while deleting the inner two strands gives the sub-braid $\sigma_1^{\pm(n+1)}.$
If $n$ is even then deleting the outer two strands in the braid $b_{n^{\pm}}$  leaves the two
strand sub-braid $\sigma_1^{\pm(n+1)},$ while deleting the inner two strands gives the sub-braid
$\sigma_1^{\pm n}.$
\end{lemma}
\begin{proof}
Consider the upper portion of the braid $b_{n^{\pm}}$ given by
$(\sigma_2\sigma_1\sigma_3\sigma_2)^{\pm n}.$ The outer two strands and the inner pair of strands
both are crossed $\pm n$ times. If $n$ is even then the upper portion is pure and so the next
generator $\sigma_2^{\pm 1}$ is applied to the inner two strands. If $n$ is odd then the upper
portion has the associated permutation which sends $\{1~ 2~3~4\}\to\{3~4~1~2\} $ and so the  next
 generator $\sigma_2^{\pm 1}$ is applied to the outer two strands. Note that the lower portion of the braid given
by $(\sigma_1\sigma_3)^{\mp n}$ contributes no further crossings to either the outer or inner
sub-braids.
\end{proof}
\begin{corollary}
The braids $b_{n^{\pm}}$ and $b_{m^{\pm}}$ are equivalent if and only if $m=n$ and the superscript
signs are the same.
\end{corollary}
\begin{proof}
For two braids to be equivalent it  is necessary that all their corresponding sub-braids be
equivalent. If $n,m$ are both odd (or both even) and the signs are the same then the implication is
clear by Lemma~\ref{uno}. Let $n,m$ be  both odd (or both even) with the signs not the same. Then
if we assume that $b_{n^{\pm}}$ and $b_{m^{\pm}}$ are equivalent then use of Lemma~\ref{uno} leads
to the absurd implication $1=-1.$
 Let $n$ be odd and $m$ be even with the superscript signs the same. Then if the sub-braids formed by the
outer strands are equal we have that $m=n+1.$ Then $m+1 = n+2 \ne n$ so the inner sub-braids are
not equal. Finally  let $n$ be odd and $m$ be even with the superscript signs not the same. If the
braids are equivalent then $m = -(n+1)$ but both $m$ and $n$ are required to be non-negative, so
this is a contradiction.
\end{proof}

Now the main result:

\begin{theorem}\label{main}
A braid $b \in B_4$ is interchanging and unital if and only if it is equivalent to one of the
braids $b_{n^{\pm}}.$
\end{theorem}

\begin{proof}
We will show:  [$b = b_{n^{\pm}}$] $\Longrightarrow$ [$b$ gives rise to a middle four interchange
for a monoidal structure on ${\cal V}$-Cat for arbitrary ${\cal V}$]$\Longrightarrow$  [$b$
interchanging and unital] $\Longrightarrow$ [$b = b_{n^{\pm}}$].

By Theorem \ref{newtimes} the middle four interchanges given by, suppressing the associators,
 $$\eta_{ABCD} =(c^{\mp n}_{CA}\otimes c^{\mp n}_{DB})\circ (1_C\otimes c^{\pm 1}_{DA}
\otimes 1_B)\circ(c^{\pm n}_{(A\otimes B)(C\otimes D)}),$$
   with underlying braid $b_{n^{\pm}},$ are indeed each a middle four interchange.

   Therefore by Lemma~\ref{underlies}
the braids $b_{n^{\pm}}$ obey $L'b_{n^{\pm}} = R'b_{n^{\pm}}$. Since the braids $b_{n^{\pm}}$ are
equal to their own 180 degree rotations, as mentioned in the introduction, this also implies that
$Lb_{n^{\pm}} = Rb_{n^{\pm}}.$ The internal unit conditions are fairly easy to verify by inspection
of the braids $b_{n^{\pm}}$; deleting the first two or the last two strands leaves the identity.
The external unit conditions are checked just as easily if we again do so using the 180 degree
rotations of $b_{n^{\pm}}.$

For the converse we assume that $b$ is interchanging  and unital and therefore by Lemma~\ref{y} it
underlies an interchange $\eta_{(b)}$ in a braided category ${\cal V}$ seen as a 2-fold monoidal
category with $\otimes_1 = \otimes_2 =\otimes.$  We focus on the two strand sub-braids of $b$
underlying $\eta_{(b)_{AIIB}}$ (the outer sub-braid) and $\eta_{(b)_{IABI}}$ (the inner sub-braid).
We will now show that a selection of these two underlying braids uniquely determines the braid $b.$

First assume that we have chosen a two-strand braid to underlie the inner sub-braid of $b.$
Consider the internal associativity
 axiom but with $U = W =  Z=I.$ Now the
top  horizontal arrow of the diagram has as its underlying braid the three-strand identity braid.
The left vertical side of the diagram has as the underlying braid formed by placing the underlying
braid of $\eta_{(b)_{IXYI}}$  (the inner sub-braid)  to the right of a single strand. The right
vertical side has the underlying braid of
 $\eta_{(b)_{I,V\otimes X,Y,I}}.$ This latter is just the choice we made for the inner two-strand sub-braid of $b$, with
the first strand doubled. The bottom horizontal arrow has the underlying braid of
$\eta_{(b)_{IVYX}}.$ This last
 three strand
sub-braid  of $b$ is thus determined by the assumption that the diagram commutes, the braided
coherence theorem, and
 the operation of taking the
inverse in the braid group $B_3.$
 Thus we have determined the three strand sub-braid of $b$ formed by
deleting the first strand.

Next we assume that we have chosen a braid to underlie $\eta_{(b)_{AIIB}}$. Then
 we again use the internal associativity diagram, this time with $V = X = Y = I $ to similarly
 determine the underlying braid
of $\eta_{(b)_{UIWZ}},$ i.e. the three strand sub-braid of $b$ formed by deleting the second
strand.

Finally we set $V = W =I$ in the internal associativity diagram. Now the top horizontal arrow has
the underlying braid formed by placing the outer two-strand sub-braid of $b$ to the left of the
two-strand identity braid. The left vertical side has the underlying braid formed by placing a
predetermined three strand sub-braid of $b$ (formed by deleting the first strand) to the right of a
single strand. The bottom horizontal arrow has the underlying braid formed by doubling the last
strand of a predetermined three-strand sub-braid of $b$ (formed by deleting the second strand).
Thus by
 the operation of taking the
inverse in the braid group $B_4$  we can determine the braid underlying the right vertical side.
This is precisely the braid $b$, underlying $\eta_{(b)_{UXYZ}}.$

Next we will limit the choices we can make for the underlying braids of $\eta_{(b)_{IABI}}$ and
$\eta_{(b)_{AIIB}}.$ We utilize Joyal and Street's result that, for any interchange $\eta$, a
braiding is given by:
$$  \diagram
    A\otimes_1 B\rrto^{\eta_{IABI}}
    && B\otimes_2 A\rrto^{\eta_{BIIA}^{-1}}
    && B\otimes_1 A.
    \enddiagram
    $$

Now our $\eta_{(b)_{IABI}}$ and $\eta_{(b)_{AIIB}}$ have underlying two-strand braids. Thus by
Lemma~\ref{dos} and braided coherence we have the equation
$$  \eta_{(b)_{BIIA}}^{-1} \circ
    \eta_{(b)_{IABI}}
    = c_{AB}^{\pm 1}~~~~ %
    $$
  or  $$  ~~~~\eta_{(b)_{IABI}}  = \eta_{(b)_{BIIA}} \circ c_{AB}^{\pm 1} .
    $$

Now in order for the permutation associated to $b$ to be $(2 ~ 3)$,
 $\eta_{(b)_{IABI}}$ must be an odd power of $c$ or $c^{-1}.$ Therefore our choice for the underlying braids of
$\eta_{(b)_{IABI}}$ and $\eta_{(b)_{AIIB}}$ is reduced respectively to a choice of an odd integer
$z$ and a choice of one of its neighboring integers $z\pm 1$. The choice of $z$ is the power of the
$c$, and thus the power of $\sigma_1$ for the inner sub-braid. The latter choice of
 $\pm 1$ is the choice of the exponent of $c$ in the above equation, and thus determines
 the power of $\sigma_1$ for the outer sub-braid.

Now by Lemma~\ref{uno} these possible choices for the underlying braids of $\eta_{(b)_{AIIB}}$ and
$\eta_{(b)_{IABI}}$ are all actually represented by one of the $b_{n^{\pm}}.$

Therefore if any braid $b$ is interchanging and unital then it is equivalent to one of the braids
$b_{n^{\pm}}.$

\end{proof}

The next item on the agenda is to investigate the equivalence of the various 2-fold monoidal
structures which can be constructed from a braiding, with differing underlying interchanging
braids.

\section{Equivalence of 2-fold Monoidal Categories}

By finding interchanges which are formed from a braiding we have actually defined a collection of
functors from the category of braided categories to the category of 2-fold monoidal categories. The
complete classification of interchanging unital braids is a well defined parameterization of this
family.
\begin{definition}
For $b$ an interchanging unital braid, the functor $F_b$ takes each braided category ${\cal V}$ to
itself, seen as a 2-fold monoidal category with interchange $\eta_{(b)}.$  A braided tensor functor
$f$ with $\phi:f(A)\otimes f(B) \to f(A\otimes B)$  is taken by $F_b$ to a 2-fold monoidal functor
$F_b(f)$ which has the same definition on objects and morphisms
  and for which $\lambda^1=\lambda^2 = \phi.$
\end{definition}
  \begin{theorem}\label{coolio}
Given an interchanging  unital braid $b$ the functor $F_b$ is naturally equivalent to either
$F_{b_{0^+}}$ or to $F_{b_{0^-}}$ but not to both.
  \end{theorem}
\begin{proof}
It is directly implied in \cite{JS} that given a 2-fold monoidal category ${\cal V}$ with
$\otimes_1 = \otimes_2$ and
 with strong interchange $\eta$ then that category is equivalent to
the 2-fold monoidal category ${\cal V'}$ with the same objects and morphisms but with interchange
given by
$$
\eta'_{ABCD} = 1_A\otimes(\eta^{-1}_{CIIB}\circ\eta_{IBCI})\otimes 1_D
$$

For ${\cal V}$ braided and in terms of an original interchange $\eta_{(b)}$ based on a braiding $c$
with $b$ interchanging  and unital,
 we have seen in the proof of Theorem~\ref{main}
that $\eta'_{ABCD} = 1_A\otimes c_{BC}^{\pm 1} \otimes 1_D.$ Thus ${\cal V'} =
F_{b_{0^{\pm}}}({\cal V}).$ The 2-fold monoidal functorial equivalence $U_{\cal
V}:F_{b_{0^{\pm}}}({\cal V})\to F_{b}({\cal V})$ is the identity on objects and morphisms.
Explicitly  $U_{\cal V}$ has $\lambda^2_{AB} = 1_{A\otimes B}$ and $\lambda^1_{AB} =
\eta_{(b)_{AIIB}}.$ This allows us to define in the target category:
 $$\eta_{U_{\cal V}(A)U_{\cal V}(B)U_{\cal V}(C)U_{\cal
V}(D)} = \eta_{(b)_{ABCD}}.$$ The required hexagonal interchange diagram commutes due to braided
coherence, using the braid equalities mentioned in Remark~\ref{foop}.

For $b$ such that $\eta^{-1}_{(b)_{CIIB}}\circ\eta_{(b)_{IBCI}} = c,$ i.e. $b \in \{ b_{n^{+}}~|~n
\text{ is even }\} \bigcup \{b_{n^{-}}~|~n \text{ is odd }\},$
 the family of functors
 $U_{\cal V}$ make
up a natural isomorphism $U:F_{b_{0^{+}}}\to F_{b}.$

For $b$ such that $\eta^{-1}_{(b)_{CIIB}}\circ\eta_{(b)_{IBCI}} = c^{-1},$ i.e. $b \in \{
b_{n^{+}}~ |~n \text{ is odd }\} \bigcup \{b_{n^{-}}~|~n \text{ is even }\},$
 the family of functors
 $U_{\cal V}$ make
up a natural isomorphism $U:F_{b_{0^{-}}}\to F_{b}.$

There is not in general a  natural isomorphism from $F_{b_{0^{-}}}$ to $F_{b_{0^{+}}}.$ If there
were then
 the hexagonal interchange diagram for 2-fold monoidal functors with $A = D =I$ would become the diagram of braided
 equivalence between  ${\cal V}$ with braiding $c$ and ${\cal V}$ with braiding $c^{-1}$.
 There is not in general a braided equivalence between ${\cal V}$ with braiding $c$ and ${\cal V}$ with
 braiding $c^{-1}$ since any $\lambda^2$ (in general based upon $c$) would have to
 satisfy $\lambda^2\circ c^{-1} = c \circ \lambda^2$ which is precluded by the braided coherence theorem
 and the fact that $B_2$ is abelian.

\end{proof}
Thus the interchanging  braids can be divided into two equivalence classes by the relation given by
$b\equiv b'$ if $F_b$ is equivalent to $F_b'.$ The two classes are canonically represented by the
braids $b_{0^{+}}$ and $b_{0^{-}}.$ It would be an interesting future study to consider the braid
groups $B_n$ for $n \ge 4$ modulo that equivalence relation on the first four strands.  With that
in mind we turn to examine some shortcuts to determining whether a given braid is interchanging and
unital.

 \section{Obstructions to being an interchange.}

  The general scheme is to find extra conditions on the interchange $\eta_{(b)}$ which
  together with the unit conditions and the associativity conditions will force
  the underlying braid $b$ to have easily checked characteristics. Then we can find
  families of unital braids in $B_4$ which cannot underlie an interchange, i.e. which are not
  equivalent to any braid $b_{n^{\pm}}.$
  \begin{theorem}
  Given an interchange candidate braid $b$ with the property that deleting either the 2nd or 3rd
  strand gives the identity braid on three strands, then $b$ is interchanging  if and only if
  $b = \sigma_2,$ the second generator of $B_4$, or its inverse.
  \end{theorem}

 \begin{proof} This follows the logic of \cite{Balt}. Letting $\eta = \eta_{(b)}$ be the interchange
  based on the braiding of ${\cal V}$
  with underlying braid $b$, note that deleting a strand in $b$ corresponds to replacing the
  respective object in the product $A\otimes B \otimes C \otimes D$ with the identity $I.$
   Now let $V=W=I$ in the internal associativity diagram to see that
  due to the hypotheses on $b$ we have that
  $\eta_{UXYZ}=1_U \otimes \eta_{IXYZ}.$ Then let $X=Y=I$ in the internal associativity diagram to
  see that $\eta_{UVWZ}=\eta_{UVWI} \otimes 1_Z.$ Together these two facts imply that
  $\eta_{ABCD} = 1_A \otimes \eta_{IBCI} \otimes 1_D.$ Then if we take $U=Z=W=0$ in the internal associativity
  law we get the first axiom of a braided category for $c'_{BC} = \eta_{IBCI},$ and letting $U=Z=X=0$ in the
  internal associativity diagram gives the other one. This then implies that either $c' = c$ or $c' = c^{-1},$
  since no other combinations of $c$ give a braiding. Therefore
  $\eta_{ABCD} = 1_A \otimes c^{\pm 1}_{BC} \otimes 1_D$  which has the underlying braid $\sigma_2^{\pm 1}.$
  The converse is also clear from this discussion, since all the implications can be reversed. Of
  course, we already have the converse since the braids $b_{0^{\pm}}$ are interchanging.
  \end{proof}

  This sort of obstruction can rule out candidate braids such as the braid $b_{(4)}$
  in the last section.
  It also rules out all but one element each of the left and
  right $\sigma^{\pm(2n-1)}_2$-cosets of the Brunnian braids in $B_4,$ where the
  Brunnian braids are those pure braids where any strand deletion gives the identity braid.
  Even more broadly
  this obstruction rules out braids such as:

  $$
  b = 
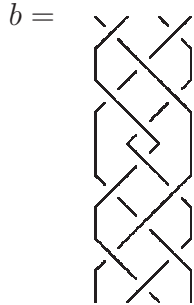
\xy 0;/r1pc/:
     ,{\xoverv}+(2,1),{\xoverv}+(-2,0)
     ,{\xcapv-@(0)}+(1,1),{\xunderv}+(2,1),{\xcapv-@(0)}+(-3,0)
     ,{\xunderv}+(2,1),{\xunderv}+(-2,0)
     ,{\xcapv-@(0)}+(1,1),{\xunderv}+(2,1),{\xcapv-@(0)}+(-3,0)
     ,{\xcapv-@(0)}+(1,1),{\xunderv}+(2,1),{\xcapv-@(0)}+(-3,0)
     ,{\xoverv}+(2,1),{\xoverv}+(-2,0)
     ,{\xcapv-@(0)}+(1,1),{\xoverv}+(2,1),{\xcapv-@(0)}+(-3,0)
     ,{\xunderv}+(2,1),{\xunderv}+(-2,0)
     ,{\xcapv-@(0)}+(1,1),{\xoverv}+(2,1),{\xcapv-@(0)}+(-3,0)
   \endxy
  $$

\begin{theorem}
Let $b$ be an interchange candidate braid with the property that deleting both the inner two
strands leaves
 the identity sub-braid
on the remaining two strands. Then  if $b$ is interchanging it follows that deleting either the
second or the third strand will result in the three strand identity sub-braid on either of the
remaining subsets of strands.
\end{theorem}
\begin{proof}
Let $\eta_{(b)}$ be the interchange based on the braiding, with underlying interchanging braid $b.$
We are given that $\eta_{(b)_{AIIB}} = 1_{A\otimes B}$ and must demonstrate that $\eta_{(b)_{AIBC}}
= \eta_{(b)_{ABIC}}
 = 1_{A\otimes B \otimes C}.$ The conclusion about the deletion of the second strand is shown by considering the internal
associativity diagram with $V = X = Y = I.$ The conclusion about the deletion of the third strand
is shown by considering the internal associativity diagram with $V = W = Y =I.$ An alternative
proof just uses the main results to check all the interchanging  braids which fit the hypothesis.
\end{proof}

This obstruction rules out all but one element each of the left and
  right $\sigma^{\pm(2n-1)}_2$-cosets of the $2$-trivial or $2$-decomposable braids in $B_4.$
These latter braids are a generalization of the Brunnian braids in which deletion of any 2 strands
results in a trivial braid.
$$$$

  Notice that the longer interchanging  braids $b_{n^{\pm}}$ for $n > 0$
give examples of interchanges that
  do not fit the conditions of the obstruction theorems so far. They also serve as examples of interchanges $\eta$ such that
  $\eta_{IBCI}$ is not a braiding. Recall however that they do give a braiding via
  $c'_{AB}= \eta^{-1}_{BIIA} \circ \eta_{IABI}$ as predicted by Joyal and Street. The latter condition
  also serves as a source of obstructions on its own. According to their theorem, any interchanging braid will
  have the property that dropping the outer two strands will give a two strand braid with one more or one less crossing of
  the same handedness
  than the two strand braid achieved by dropping the inner two strands. Indeed this condition rules out some
  of the same braids just mentioned, namely the Brunnian cosets of higher powers of $\sigma_2$ in $B_4.$

  The next sort of obstruction is found by slightly weakening the extra conditions.
  This will allow us to rule out a larger, different class of candidates, but they will be a little
  bit harder to recognize.

  \begin{theorem}
  Let $b$ be an interchange candidate braid  with the property that deleting either the first or the fourth
  strand results in a 3-strand braid that is just
  a power of the braid generator on what were the middle
   two strands: $\sigma_i^{\pm n} $;  $i= 2$ or $i=1$ respective of whether the first or fourth
  strand was deleted.
  Then $b$ is interchanging  implies that $n =1.$
  \end{theorem}

  \begin{proof}
  The strand deletion conditions on the underlying braid $b$ of $\eta$ are equivalent to assuming
  that $\eta_{IBCD} = \eta_{IBCI} \otimes 1_D$ and that $\eta_{ABCI} = 1_A \otimes \eta_{IBCI}$.
   Of course the  power of the generator $\sigma_i$ being $\pm 1$  is equivalent to saying that $\eta_{IBCI}$ is
  the braiding $c$ or its inverse.
  Hence we need only show that the assumptions imply that $\eta_{IBCI}$ is a braiding. This is seen
  immediately upon letting $U=Z=W=0$ in the internal associativity axiom to get the first axiom of a braiding and
  letting $U=Z=X=0$ to get the other one.
  \end{proof}

  This theorem can  directly rule out candidates
  which satisfy the Joyal and Street condition that  $c_{AB}= \eta^{-1}_{BIIA} \circ \eta_{IABI}$ and the first or
  last strand deletion condition given here, but which fail to give a single crossing braid upon that removal.
  The simplest example is this braid:
  $$
  
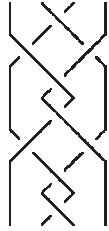
\xy 0;/r1pc/:
     ,{\xcapv-@(0)}+(1,1),{\xunderv}+(2,1),{\xcapv-@(0)}+(-3,0)
     ,{\xunderv}+(2,1),{\xoverv}+(-2,0)
     ,{\xcapv-@(0)}+(1,1),{\xunderv}+(2,1),{\xcapv-@(0)}+(-3,0)
     ,{\xcapv-@(0)}+(1,1),{\xunderv}+(2,1),{\xcapv-@(0)}+(-3,0)
     ,{\xoverv}+(2,1),{\xunderv}+(-2,0)
     ,{\xcapv-@(0)}+(1,1),{\xunderv}+(2,1),{\xcapv-@(0)}+(-3,0)
     ,{\xcapv-@(0)}+(1,1),{\xunderv}+(2,1),{\xcapv-@(0)}+(-3,0)
   \endxy
  $$
  It is also true that a candidate braid
  which yields a single crossing after deletion of the first and fourth strands, if interchanging,
  must then obey
  the condition that deleting the first or last strand frees the other of those two from any crossings.
  This can be most easily seen by use of the main result; we simply check all four examples of
  interchanging unital braids which have inner two strand sub-braids a single crossing. They are $b_{0^{\pm}}$
  and $b_{1^{\pm}}.$

\section{Obstructions to braiding in ${\cal V}$-Cat .}

Notice that in the case of symmetric ${\cal V}$ the axioms of enriched categories for ${\cal
A}\otimes {\cal B}$ and
  the existence of a coherent 2-natural associator follow from the coherence
  of symmetric categories and the
  enriched axioms for ${\cal A}$ and ${\cal B}.$
 It remains to consider just why it is
  that ${\cal V}$-Cat is braided if and only if ${\cal V}$ is symmetric, and that if so then ${\cal V}$-Cat
  is symmetric as well. This fact is stated in \cite{JS}. We choose to give a proof here which
  covers all possible interchanging braids explicitly, and all potential braidings on ${\cal V}$-Cat  based on any
  odd power of the braiding on ${\cal V},$ by appealing to information from the theory of knots and
  links. This is opposed to arguments based on the fact that a braiding transports over a tensor
  equivalence, and on Theorem~\ref{coolio}.

  Our choice allows us to demonstrate how low dimensional
  topology
  can inform category theory as well as vice versa.
  A braiding $c^{(1)}$ on ${\cal V}$-Cat is a 2-natural transformation so
  $c^{(1)}_{{\cal A}{\cal B}}$ is a ${\cal V}$-functor
  ${\cal A}\otimes {\cal B} \to {\cal B}\otimes {\cal A}$. On objects $c^{(1)}_{{\cal A}{\cal B}}((A,B)) = (B,A).$
  Now to be precise we define $c^{(1)}$ to be based upon $c$ to mean that
  $$c^{(1)}_{{\cal A}{\cal B}_{(A,B)(A',B')}}: ({\cal A}\otimes {\cal B})((A,B),(A',B')) \to ({\cal B}\otimes {\cal A})((B,A),(B',A'))$$
  is defined to be:   $$ c_{{\cal A}(A,A'){\cal B}(B,B')}:{\cal A}(A,A')\otimes {\cal B}(B,B') \to {\cal B}(B,B')\otimes {\cal A}(A,A')$$
  This potential braiding must be checked for ${\cal V}$-functoriality.
  Again the unit axioms are trivial and we consider the
  more interesting
  associativity of hom-object morphisms property. The following diagram must commute
  \begin{footnotesize}
  $$
  \xymatrix{
  ({\cal A}\otimes {\cal B})((A',B'),(A'',B''))\otimes({\cal A}\otimes {\cal B})((A,B),(A',B'))
  \ar[rr]^-{M}
  \ar[d]^{c^{(1)}\otimes c^{(1)}}
  &&({\cal A}\otimes {\cal B})((A,B),(A'',B''))
  \ar[d]^{c^{(1)}}
  \\
  ({\cal B}\otimes {\cal A})((B',A'),(B'',A''))\otimes ({\cal B}\otimes {\cal A})((B,A),(B',A'))
  \ar[rr]^-M
  &&({\cal B}\otimes {\cal A})((B,A),(B'',A''))
  }
  $$
  \end{footnotesize}
  Let $X = {\cal A}(A',A'')$, $Y = {\cal B}(B',B'')$, $Z = {\cal A}(A,A')$ and $W = {\cal B}(B,B')$
  Then expanding the above diagram using the composition defined as above (denoting various composites
  of $\alpha$ by unlabeled arrows) we have
  $$
  \xymatrix{
  &(X\otimes Y) \otimes (Z\otimes W)
  \ar[dr]
  \ar[dl]|{c_{XY} \otimes  c_{ZW}}\\
  (Y\otimes X) \otimes (W\otimes Z)
  \ar[d]
  &&X\otimes ((Y \otimes Z)\otimes W)
  \ar[d]^{1 \otimes (c_{YZ} \otimes 1)}\\
  Y\otimes ((X \otimes W)\otimes Z)
  \ar[d]^{1 \otimes (c_{XW} \otimes 1)}
  &&X\otimes ((Z \otimes Y)\otimes W)
  \ar[d]
  \\
  Y\otimes ((W \otimes X)\otimes Z)
  \ar[d]
  &&(X\otimes Z) \otimes (Y\otimes W)
  \ar[dll]|{c_{(X\otimes Z)(Y\otimes W)}}
  \ar[d]^{M_{AA'A''}\otimes M_{BB'B''}}\\
  (Y\otimes W) \otimes (X\otimes Z)
  \ar[dr]|{M_{BB'B''}\otimes M_{AA'A''}}
  &&{\cal A}(A,A'')\otimes {\cal B}(B,B'')
  \ar[dl]^{c}\\
  &{\cal B}(B,B'')\otimes {\cal A}(A,A'')
  }
  $$
  The bottom quadrilateral commutes by naturality of $c$. The top region must then commute for the diagram to commute, but
  the left and right legs have the following underlying braids
  $$
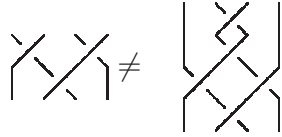
\xy (0,0) *{
  \xy 0;/r1pc/:
  ,{\xoverv}+(2,1),{\xoverv}+(-2,0)
  ,{\xcapv-@(0)}+(1,1),{\xoverv}+(2,1),{\xcapv-@(0)}+(-3,0)
  \endxy}\endxy
  \xy (0,0) *{~\ne~}\endxy
  \xy (0,0) *{
  \xy 0;/r1pc/:
  ,{\xcapv-@(0)}+(1,1),{\xoverv}+(2,1),{\xcapv-@(0)}+(-3,0)
  ,{\xcapv-@(0)}+(1,1),{\xoverv}+(2,1),{\xcapv-@(0)}+(-3,0)
  ,{\xoverv}+(2,1),{\xoverv}+(-2,0)
  ,{\xcapv-@(0)}+(1,1),{\xoverv}+(2,1),{\xcapv-@(0)}+(-3,0)
  \endxy}\endxy
  $$

  Thus as noted in \cite{JS}
  neither braid $b_{(1)}$ nor its inverse can in general give a monoidal structure
  with a braiding based on the original braiding. In fact, it is easy to show more.
  \begin{theorem}
  Let ${\cal V}$ be a braided category with braiding $c$. Let the tensor product on ${\cal V}$-Cat
  be given by the canonical choices for the objects, hom-objects, unit morphisms, unit enriched
  category, and associator, and let $b$ be the underlying braid of the composition morphisms for the
  tensor product of enriched categories. Then in general there will not be a
    braiding in ${\cal V}$-Cat based upon the braiding $c$ in
  ${\cal V}$. Moreover, this failure will also be the case for attempts to produce a braiding in
${\cal V}$-Cat based upon any (odd) power $c^{2n+1}.$
  \end{theorem}

  \begin{proof} Notice that in the above braid
  inequality each side of the
  inequality consists of the braid which underlies the definition of the composition morphism,
  in this case $b_{(1)}$, and
  an additional braid which underlies the segment of the preceding diagram that corresponds to a composite of $c^{(1)}.$
  In terms of braid generators the left side of the braid inequality begins with $\sigma_1\sigma_3$
  corresponding to $c_{XY} \otimes  c_{ZW}$ and the right side of the braid inequality
  ends with $\sigma_2\sigma_1\sigma_3\sigma_2$ corresponding to $c_{(X\otimes Z)(Y\otimes W)}.$ Since the same braid
  $b$ must
  end the left side as begins the right side, then for the diagram to commute we require $b\sigma_1\sigma_3 = \sigma_2\sigma_1\sigma_3\sigma_2b.$
  This implies $\sigma_1\sigma_3 = b^{-1}\sigma_2\sigma_1\sigma_3\sigma_2b$, or that the braids
  $\sigma_1\sigma_3$ and $\sigma_2\sigma_1\sigma_3\sigma_2$ are conjugate in $B_4$. Conjugate braids have precisely the same
  link as their closures, but the closure of $\sigma_1\sigma_3$ is an unlinked pair of circles whereas
  the closure of $\sigma_2\sigma_1\sigma_3\sigma_2$ is the Hopf link.

 $$ 
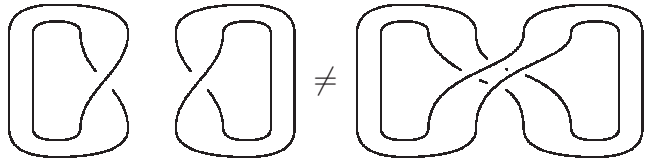
\xy 0;/r.15pc/:
  (-25,-10)*{}="tc3"; (-30,-10)*{}="tc4";
  (25,-10)*{}="tc2"; (30,-10)*{}="tc1";
  (-20,-12.5)*{}="td4"; (-20,-16)*{}="td3";
  (20,-12.5)*{}="td2"; (20,-16)*{}="td1";
(-25,10)*{}="c3"; (-30,10)*{}="c4";
  (25,10)*{}="c2"; (30,10)*{}="c1";
  (-20,12.5)*{}="d4"; (-20,16)*{}="d3";
  (20,12.5)*{}="d2"; (20,16)*{}="d1";
  (15,10)*{}="b1"; (5,10)*{}="b2";
  (-5,10)*{}="b3"; (-15,10)*{}="b4";
  (15,-10)*{}="T1"; (5,-10)*{}="T2";
  (-5,-10)*{}="T3"; (-15,-10)*{}="T4";
 "d1"; "c1" **\crv{(30,16) };
 "d1"; "b2" **\crv{(5,16) };
 "d2"; "c2" **\crv{(25,12.5) };
 "d2"; "b1" **\crv{(15,12.5) };
"d3"; "c4" **\crv{(-30,16)};
 "d3"; "b3" **\crv{(-5,16)};
 "d4"; "c3" **\crv{(-25,12.5) };
 "d4"; "b4" **\crv{(-15,12.5) };
  "c3"; "tc3" **\crv{"c3"};
  "c2"; "tc2" **\crv{"c2"};
  "c1"; "tc1" **\crv{"c1"};
  "c4"; "tc4" **\crv{"c4"};
  "td1"; "tc1" **\crv{(30,-16) };
 "td1"; "T2" **\crv{(5,-16) };
 "td2"; "tc2" **\crv{(25,-12.5) };
 "td2"; "T1" **\crv{(15,-12.5) };
"td3"; "tc4" **\crv{(-30,-16)};
 "td3"; "T3" **\crv{(-5,-16)};
 "td4"; "tc3" **\crv{(-25,-12.5) };
 "td4"; "T4" **\crv{(-15,-12.5) };
  "b1"; "T2" **\crv{(15,2) & (5,-2)}
   \POS?(.5)*{\hole}="2x";
  "b3"; "T4" **\crv{(-5,2) & (-15,-2)}
  \POS?(.5)*{\hole}="2y";
  "b2"; "2x" **\crv{(5,5) };
  "b4"; "2y" **\crv{(-15,5) };
  "2x"; "T1" **\crv{(15,-5)};
  "2y"; "T3" **\crv{(-5,-5)};
  \endxy
~~\ne~~
 \xy 0;/r.15pc/:
  (-25,-10)*{}="tc3"; (-30,-10)*{}="tc4";
  (25,-10)*{}="tc2"; (30,-10)*{}="tc1";
  (-20,-12.5)*{}="td4"; (-20,-16)*{}="td3";
  (20,-12.5)*{}="td2"; (20,-16)*{}="td1";
(-25,10)*{}="c3"; (-30,10)*{}="c4";
  (25,10)*{}="c2"; (30,10)*{}="c1";
  (-20,12.5)*{}="d4"; (-20,16)*{}="d3";
  (20,12.5)*{}="d2"; (20,16)*{}="d1";
  (15,10)*{}="b1"; (5,10)*{}="b2";
  (-5,10)*{}="b3"; (-15,10)*{}="b4";
  (15,-10)*{}="T1"; (5,-10)*{}="T2";
  (-5,-10)*{}="T3"; (-15,-10)*{}="T4";
 "d1"; "c1" **\crv{(30,16) };
 "d1"; "b2" **\crv{(5,16) };
 "d2"; "c2" **\crv{(25,12.5) };
 "d2"; "b1" **\crv{(15,12.5) };
"d3"; "c4" **\crv{(-30,16)};
 "d3"; "b3" **\crv{(-5,16)};
 "d4"; "c3" **\crv{(-25,12.5) };
 "d4"; "b4" **\crv{(-15,12.5) };
  "c3"; "tc3" **\crv{"c3"};
  "c2"; "tc2" **\crv{"c2"};
  "c1"; "tc1" **\crv{"c1"};
  "c4"; "tc4" **\crv{"c4"};
  "td1"; "tc1" **\crv{(30,-16) };
 "td1"; "T2" **\crv{(5,-16) };
 "td2"; "tc2" **\crv{(25,-12.5) };
 "td2"; "T1" **\crv{(15,-12.5) };
"td3"; "tc4" **\crv{(-30,-16)};
 "td3"; "T3" **\crv{(-5,-16)};
 "td4"; "tc3" **\crv{(-25,-12.5) };
 "td4"; "T4" **\crv{(-15,-12.5) };
  "b2"; "T4" **\crv{(5,2) & (-15,2.5)}
  \POS?(.35)*{\hole}="2x" \POS?(.54)*{\hole}="2y";
  "b4"; "2y" **\crv{(-15,6) };
  "b3"; "2x" **\crv{(-5,5) };
  "b1"; "T3" **\crv{(15,3.5) & (-5,3.5)}
  \POS?(.55)*{\hole}="4x" \POS?(.7)*{\hole}="4y";
    "4y";"2y" **\crv{(-2,-1)};
    "2x";"4x" **\crv{(1.75,1.75)};
  "4y";"T2" **\crv{(5,-3)};
    "T1";"4x" **\crv{(15,-2)};
\endxy $$

If we instead let
$$c^{(1)}_{{\cal A}{\cal B}_{(A,B)(A',B')}} = c^{2n+1}_{{\cal A}(A,A'){\cal B}(B,B')}$$
then the requirement becomes that the braids $(\sigma_1\sigma_3)^{2n+1}$ and
$(\sigma_2\sigma_1\sigma_3\sigma_2)^{2n+1}$ are conjugate in $B_4.$ Both braids have as closure a
link of two components--two copies of the $(2n+1,2)$-torus knot. However the first closure is two
unlinked copies of the knot while in  the second closure the two (cabled) copies are
 linked with linking number
$2n+1.$ Thus the braids cannot be conjugate, and so the braids underlying the legs of the
functoriality
 diagram will not be equal for any choice of middle four interchange.
\end{proof}

  \begin{corollary}
  ~It is also interesting to note that the braid inequality above is the 180 degree rotation of the one which
  implies that in general $({\cal A}\otimes {\cal B})^{op} \ne {\cal A}^{op}\otimes {\cal B}^{op}.$
  Thus the proof also implies that the latter inequality holds in general for a
  tensor products of enriched categories with any braid $b$
  underlying their composition morphisms, as well as any power of $op$ as the exponent.
  \end{corollary}

  \begin{remark}
    ~It is quickly seen that if $c$ is a symmetry then in the second half of the braid inequality the upper portion
    of the braid consists of $c_{YZ}$ and $c_{ZY} = c_{YZ}^{-1}$ so in fact equality holds. In that case then the
  derived braiding $c^{(1)}$ is a symmetry simply due to the definition.
  \end{remark}

\section{Implications for operads.}

 So far
  herein we have completely characterized
   families of interchanges based on a braiding which can define either a 2-fold monoidal
  structure on a category or a monoidal structure on a 2-category.  Another common use of a braiding
  is  to define a
  monoidal structure on a category of collections, as in the theory of operads. Operads in a 2-fold
  monoidal category are defined as monoids in a certain category of collections in \cite{forcey2}.
  Here we repeat the basic ideas and the expanded definition in terms of commuting diagrams.
  The two principle components of an operad are a collection,
historically a sequence, of objects in a monoidal category and a family of composition maps.
Operads are often described as parameterizations  of $n$-ary operations.  Peter May's original
definition of operad in a symmetric (or braided) monoidal category \cite{May} has a composition
$\gamma$ that takes the tensor product of the $n\th$ object ($n$-ary operation) and $n$ others (of
various arity) to a resultant that sums the arities of those others.  The $n\th$ object or $n$-ary
operation is often pictured as a tree with $n$ leaves, and the composition appears like this:

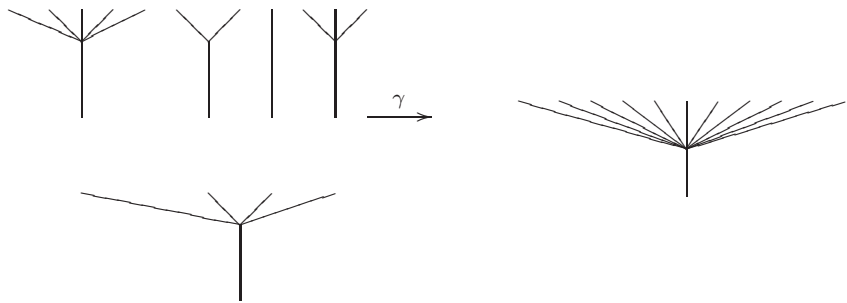
\begin{tabular}{ll}
$$
\xymatrix@W=0pc @H=2.2pc @R=1pc @C=1pc{
*=0{}\ar@{-}[drr] &*=0{}\ar@{-}[dr] &*=0{}\ar@{-}[d] &*=0{}\ar@{-}[dl] &*=0{}\ar@{-}[dll] &*=0{}\ar@{-}[dr] &*=0{} &*=0{}\ar@{-}[dl] &*=0{}\ar@{-}[d] &*=0{}\ar@{-}[dr] &*=0{}\ar@{-}[d] &*=0{}\ar@{-}[dl]\\
*=0{} &*=0{} &*=0{}\ar@{-}[d] &*=0{} &*=0{} &*=0{} &*=0{}\ar@{-}[d] &*=0{} &*=0{}\ar@{-}[d] &*=0{}  &*=0{}\ar@{-}[d] &*=0{} &*=0{}&*=0{}&*=0{}&*=0{}&*=0{}&*=0{}\\&*=0{}&*=0{}&*=0{}&*=0{}&*=0{}&*=0{}&*=0{}&*=0{}&*=0{}&*=0{}&*=0{}\ar[rr]^{\gamma}&*=0{}&*=0{}&*=0{}&*=0{}\\
*=0{} &*=0{} &*=0{}\ar@{-}[drrrrr] &*=0{} &*=0{} &*=0{}&*=0{}\ar@{-}[dr] &*=0{} &*=0{}\ar@{-}[dl] &*=0{} &*=0{}\ar@{-}[dlll] &*=0{} &*=0{}&*=0{}&*=0{}&*=0{}&*=0{}&*=0{}&*=0{}&*=0{}&*=0{}\\
*=0{} &*=0{} &*=0{} &*=0{}      &*=0{} &*=0{} &*=0{}&*=0{}\ar@{-}[d] &*=0{} &*=0{}\\ &*=0{} &*=0{} &*=0{} &*=0{} &*=0{} &*=0{} &*=0{} &*=0{}&*=0{}&*=0{}&*=0{}
}
$$ & \hspace{-6pc}
$$
\xymatrix@W=0pc @H=2.2pc @R=1.5pc @C=1pc{\\
*=0{}\ar@{-}[drrrrr]&*=0{}\ar@{-}[drrrr]&*=0{}\ar@{-}[drrr]&*=0{}\ar@{-}[drr]&*=0{}\ar@{-}[dr]&*=0{}\ar@{-}[d]&*=0{}\ar@{-}[dl]&*=0{}\ar@{-}[dll]&*=0{}\ar@{-}[dlll]&*=0{}\ar@{-}[dllll]&*=0{}\ar@{-}[dlllll]\\
*=0{}&*=0{}&*=0{}&*=0{}&*=0{}&*=0{}\ar@{-}[d]\\
*=0{}*=0{}*=0{}*=0{}*=0{}*=0{}*=0{}&*=0{}&*=0{}&*=0{}&*=0{}&*=0{}&*=0{}&*=0{}&*=0{}&*=0{}\\
}
$$
\end{tabular}

By requiring this composition to be associative we mean that it obeys this sort of pictured
commuting diagram:

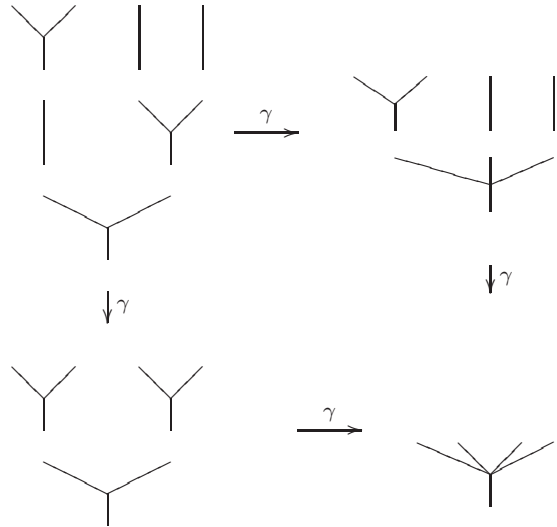
\begin{center}
\begin{tabular}{ll}
$$
\xymatrix@W=0pc @H=2.2pc @R=1pc @C=1pc{
*=0{}\ar@{-}[dr] &*=0{} &*=0{}\ar@{-}[dl] &*=0{} &*=0{}\ar@{-}[d] &*=0{} &*=0{}\ar@{-}[d] \\
*=0{} &*=0{}\ar@{-}[d] &*=0{} &*=0{} &*=0{}\ar@{-}[d] &*=0{} &*=0{}\ar@{-}[d] \\
*=0{}&*=0{}&*=0{}&*=0{}&*=0{}&*=0{}&*=0{}&*=0{}&*=0{}\\
*=0{} &*=0{}\ar@{-}[d] &*=0{} &*=0{} &*=0{}\ar@{-}[dr] &*=0{} &*=0{}\ar@{-}[dl] \\
*=0{} &*=0{}\ar@{-}[d] &*=0{} &*=0{} &*=0{} &*=0{}\ar@{-}[d] &*=0{}&*=0{}\ar[rr]^{\gamma}&*=0{}&*=0{}&*=0{}&*=0{} \\
*=0{}&*=0{}&*=0{}&*=0{}&*=0{}&*=0{}&*=0{}&*=0{}&*=0{}\\
*=0{} &*=0{}\ar@{-}[drr]  &*=0{} &*=0{} &*=0{} &*=0{}\ar@{-}[dll]\\
*=0{} &*=0{} &*=0{} &*=0{}\ar@{-}[d] \\
*=0{} &*=0{} &*=0{} &*=0{} \\
*=0{} &*=0{} &*=0{} &*=0{}\ar[d]^{\gamma}&*=0{}\\
*=0{} &*=0{} &*=0{} &*=0{} &*=0{}
}
$$ & \hspace{-4pc}
$$
\xymatrix@W=0pc @H=2.2pc @R=.85pc @C=1pc{\\
*=0{}\ar@{-}[dr] &*=0{} &*=0{}\ar@{-}[dl] &*=0{} &*=0{}\ar@{-}[d] &*=0{} &*=0{}\ar@{-}[d] \\
*=0{} &*=0{}\ar@{-}[d] &*=0{} &*=0{} &*=0{}\ar@{-}[d] &*=0{} &*=0{}\ar@{-}[d] \\
*=0{}&*=0{}&*=0{}&*=0{}&*=0{}&*=0{}&*=0{}&*=0{}&*=0{}\\
*=0{} &*=0{}\ar@{-}[drrr] &*=0{} &*=0{} &*=0{}\ar@{-}[d] &*=0{} &*=0{}\ar@{-}[dll] \\
*=0{} &*=0{} &*=0{} &*=0{} &*=0{}\ar@{-}[d] &*=0{} &*=0{} \\
*=0{}&*=0{}&*=0{}&*=0{}&*=0{}&*=0{}&*=0{}\\
*=0{} \\
*=0{} &*=0{} &*=0{} &*=0{} &*=0{}\ar[d]^{\gamma}&*=0{}\\
*=0{} &*=0{} &*=0{} &*=0{} &*=0{} &*=0{}
}
$$ \\
$$
\xymatrix@W=0pc @H=2.2pc @R=1pc @C=1pc{
*=0{}\ar@{-}[dr] &*=0{} &*=0{}\ar@{-}[dl] &*=0{} &*=0{}\ar@{-}[dr] &*=0{} &*=0{}\ar@{-}[dl] \\
*=0{} &*=0{}\ar@{-}[d] &*=0{} &*=0{} &*=0{} &*=0{}\ar@{-}[d] &*=0{} \\
*=0{}&*=0{}&*=0{}&*=0{}&*=0{}&*=0{}&*=0{}&*=0{}&*=0{}&*=0{}\ar[rr]^{\gamma}&*=0{}&*=0{}&*=0{}&*=0{} \\
*=0{} &*=0{}\ar@{-}[drr]  &*=0{} &*=0{} &*=0{} &*=0{}\ar@{-}[dll]\\
*=0{} &*=0{} &*=0{} &*=0{}\ar@{-}[d] \\
*=0{} &*=0{} &*=0{} &*=0{} \\
}
$$
& \hspace{-2pc}
$$
\xymatrix@W=0pc @H=2.2pc @R=1pc @C=1pc{\\
*=0{}\ar@{-}[drr] &*=0{}\ar@{-}[dr] &*=0{} &*=0{}\ar@{-}[dl] &*=0{}\ar@{-}[dll]\\
*=0{} &*=0{} &*=0{}\ar@{-}[d] &*=0{} \\
*=0{} &*=0{} &*=0{} &*=0{}
}
$$
\end{tabular}
\end{center}

In the above pictures the tensor products are shown just by juxtaposition, but now we would like to
think about the products more explicitly. If the monoidal category is not strict, then there is
actually required another leg of the associativity diagram, where the tensoring is reconfigured  so
that the composition can operate in an alternate order.  Here is how that rearranging looks in a
symmetric (braided) category, where the shuffling is accomplished by use of the symmetry
(braiding):

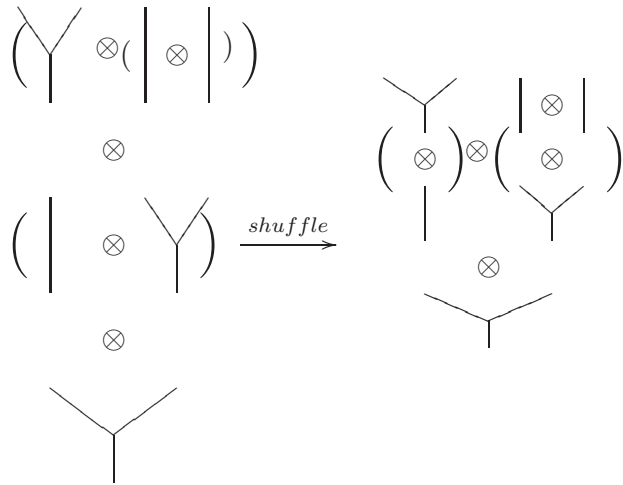
\begin{center}
\begin{tabular}{ll}
$$
\xymatrix@W=0pc @H=2.2pc @R=1.5pc @C=1pc{
*=0{}\ar@{-}[dr] &*=0{} &*=0{}\ar@{-}[dl] &*=0{} &*=0{}\ar@{-}[d] &*=0{} &*=0{}\ar@{-}[d] \\
*=0{\txt{\huge(}} &*=0{}\ar@{-}[d] &*=0{} &*=0{\otimes\txt{(}} &*=0{}\ar@{-}[d] &*=0{\otimes} &*=0{}\ar@{-}[d] &*=0{\left.\right)\txt{\huge)}}\\
*=0{}&*=0{}&*=0{}&*=0{}&*=0{}&*=0{}&*=0{}&*=0{}&*=0{}\\
*=0{}&*=0{}&*=0{}&*=0{\otimes}&*=0{}&*=0{}&*=0{}&*=0{}&*=0{}\\
*=0{} &*=0{}\ar@{-}[d] &*=0{} &*=0{} &*=0{}\ar@{-}[dr] &*=0{} &*=0{}\ar@{-}[dl] \\
*=0{\txt{\huge(}} &*=0{}\ar@{-}[d] &*=0{} &*=0{\otimes} &*=0{} &*=0{}\ar@{-}[d] &*=0{\txt{\huge)}}&*=0{}\ar[rrr]^{shuffle}&*=0{}&*=0{}&*=0{}&*=0{} \\
*=0{}&*=0{}&*=0{}&*=0{}&*=0{}&*=0{}&*=0{}&*=0{}&*=0{}\\
*=0{}&*=0{}&*=0{}&*=0{\otimes}&*=0{}&*=0{}&*=0{}&*=0{}&*=0{}\\
*=0{} &*=0{}\ar@{-}[drr]  &*=0{} &*=0{} &*=0{} &*=0{}\ar@{-}[dll]\\
*=0{} &*=0{} &*=0{} &*=0{}\ar@{-}[d] \\
*=0{} &*=0{} &*=0{} &*=0{} \\
}
$$ & \hspace{-1.25pc}
$$
\xymatrix@W=0pc @H=2.2pc @R=.85pc @C=1pc{\\
*=0{}\ar@{-}[dr] &*=0{} &*=0{}\ar@{-}[dl] &*=0{} &*=0{}\ar@{-}[d] &*=0{} &*=0{}\ar@{-}[d] \\
*=0{} &*=0{}\ar@{-}[d] &*=0{} &*=0{} &*=0{}\ar@{-}[d] &*=0{\otimes} &*=0{}\ar@{-}[d] &*=0{}\\
*=0{}&*=0{}&*=0{}&*=0{}&*=0{}&*=0{}&*=0{}&*=0{}&*=0{}\\
*=0{\txt{\Huge(}}&*=0{\otimes}&*=0{\txt{\Huge)}}&*=0{\otimes\txt{\Huge(}}&*=0{}&*=0{\otimes}&*=0{}&*=0{\txt{\Huge)}}&*=0{}\\
*=0{} &*=0{}\ar@{-}[d] &*=0{} &*=0{} &*=0{}\ar@{-}[dr] &*=0{} &*=0{}\ar@{-}[dl] \\
*=0{} &*=0{}\ar@{-}[d] &*=0{} &*=0{} &*=0{} &*=0{}\ar@{-}[d] &*=0{}&*=0{} \\
*=0{}&*=0{}&*=0{}&*=0{}&*=0{}&*=0{}&*=0{}&*=0{}&*=0{}\\
*=0{}&*=0{}&*=0{}&*=0{\otimes}&*=0{}&*=0{}&*=0{}&*=0{}&*=0{}\\
*=0{} &*=0{}\ar@{-}[drr]  &*=0{} &*=0{} &*=0{} &*=0{}\ar@{-}[dll]\\
*=0{} &*=0{} &*=0{} &*=0{}\ar@{-}[d] \\
*=0{} &*=0{} &*=0{} &*=0{} \\
}
$$
\end{tabular}
\end{center}

We now foreshadow our definition of operads in an iterated monoidal category with the same picture
as above but using two tensor products, $\otimes_1$ and $\otimes_2.$ It becomes clear that the true
nature of the shuffle is in fact that of an interchange transformation.

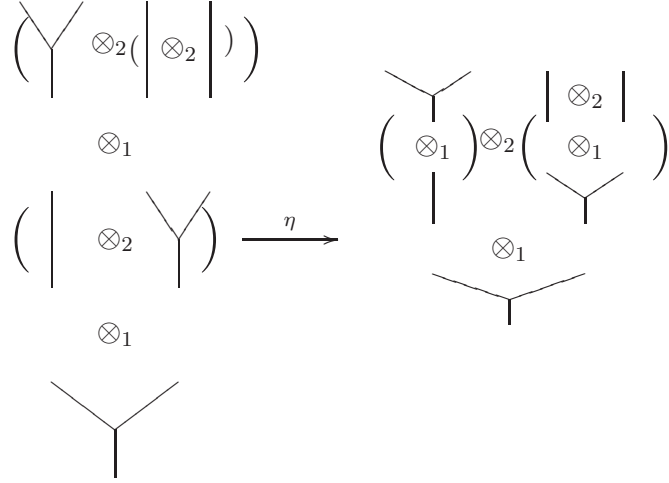
\begin{center}
\begin{tabular}{ll}
$$
\xymatrix@W=0pc @H=2.2pc @R=1.5pc @C=1pc{
*=0{}\ar@{-}[dr] &*=0{} &*=0{}\ar@{-}[dl] &*=0{} &*=0{}\ar@{-}[d] &*=0{} &*=0{}\ar@{-}[d] \\
*=0{\txt{\huge(}} &*=0{}\ar@{-}[d] &*=0{} &*=0{\otimes_2\txt{(}} &*=0{}\ar@{-}[d] &*=0{\otimes_2} &*=0{}\ar@{-}[d] &*=0{\left.\right)\txt{\huge)}}\\
*=0{}&*=0{}&*=0{}&*=0{}&*=0{}&*=0{}&*=0{}&*=0{}&*=0{}\\
*=0{}&*=0{}&*=0{}&*=0{\otimes_1}&*=0{}&*=0{}&*=0{}&*=0{}&*=0{}\\
*=0{} &*=0{}\ar@{-}[d] &*=0{} &*=0{} &*=0{}\ar@{-}[dr] &*=0{} &*=0{}\ar@{-}[dl] \\
*=0{\txt{\huge(}} &*=0{}\ar@{-}[d] &*=0{} &*=0{\otimes_2} &*=0{} &*=0{}\ar@{-}[d] &*=0{\txt{\huge)}}&*=0{}\ar[rrr]^{\eta}&*=0{}&*=0{}&*=0{}&*=0{} \\
*=0{}&*=0{}&*=0{}&*=0{}&*=0{}&*=0{}&*=0{}&*=0{}&*=0{}\\
*=0{}&*=0{}&*=0{}&*=0{\otimes_1}&*=0{}&*=0{}&*=0{}&*=0{}&*=0{}\\
*=0{} &*=0{}\ar@{-}[drr]  &*=0{} &*=0{} &*=0{} &*=0{}\ar@{-}[dll]\\
*=0{} &*=0{} &*=0{} &*=0{}\ar@{-}[d] \\
*=0{} &*=0{} &*=0{} &*=0{} \\
}
$$ & \hspace{-1.25pc}
$$
\xymatrix@W=0pc @H=2.2pc @R=.8pc @C=1.2pc{\\
*=0{}\ar@{-}[dr] &*=0{} &*=0{}\ar@{-}[dl] &*=0{} &*=0{}\ar@{-}[d] &*=0{} &*=0{}\ar@{-}[d] \\
*=0{} &*=0{}\ar@{-}[d] &*=0{} &*=0{} &*=0{}\ar@{-}[d] &*=0{\otimes_2} &*=0{}\ar@{-}[d] &*=0{}\\
*=0{}&*=0{}&*=0{}&*=0{}&*=0{}&*=0{}&*=0{}&*=0{}&*=0{}\\
*=0{\txt{\Huge(}}&*=0{\otimes_1}&*=0{\txt{\Huge)}}&*=0{\otimes_2\txt{\Huge(}}&*=0{}&*=0{\otimes_1}&*=0{}&*=0{\txt{\Huge)}}&*=0{}\\
*=0{} &*=0{}\ar@{-}[d] &*=0{} &*=0{} &*=0{}\ar@{-}[dr] &*=0{} &*=0{}\ar@{-}[dl] \\
*=0{} &*=0{}\ar@{-}[d] &*=0{} &*=0{} &*=0{} &*=0{}\ar@{-}[d] &*=0{}&*=0{} \\
*=0{}&*=0{}&*=0{}&*=0{}&*=0{}&*=0{}&*=0{}&*=0{}&*=0{}\\
*=0{}&*=0{}&*=0{}&*=0{\otimes_1}&*=0{}&*=0{}&*=0{}&*=0{}&*=0{}\\
*=0{} &*=0{}\ar@{-}[drr]  &*=0{} &*=0{} &*=0{} &*=0{}\ar@{-}[dll]\\
*=0{} &*=0{} &*=0{} &*=0{}\ar@{-}[d] \\
*=0{} &*=0{} &*=0{} &*=0{} \\
}
$$
\end{tabular}
\end{center}

To see this just focus on the actual domain and range of $\eta$ which are the upper two levels of
trees in the pictures, with the tensor product $\left({\mathbf|}\otimes_2{\mathbf|}\right)$
considered as a single object.

Now we are ready to give the technical definitions. We begin with the definition of 2-fold operad
in an $n$-fold monoidal category, as in the above picture, and then mention how it generalizes the
case of operad in a braided category.

\begin{definition}
Let ${\cal V}$ be a strict $2$-fold monoidal category. A 2-fold operad ${\cal C}$ in ${\cal V}$
consists of objects ${\cal C}(j)$, $j\ge 0$, a unit map ${\cal J}:I\to {\cal C}(1)$, and
composition maps in ${\cal V}$
\[
\gamma^{12}:{\cal C}(k) \otimes_1 ({\cal C}(j_1) \otimes_2 \dots \otimes_2 {\cal C}(j_k))\to {\cal
C}(j)
\]
for $k\ge 1$, $j_s\ge0$ for $s=1\dots k$ and $\smash{\sum\limits_{n=1}^k j_n = j}$. The composition
maps obey the following axioms:
\begin{enumerate}
\item Associativity: The following diagram is required to commute
for all  $k\ge 1$, $j_s\ge 0$ and $i_t\ge 0$, and where $\sum\limits_{s=1}^k j_s = j$ and
$\sum\limits_{t=1}^j i_t = i.$ Let $g_s= \sum\limits_{u=1}^s j_u$ and let
$h_s=\sum\limits_{u=1+g_{s-1}}^{g_s} i_u$.  The $\eta$ labeling the leftmost arrow actually stands
for a variety of equivalent maps which factor into instances of the interchange.
\[
\xymatrix{ {\cal C}(k)\otimes_1\left(\bigotimes\limits_{s=1}^k{}_2 {\cal C}(j_s)\right)\otimes_1
\left(\bigotimes\limits_{t=1}^j{}_2 {\cal C}(i_t)\right) \ar[rr]^>>>>>>>>>>>>{\gamma^{12} \otimes_1
\text{id}} \ar[dd]_{\text{id} \otimes_1 \eta} &&{\cal C}(j)\otimes_1
\left(\bigotimes\limits_{t=1}^j{}_2 {\cal C}(i_t)\right) \ar[d]^{\gamma^{12}}
\\
&&{\cal C}(i)
\\
{\cal C}(k)\otimes_1 \left(\bigotimes\limits_{s=1}^k{}_2 \left({\cal C}(j_s)\otimes_1
\left(\bigotimes\limits_{u=1}^{j_s}{}_2 {\cal C}(i_{u+g_{s-1}})\right)\right)\right)
\ar[rr]_>>>>>>>>>>{\text{id} \otimes_1(\otimes^k_2\gamma^{12})} &&{\cal
C}(k)\otimes_1\left(\bigotimes\limits_{s=1}^k{}_2 {\cal C}(h_s)\right) \ar[u]_{\gamma^{12}} }
\]

\item Respect of units is required just as in the symmetric case.
The following unit diagrams commute.
\[
\xymatrix{ {\cal C}(k)\otimes_1 (\otimes_2^k I) \ar[d]_{1\otimes_1(\otimes_2^k {\cal J})}
\ar@{=}[r]^{}
&{\cal C}(k)\\
{\cal C}(k)\otimes_1(\otimes_2^k {\cal C}(1)) \ar[ur]_{\gamma^{12}} } \xymatrix{ I\otimes_1 {\cal
C}(k) \ar[d]_{{\cal J}\otimes_1 1} \ar@{=}[r]^{}
&{\cal C}(k)\\
{\cal C}(1)\otimes_1 {\cal C}(k) \ar[ur]^{\gamma^{12}} }
\]
\end{enumerate}
\end{definition}
Now the problem of  describing the various
  sorts of operads in a braided monoidal category becomes more clear, as a special case.
  Here again we let
  $\otimes = \otimes_1= \otimes_2.$  The family of 2-fold structures
  based on interchanging  braids gives rise to a family of monoidal structures on the category of collections, and
  thus to a family of operad structures.

  In the operad picture the underlying braid of an operad structure
  only becomes important when we inspect the various ways of
  composing a product with 4 levels of trees in the heuristic diagram, such as
  ${\cal C}(2)\otimes({\cal C}(1)\otimes{\cal C}(1))\otimes ({\cal C}(1)\otimes{\cal C}(1))\otimes({\cal C}(1)\otimes{\cal C}(1)).$
  For this composition to be well defined we require the internal associativity of the interchange
  that is used to rearrange
  the terms. When we consider composing a product with 3 levels of trees in the heuristic diagram,
  but with a base term  ${\cal C}(n)$ with  $n \ge 3$, such as:
  ${\cal C}(3)\otimes({\cal C}(1)\otimes{\cal C}(1)\otimes{\cal C}(1))\otimes ({\cal C}(1)\otimes{\cal C}(1)\otimes{\cal C}(1)),$
  then we see that the external associativity of $\eta$ is also required.

  Thus the same theorems proven above for interchanging and non-interchanging families
  of braids  apply here as well, in deciding whether a certain braid based shuffling of the terms
  in an operad product
  is allowable. The point is that not all shuffles  using a braiding make sense, and
  the viewpoint of
  the 2-fold monoidal structure is precisely what is needed to see which shuffles
  do make sense. By seeing various shuffles as being
  interchanges on a fourfold product rather than braidings on a simple binary product,
  we are able to describe  an infinite family of distinct
  compositions of the braiding each leading to well defined operad structure. The underlying braids are
  precisely those we denoted $b_{n^{\pm}}$.  In summary, structures based on a braiding  are
  at worst ill-defined, at best defined up to equivalence,
  unless a 2-fold monoidal structure is chosen. Often in the literature
  the default is understood to be the simplest such structure where
  $\eta_{ABCD} = 1_A\otimes c_{BC} \otimes 1_D,$ but to be careful this choice should be made
  explicit.
We have directly addressed operads and tensor products of enriched categories. The results herein
should also be applied to ${\cal V}$-Act, the category of categories with an action of a monoidal
category as described in \cite{pad}, as well as to ${\cal V}$-Mod, the bicategory of enriched
categories and modules as described in \cite{street2}.


\newpage
      \begin{thebibliography}{99}
      \bibitem[Baez, 1997]{Baez}{John C. Baez, Higher-Dimensional Algebra II. 2-Hilbert Spaces,
      Advances in Mathematics 127, 125-189 (1997)}
    \bibitem[Baez and Dolan, 1998]{Baez1}{J. C. Baez and J. Dolan, Categorification, in ``Higher Category Theory'',
        eds. E. Getzler and M. Kapranov, Contemp. Math. 230 , American Mathematical Society, 1-36, (1998).}
    \bibitem[Balteanu et.al, 2003]{Balt}{C. Balteanu, Z. Fiedorowicz, R. Schw${\rm \ddot a}$nzl, R. Vogt,
        Iterated Monoidal Categories,
        Adv. Math. 176 (2003), 277-349.}
   \bibitem[Borceux, 1994]{Borc}{F. Borceux, Handbook of Categorical Algebra 1:
        Basic Category Theory, Cambridge University Press, 1994.}
   \bibitem[Day, 1970]{Day}{B.J. Day, On closed categories of functors, Lecture Notes in
        Math 137 (Springer, 1970) 1-38}
        \bibitem[Day and Street, 1997]
        {street2}{Brian Day and Ross Street, Monoidal Bicategories and Hopf Algebroids,
        Advances in Mathematics 129, 99-157 (1997)}
   \bibitem[Eilenberg and Kelly, 1965]{EK1}{S. Eilenberg and G. M. Kelly, Closed Categories,
        Proc. Conf. on Categorical Algebra,
        Springer-Verlag (1965), 421-562. }
   \bibitem[Forcey, 2004]{forcey1}{S. Forcey, Enrichment Over Iterated Monoidal Categories,
        Algebraic and Geometric Topology 4 (2004), 95-119.}
   \bibitem[Forcey, 2007]{forcey2}{S. Forcey ,J.Siehler, and E.Seth Sowers, Operads in iterated monoidal categories,
        \emph{to appear in} Journal of Homotopy and Related Structures 2 (2007), }
   \bibitem[Garside, 1969]{Gar}{F. Garside, The Braid Group and Other Groups, Quart. J. Math Oxford 20
(1969), 235-254.}
   \bibitem
       [Joyal and Street, 1993]
       {JS}{A. Joyal and R. Street, Braided tensor categories, Advances in Math. 102(1993), 20-78.}
    \bibitem[Kelly, 1982]{Kelly}{G. M. Kelly, Basic Concepts of Enriched Category Theory, London Math.
        Society Lecture Note Series 64, 1982.}
    \bibitem[Mac Lane, 1998]{MacLane}{S. Mac Lane, Categories for the
        Working Mathematician 2nd. edition, Grad. Texts in Math. 5, 1998.}
\bibitem[Mac Lane, 1965]
       {Mac}
        {S. MacLane, Categorical algebra, Bull. A. M. S. 71(1965), 40-106.}
   \bibitem[May, 1972]
       {May}{J. P. May, The geometry of iterated loop spaces, Lecture Notes in
        Mathematics, Vol. 271, Springer, 1972}
        \bibitem[McCrudden,2000]
        {pad}{Paddy McCrudden, Balanced Coalgebroids,
        Theory and Applications of Categories, Vol. 7, No. 6, 2000, pp. 71–147.}
\bibitem[Shum, 1994]
     {Shum}
       {M.C. Shum, Tortile tensor categories, J.of Pure and Apl. Alg. 193(194),no.1, 57-110.}
   \end{thebibliography}
\end{document}